# Tessellated Isotropic Elastic Lattice Spring Model for Quasi-Brittle Fracture

D. M. LI [a, b, c *], Meng-Cheng HE [a]

a. *School of Civil Engineering and Architecture, Wuhan University of Technology, Wuhan, 430070, China*

b. *Embodied Physical Intelligence Center in Civil Engineering, Wuhan University of Technology, Wuhan, 430070, China*

c. *Key Laboratory of Low-Altitude Technology and Smart Urban Renewal, Department of Housing and Urban-Rural Development of Hubei Province, Wuhan University of Technology, Wuhan, 430070, China*

**Abstract:** Quasi-brittle fracture is a prevalent failure mechanism in a wide range of engineering materials and structures, including concrete, rock, ceramics, composites, and masonry. Its numerical simulation faces a persistent trade-off among accuracy, efficiency, and implementation simplicity. The classical Lattice Spring Model (LSM) captures crack initiation and propagation through discrete bond breakage without remeshing, crack tracking, or enrichment, but its element construction remains empirical and is typically limited to a few seamlessly tessellable shapes with fixed Poisson's ratios. This paper proposes a tessellated Isotropic Elastic Lattice Spring Model (IELSM), which discretizes the continuum into polygonal tessellation elements with axial springs and an additional volumetric constraint, achieving isotropic elasticity on arbitrary polygonal tessellations. Macroscopic isotropy is reduced to a set of governing equations whose solvability provides a theoretical criterion for element admissibility. This criterion proves the admissibility of conventional tessellable elements and extends it to arbitrary regular *N*-gons and concave elements. Exploiting boundary interpolation compatibility with standard finite elements, IELSM is assembled through direct node sharing, without interface elements, overlapping regions, or kinematic constraints. IELSM is further coupled with an isotropic damage model for fracture simulation. A pure bending test and four fracture benchmarks show that the coupling preserves displacement accuracy, yields crack paths and load-displacement curves that agree well with experiments and outperform those obtained by standard FEM, and remains insensitive to mesh refinement. This compatibility also allows IELSM to be restricted to the damage-prone region, while the remainder of the domain is modeled by conventional finite elements through direct node sharing. For the considered benchmarks, this restriction reduces the number of nodes by 41.2% – 80.9% and the CPU time by 34.2% – 85.2% compared with full-domain IELSM. The framework advances IELSM element construction from empirical trial and error to theoretical determination and simplifies finite element coupling to node sharing, offering a balanced route for complex quasi-brittle fracture analysis.

**Keywords:** isotropic elastic lattice spring model; tessellation; coupling with finite elements; isotropic damage model; quasi-brittle fracture

[*]Corresponding to D. M. LI (dongmli2-c@my.cityu.edu.hk, domili@whut.edu.cn).

## 1. Introduction

Quasi-brittle fracture is the dominant failure mechanism in engineering materials such as concrete, rock, and masonry, and its numerical simulation has long been one of the central challenges in computational solid mechanics [1, 2]. Among the numerous numerical methods available, the lattice spring model (LSM) occupies a distinctive position [3]. In the classical LSM, a continuum is represented as a network of interconnected axial springs (trusses). In fracture simulations, cracks emerge naturally from the breakage of discrete bonds, without requiring remeshing, crack-tracking algorithms, or special enrichment [4–8]. However, since Hrennikoff's pioneering work [9], the construction of lattice geometries has long relied on empirical choices. This construction has mainly been limited to three shapes that admit seamless monohedral tessellation: equilateral triangles, squares, and regular hexagons [10–13]. Consequently, the equivalent continua of these lattices are restricted to fixed Poisson's ratios (1/3 under plane stress and 1/4 under plane strain). To remove this Poisson's ratio restriction, subsequent studies introduced beam elements [14], nonlocal interactions [15], angular springs [16], and shear springs [17], thereby enriching the LSM framework. Nevertheless, purely axial spring systems are still used today owing to their simplicity and compatibility. In particular, the recently proposed Isotropic Elasticity Lattice Spring Model (IELSM) introduces an additional bulk-modulus term into axial springs to remove the Poisson's ratio restriction [13]; this method has clear physical meaning and a concise formulation. However, among the many lattice models, there is still a lack of theoretical explanation for which element shapes are "admissible"; nor is there a clear criterion for whether admissible shapes other than those capable of monohedral tessellation exist. In mathematics, Grünbaum and Shephard [18] systematically established tessellation theory, which aims to study how the plane can be tessellated by polygons without overlaps or gaps and to provide a theoretical characterization of their classification and construction. Similarly, an element shape usable in the LSM must not only satisfy the geometric condition of seamless tessellation but also ensure that the spring system arranged within it satisfies an isotropic elastic constitutive law. In other words, the admissibility of an LSM element is subject to the dual constraints of geometric tessellation and isotropic elasticity. This theoretical gap is precisely the scientific question addressed in this paper.

To simulate material fracture using lattice models, researchers have proposed various failure criteria and fracture algorithms. The most traditional approach is to directly define failure criteria to identify critical elements to be removed from the model. Most of these criteria can be classified as critical bond force (stress) criteria and critical displacement (strain) criteria [19–22]. However, such methods have several inherent limitations: oversimplified failure criteria struggle to capture complex stress states [4]; the criterion parameters require calibration for specific loading conditions and are difficult to generalize [23]; and, more importantly, the simulation results depend on the mesh size [24]. To address these issues, researchers have proposed fracture-energy-based piecewise spring constitutive laws, such as bilinear and trilinear ones [24–27]. Such criteria do not directly remove springs that exceed a threshold; instead, they gradually degrade the spring stiffness through a softening constitutive law. A common feature of the above softening constitutive laws is that the area enclosed by the single-bond force-displacement curve and the horizontal axis is taken as the bond fracture work, which allows the constitutive laws of springs of different lengths to be calibrated using a unified fracture energy, thereby alleviating mesh-size dependence. However, their critical failure strain remains an empirical parameter that requires repeated calibration. In addition to piecewise softening constitutive laws, another approach is to combine lattice models with isotropic damage models, e.g., Yue et al. [28], Grassl et al. [29], Berton and Bolander [30], and Cusatis et al. [26, 27]. Recently, the Lattice Particle Method (LPM) [28, 29], combined with an isotropic damage model, has effectively addressed the above issues: it uses a single damage variable to characterize stiffness degradation and drive crack initiation. This isotropic damage model originates from the crack band theory (CBT) proposed by Bažant and Oh [35], in which the modeled crack band width is related to the finite element size; in LPM, the characteristic length is recalibrated for square packing. It should be noted that when CBT is used in FEM, it can eliminate mesh-size dependence, but the predicted crack paths still depend on mesh alignment [36]. The reason for this bias in standard FEM is that the stress/strain fields that drive damage evolution do not converge [32, 33]. In view of this, Cervera et al. [39] proposed a mixed finite element formulation to ensure the convergence of stress/strain fields in smeared crack methods and to obtain mesh-bias-free results in crack trajectories; however, the implementation of mixed finite element methods can be quite complex. As for LPM, its use of a particle method effectively avoids mesh dependence, but its theory is built on a regular square packing, so a uniform particle radius can only be used throughout the domain [28, 29]. This significantly increases unnecessary computational effort, and the computational cost grows sharply with mesh density.

Conventional LSMs are typically limited to uniform meshes, making it difficult to directly generate nodal density gradients. One possible solution is to establish an LSM theory for arbitrary element shapes so as to generate meshes with graded nodal densities; however, schemes introduced to overcome shape restrictions, such as rotational degrees of freedom and shear springs, often come at the cost of increased computational effort or reduced accuracy [13]. Another approach is to couple LSM with FEM: LSM is used only in regions of interest, while FEM continuum elements are used in non-critical regions to reduce nodal density, thereby maintaining accuracy while improving computational efficiency. Although this strategy is attractive, several models still face different difficulties in their coupling implementation, as illustrated by the following representative works. Schlangen [40] mixed beam-element lattices with plane-stress continuum elements. Since beam-element nodes possess rotational degrees of freedom whereas plane-stress elements have only translational degrees of freedom, and dependency relations can usually associate only translational degrees of freedom, interfacial bending moments cannot be transferred naturally. To avoid direct contact detection between distinct-LSM particles and elements in the Numerical Manifold Method, Zhao et al. [21] introduced an intermediate layer, the Particle based Manifold Method; this comes at the cost of accuracy in high-frequency responses and implementation simplicity. Zope [41] proposed a nonlocal method called the volumetric-compensated particle model (VCPM). When coupled with FEM, virtual VCPM particles must be attached to FEM nodes to transfer nonlocal effects. The calculation of interfacial forces by virtual particles and transition elements is an approximation; some spring energy may be double-counted or omitted by adjacent FEM elements, and its accuracy and convergence have not been fully verified. In addition, another widely studied direction is the quasicontinuum method [42–44]. Quasicontinuum methods retain the full discrete lattice, build continuum elements from representative nodes, derive their constitutive laws by summation, and interpolate the lattice displacements. However, studies on quasicontinuum methods [42–44] generally do not address the Poisson's ratio restriction of discrete lattices. Instead, these studies focus on crack initiation and propagation in discrete lattices, which is beyond the scope of this paper.

The objective of this study is to generalize the IELSM, based on the previous constitutive equivalence framework [13], into a general theory applicable to arbitrary polygonal elements, and to couple it with standard finite elements and a node-based isotropic damage model for quasi-brittle fracture simulation. To this end, this paper is organized into the following three aspects. First, a constitutive equivalence theory of IELSM applicable to arbitrary polygonal elements is established. This theory reduces the macroscopic isotropy requirement to five governing equations, namely four rotational invariance conditions and one elastic parameter matching condition. Given the element geometry and spring topology, the spring stiffnesses are the unknowns, and the solvability of the equation system determines whether the element shape is admissible. Depending on whether the number of unknowns is smaller than, equal to, or larger than the number of equations, the system is overdetermined, determined, or underdetermined, respectively. In the overdetermined case, a solution exists only when the element geometry satisfies specific symmetry conditions; this fundamentally explains why the classical LSM has long been limited to a few regular shapes such as equilateral triangles. In the determined case, the system can be solved directly. In the underdetermined case, additional constraints are required; in this paper, symmetry and equal-stiffness conditions are imposed to reduce the number of unknowns, yielding analytical solutions. On this basis, under the symmetry and equal-stiffness conditions, this paper proves that equilateral triangles, rectangles satisfying positive-definiteness constraints, regular hexagons, and arbitrary regular $N$-gons ($N \geq 5$) can all serve as admissible elements. An L-shaped concave polygonal element is used as an example to illustrate the applicability of the theory to concave elements. This extends the element types permitted by IELSM from the three elements capable of monohedral tessellation to arbitrary regular polygons (multiple element types can be combined to form a seamless tessellation).

The above governing equations do not depend on a specific element geometry, spring topology, or number of nodes, thus providing a unified framework for determining the admissibility of arbitrary polygonal shapes and enabling mesh refinement. It should be noted that volumetric compensation is not unique to IELSM: the LPM, also known as the VCPM [45], also includes a volumetric term. However, in LPM, volumetric compensation is incorporated into the interaction law between particle pairs and treated as a nonlocal term, making the model closely tied to the pairwise connectivity of a given packing. In contrast, the additional volumetric constraint in IELSM acts on the entire element through the local volumetric strain of the element and is treated independently in the theoretical derivation and numerical discretization. This not only simplifies the imposition of constraints but also confines both the springs and the volumetric constraint to the interior of the element, thereby ensuring locality. Second, an interface-free direct coupling between IELSM and standard finite elements is achieved. By exploiting the equivalence between the boundary displacement interpolation of IELSM elements and that of constant strain triangular (CST) and four-node isoparametric (Q4) elements, IELSM and standard finite elements can be assembled through direct node sharing, without the need for interface elements, overlapping regions, or kinematic constraints. Here, the damage-prone region is predefined according to the expected crack path and stress concentration. Consequently, IELSM is restricted to this region while the remaining region is covered by standard finite elements. This confines the fine-scale discrete description to the necessary region and avoids discretizing the entire domain with IELSM, thereby reducing the computational cost. Third, IELSM is combined with a node-based isotropic damage model for quasi-brittle fracture simulation. The constructed framework has been verified through a higher-order displacement field convergence test (pure bending test) and four fracture benchmarks, covering typical cases such as symmetric mode-I fracture with size effect, asymmetric mixed-mode fracture, complex geometry, and notched configurations. In summary, the core contributions of this paper are as follows: advancing the element construction of IELSM from empirical trial and error to theoretical determination, and simplifying its coupling with standard finite elements from interface engineering to node sharing. As a result, high-fidelity discretization can be applied only to the damage-prone region while material parameters retain their physical meaning. This framework provides an approach for complex quasi-brittle fracture analysis that balances accuracy, efficiency, and implementation simplicity.

The rest of this paper is organized as follows. Section 2 establishes the IELSM theory for arbitrary polygonal elements, derives the governing equations for isotropic elasticity, discusses the admissibility of representative element shapes that admit monohedral tessellation, and presents an extension to element shapes that do not admit monohedral tessellation. Section 3 presents the numerical implementation, including the derivation of element stiffness matrices, a mixed assembly strategy for IELSM and finite elements, stress/strain field computation for IELSM, and eigenvalue analysis. Section 4 combines IELSM with an isotropic damage model for quasi-brittle fracture simulation. Section 5 validates the proposed framework through five benchmark examples, namely a bending test, a symmetric three-point bending test, an asymmetric three-point bending test, an L-shaped panel, and a double-edge-notched specimen. Section 6 summarizes the paper. Appendices A, B, and C provide supplementary derivations and construction methods.

## 2. Theory of IELSM for Arbitrary Polygonal Tessellations

### 2.1 Governing Equations for Isotropic Elasticity of IELSM Elements

Building on the previously proposed IELSM for square representative volume elements (RVEs) [13], this work extends the formulation to arbitrary polygonal tessellation elements. The solution domain is regarded as a seamless tessellation composed of polygonal elements, with axial springs and an additional volumetric constraint arranged inside each element. If an element satisfies the macroscopic isotropy condition, it is referred to as an admissible tessellation element. Conventional LSMs mostly construct representative volume elements centered on a single node or particle; in contrast, this paper builds the model starting from an arbitrary polygonal tessellation element and does not require the element to be periodically repeatable. As shown in FIGURE 1, assuming that the strain field within the element is uniform and that the element contains $n$ axial springs and a uniform additional volumetric constraint, the strain energy of the $i$-th spring is:

$$U(\theta_i, L_i, k_i) = \frac{k_i \Delta l_i^2}{2} = \frac{k_i L_i^{\,2} \varepsilon_n^2(\theta_i)}{2}, \tag{1}$$

where $\theta_i$, $L_i$ and $k_i$ are the angle between the $i$-th axial spring and the $x$-axis (positive counterclockwise), the original length, and the stiffness, respectively; $\Delta l$ is the change in spring length, and $\varepsilon_n(\theta_i)$ is the axial strain of the spring, with:

$$\varepsilon_n(\theta_i) = \frac{\Delta l_i}{L_i} = \cos^2\theta_i \;\varepsilon_{xx} + \sin\theta_i \cos\theta_i (\varepsilon_{xy} + \varepsilon_{yx}) + \sin^2\theta_i \;\varepsilon_{yy}, \tag{2}$$

where $\varepsilon_{ij}$ represent global strain components. Although the strain energy of a single spring is independent of rigid-body motion (translation and rotation), once embedded in an element, its equivalent macroscopic strain energy density is closely related to the spring orientation. Substituting Eq. (2) into Eq. (1) and rearranging, the strain energy can be separated into terms independent of $\theta_i$:

$$\begin{aligned} U(\theta_i, L_i, k_i) &= \frac{k_i L_i^{\,2}}{2}\left(\cos^2\theta_i \;\varepsilon_{xx} + \sin\theta_i \cos\theta_i (\varepsilon_{xy} + \varepsilon_{yx}) + \sin^2\theta_i \;\varepsilon_{yy}\right)^2 \\ &= \frac{k_i L_i^{\,2}}{2}\left(\begin{array}{c} A^2 + \frac{1}{2}B^2 + \frac{1}{2}C^2 + \\ 2AB\cos(2\theta_i) + 2AC\sin(2\theta_i) + BC\sin(4\theta_i) + \frac{(B^2 - C^2)}{2}\cos(4\theta_i) \end{array}\right), \end{aligned} \tag{3}$$

where $A = (\varepsilon_{xx} + \varepsilon_{yy})/2$, $B = (\varepsilon_{xx} - \varepsilon_{yy})/2$, and $C = \gamma_{xy}/2$. Then, the total strain energy of all springs in the element is:

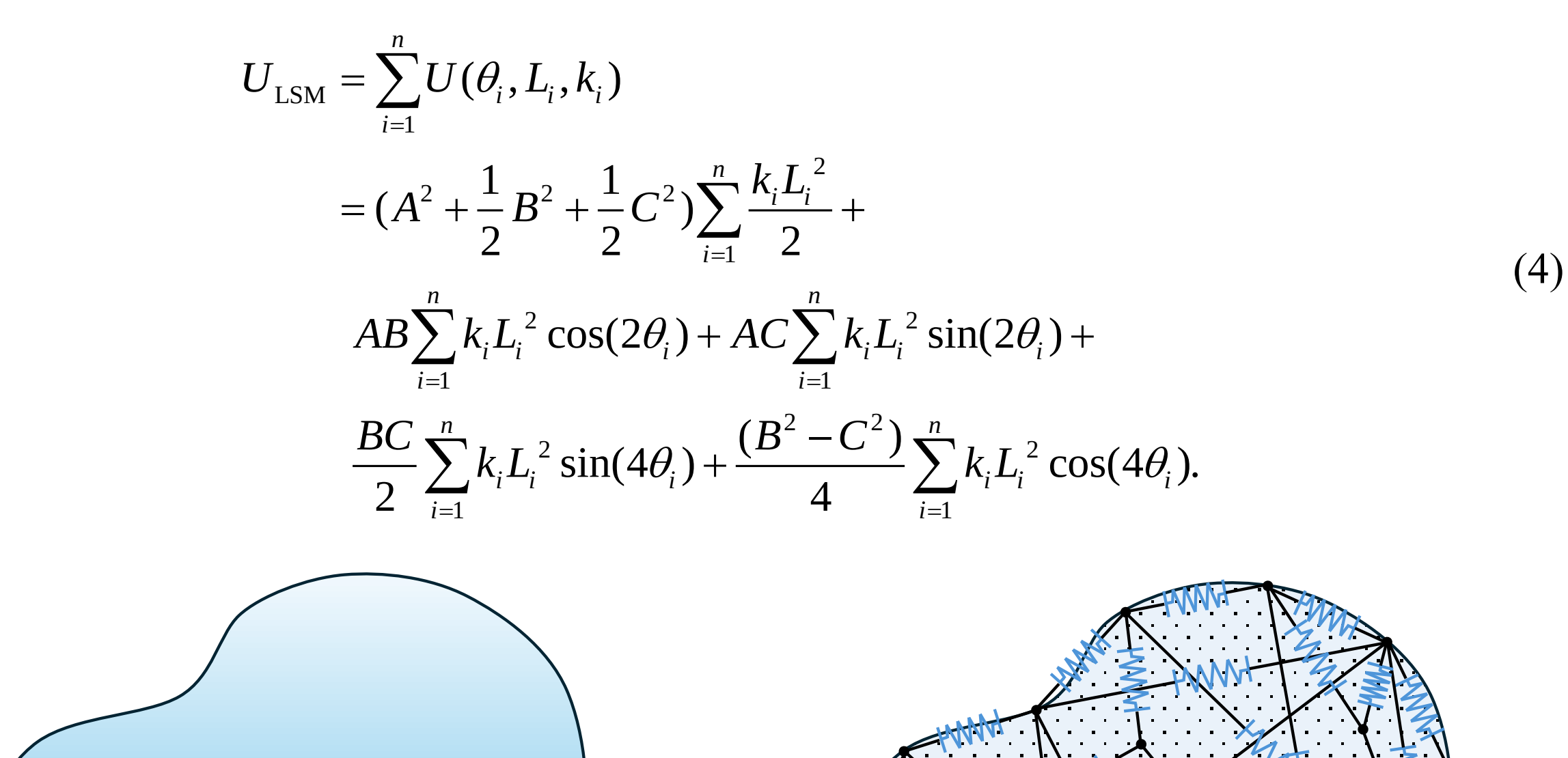

$$\begin{aligned} U_{\mathrm{LSM}} &= \sum_{i=1}^{n} U(\theta_i, L_i, k_i) \\ &= (A^2 + \frac{1}{2}B^2 + \frac{1}{2}C^2)\sum_{i=1}^{n} \frac{k_i L_i^{\,2}}{2} + \\ &\quad AB\sum_{i=1}^{n} k_i L_i^{\,2}\cos(2\theta_i) + AC\sum_{i=1}^{n} k_i L_i^{\,2}\sin(2\theta_i) + \\ &\quad \frac{BC}{2}\sum_{i=1}^{n} k_i L_i^{\,2}\sin(4\theta_i) + \frac{(B^2 - C^2)}{4}\sum_{i=1}^{n} k_i L_i^{\,2}\cos(4\theta_i). \end{aligned} \tag{4}$$

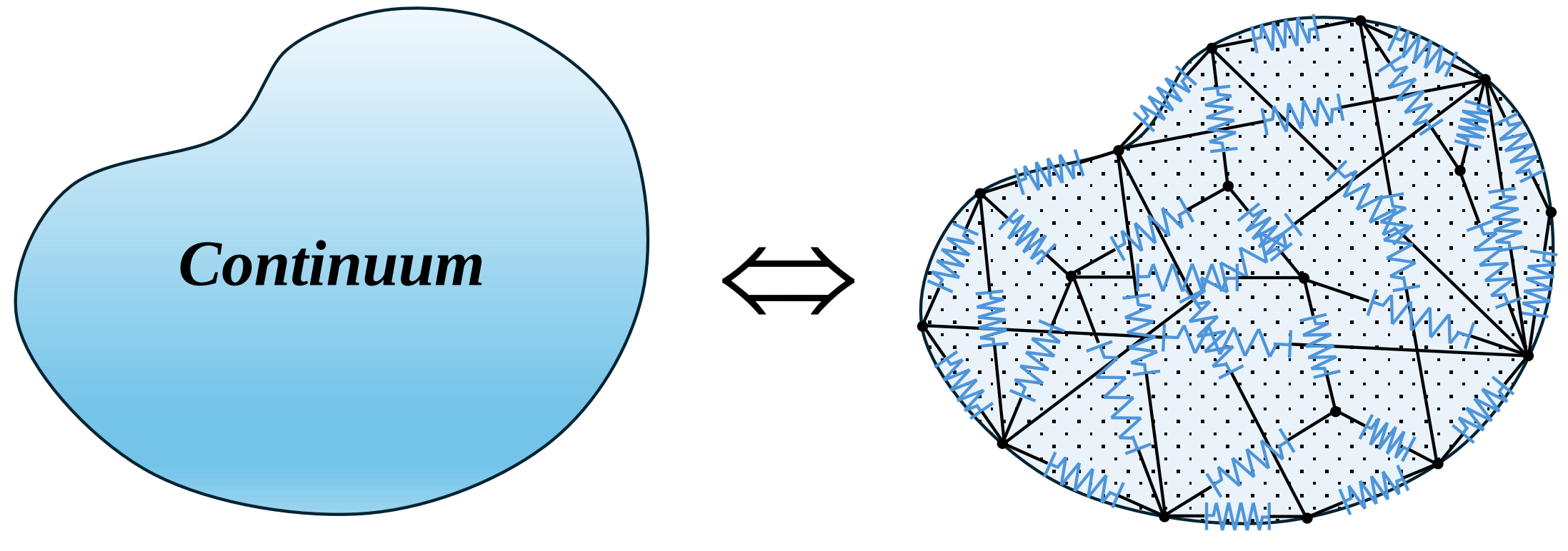

FIGURE 1 Conceptual illustration of the equivalence between IELSM and a continuum.

To make this collection of axial springs isotropic, a rigid-body rotation $\theta$ can be applied to the element. The directions of all springs then become $\theta_i + \theta$, and the strain state transforms accordingly in the rotated coordinate system. Isotropy requires that the total strain energy remain unchanged for any rotation angle $\theta$, i.e.,

$$\sum_{i=1}^{n} U(\theta_i, L_i, k_i) = \sum_{i=1}^{n} U(\theta_i + \theta,\ L_i,\ k_i), \tag{5}$$

expanding gives:

$$\begin{aligned} &AB\sum_{i=1}^{n} k_i L_i^2 \left(\cos(2\theta_i) - \cos(2\theta_i + 2\theta)\right) + AC\sum_{i=1}^{n} k_i L_i^2 \left(\sin(2\theta_i) - \sin(2\theta_i + 2\theta)\right) + \\ &\frac{BC}{2}\sum_{i=1}^{n} k_i L_i^2 \left(\sin(4\theta_i) - \sin(4\theta_i + 4\theta)\right) + \frac{(B^2 - C^2)}{4}\sum_{i=1}^{n} k_i L_i^2 \left(\cos(4\theta_i) - \cos(4\theta_i + 4\theta)\right) = 0. \end{aligned} \tag{6}$$

Since $\varepsilon_{xx}$, $\varepsilon_{yy}$, $\gamma_{xy}$ are independent variables and the transformation to $A$, $B$, $C$ is invertible, $A$, $B$, $C$ are independent variables. Hence the quadratic monomials $AB$, $AC$, $BC$, $B^2$, $C^2$ are linearly independent. Consequently, $AB$, $AC$, $BC$, $B^2$ - $C^2$ are also linearly independent. Therefore, from Eq. (6), we have:

$$\begin{cases} \sum_{i=1}^{n} k_i L_i^2 \left(\cos(2\theta_i) - \cos(2\theta_i + 2\theta)\right) = 0, \\ \sum_{i=1}^{n} k_i L_i^2 \left(\sin(2\theta_i) - \sin(2\theta_i + 2\theta)\right) = 0, \\ \sum_{i=1}^{n} k_i L_i^2 \left(\sin(4\theta_i) - \sin(4\theta_i + 4\theta)\right) = 0, \\ \sum_{i=1}^{n} k_i L_i^2 \left(\cos(4\theta_i) - \cos(4\theta_i + 4\theta)\right) = 0, \end{cases} \tag{7}$$

expanding gives:

$$\begin{cases} (1 - \cos(2\theta))\sum_{i=1}^{n} k_i L_i^2 \cos(2\theta_i) + \sin(2\theta)\sum_{i=1}^{n} k_i L_i^2 \sin(2\theta_i) = 0, \\ (1 - \cos(2\theta))\sum_{i=1}^{n} k_i L_i^2 \sin(2\theta_i) - \sin(2\theta)\sum_{i=1}^{n} k_i L_i^2 \cos(2\theta_i) = 0, \\ (1 - \cos(4\theta))\sum_{i=1}^{n} k_i L_i^2 \cos(4\theta_i) + \sin(4\theta)\sum_{i=1}^{n} k_i L_i^2 \sin(4\theta_i) = 0, \\ (1 - \cos(4\theta))\sum_{i=1}^{n} k_i L_i^2 \sin(4\theta_i) - \sin(4\theta)\sum_{i=1}^{n} k_i L_i^2 \cos(4\theta_i) = 0, \end{cases} \tag{8}$$

since $1 - \cos(2\theta)$ and $\sin(2\theta)$, and likewise $1 - \cos(4\theta)$ and $\sin(4\theta)$, are linearly independent, and Eq. (8) must hold for arbitrary $\theta$, we obtain:

$$
\begin{cases}
\sum_{i=1}^{n} k_i L_i^2 \cos(2\theta_i) = 0, \\
\sum_{i=1}^{n} k_i L_i^2 \sin(2\theta_i) = 0, \\
\sum_{i=1}^{n} k_i L_i^2 \cos(4\theta_i) = 0, \\
\sum_{i=1}^{n} k_i L_i^2 \sin(4\theta_i) = 0.
\end{cases} \tag{9}
$$

Equation system (9) is a set of necessary conditions for ensuring the macroscopic isotropy of the element. Substituting Eq. (9) back into Eq. (4), the total strain energy is simplified to an isotropic form independent of the spring orientation:

$$
\begin{aligned}
U_{\text{LSM}} &= (A^2 + \frac{1}{2}B^2 + \frac{1}{2}C^2)\sum_{i=1}^{n}\frac{k_i L_i^2}{2} \\
&= \left(\frac{3\varepsilon_{xx}^2 + 3\varepsilon_{yy}^2 + 2\varepsilon_{xx}\varepsilon_{yy} + \gamma_{xy}^2}{16}\right)\sum_{i=1}^{n} k_i L_i^2 \\
&= \frac{S}{2}K_0(\varepsilon_{xx} + \varepsilon_{yy})^2 + G\left(\varepsilon_{xx}^2 + \varepsilon_{yy}^2 - \frac{1}{3}(\varepsilon_{xx} + \varepsilon_{yy})^2 + \frac{1}{2}\gamma_{xy}^2\right)S,
\end{aligned} \tag{10}
$$

where $S$ is the area of the element, $G$ is the shear modulus, with $G = \sum_{i=1}^{n} k_i L_i^2 \Big/ 8S$, and $K_0$ is the bulk modulus under plane strain, with $K_0 = 5\sum_{i=1}^{n} k_i L_i^2 \Big/ 24S$. Under plane stress, $K_0$ is replaced by $K_0' = 9K_0G/(3K_0 + 4G)$, while $G$ remains unchanged.

Equation (10) is the total strain energy form of the classical LSM, indicating that its bulk modulus is a fixed value. Therefore, IELSM aims to directly adjust its bulk modulus, thereby modifying the strain energy. By imposing an additional volumetric constraint on the basis of the classical LSM total strain energy, the total strain energy of the IELSM RVE is obtained:

$$
\begin{aligned}
U_{\text{IELSM}} &= U_{\text{LSM}} + U_{\text{volume}} \\
&= \frac{S}{2}(K_0 + k_v)(\varepsilon_{xx} + \varepsilon_{yy})^2 + G\left(\varepsilon_{xx}^2 + \varepsilon_{yy}^2 - \frac{1}{3}(\varepsilon_{xx} + \varepsilon_{yy})^2 + 2\varepsilon_{xy}^2\right)S,
\end{aligned} \tag{11}
$$

where

$$
U_{\text{volume}} = \frac{1}{2}k_v\varepsilon_v^2 S =
\begin{cases}
\frac{1}{2}k_v(\varepsilon_{xx} + \varepsilon_{yy})^2 S & \text{, (plane strain)} \\
\frac{1}{2}k_v'\left(\frac{1-2\nu}{1-\nu}\right)^2(\varepsilon_{xx} + \varepsilon_{yy})^2 S & \text{, (plane stress)}
\end{cases} \tag{12}
$$

where $\varepsilon_v$ denotes volumetric strain. Under small-deformation assumptions: $\varepsilon_v = \varepsilon_{xx} + \varepsilon_{yy} + \varepsilon_{zz}$, $\varepsilon_{zz} = 0$ for plane strain conditions, thus $\varepsilon_v = \varepsilon_{xx} + \varepsilon_{yy}$; and $\varepsilon_{zz} = -\nu(\varepsilon_{xx} + \varepsilon_{yy})/(1 - \nu)$ for plane stress conditions, thus $\varepsilon_v = (1 - 2\nu)(\varepsilon_{xx} + \varepsilon_{yy})/(1 - \nu)$. To unify these formulations, we define $k_v = k_v'(1 - 2\nu)^2/(1 - \nu)^2$ as the normalized additional bulk modulus. The bulk modulus of the IELSM is given by $K = K_0 + k_v$. Therefore, the additional bulk modulus $k_v$ of the present model can be either positive or negative in a physical sense, and its purpose is to modify the bulk modulus of the LSM, thereby achieving an adjustable Poisson's ratio.

Dividing Eq. (11) by the element area S gives the strain energy density of the element. Then, taking its second-order partial derivatives with respect to the strain components gives the elastic coefficient matrix C:

$$\mathbf{C} = \begin{bmatrix} k_v + \frac{3}{8S}\sum_{i=1}^{n} k_i L_i^2 & k_v + \frac{1}{8S}\sum_{i=1}^{n} k_i L_i^2 & 0 \\ k_v + \frac{1}{8S}\sum_{i=1}^{n} k_i L_i^2 & k_v + \frac{3}{8S}\sum_{i=1}^{n} k_i L_i^2 & 0 \\ 0 & 0 & \frac{1}{8S}\sum_{i=1}^{n} k_i L_i^2 \end{bmatrix}. \tag{13}$$

Matching **C** with the linear elastic constitutive tensor $\mathbf{C}_0$ yields:

$$\sum_{i=1}^{n} k_i L_i^2 = \frac{4SE}{1+\nu}, \tag{14}$$

and

$$k_v = \begin{cases} \dfrac{E(4\nu - 1)}{2(1-2\nu)(1+\nu)} & \text{, (plane strain)} \\ \dfrac{E(3\nu - 1)}{2(1-\nu^2)} & \text{. (plane stress)} \end{cases} \tag{15}$$

From Eqs. (14) and (15), one can also obtain:

$$\begin{cases} E = \dfrac{\psi(24k_v S + 5\psi)}{16(4k_v S + \psi)S} \\ \nu = \dfrac{8k_v S + \psi}{4(4k_v S + \psi)} \end{cases} \text{, (plane strain)} \qquad \begin{cases} E = \dfrac{\psi(4k_v S + \psi)}{(8k_v S + 3\psi)S} \\ \nu = \dfrac{8k_v S + \psi}{8k_v S + 3\psi} \end{cases} \text{, (plane stress)} \tag{16}$$

where $\psi = \sum_{i=1}^{n} k_i L_i^2$ .

Equations (9), (14), and (15) are the governing equations for the element to satisfy the isotropic elastic constitutive law, where Eqs. (9) and (14) do not distinguish between plane stress and plane strain. For non-negative axial stiffness ($k_i \geq 0$ and $\psi \geq 0$), the lower bound of Poisson's ratio is $\nu$ > -1. As $k_v \gg \psi$ (dominant additional volumetric constraint), Poisson's ratio asymptotically approaches $\nu \rightarrow 1$ under plane stress and $\nu \rightarrow 1/2$ under plane strain. When $k_v = 0$ (negligible additional volumetric constraint), IELSM reverts to classical LSM with $\nu = 1/3$ (plane stress) and $\nu = 1/4$ (plane strain). When $k_v < 0$ (promote compression/expansion), the Poisson's ratio can range from -1 to 1/3 under plane stress and range from -1 to 1/4 under plane strain. Consequently, IELSM spans the physically admissible ranges: $-1 < \nu < 1$ for plane stress and $-1 < \nu < 1/2$ for plane strain, while maintaining exactly two parameters matching isotropic elasticity theory.

The above derivation starts from an abstract element that contains axial springs of arbitrary orientation and length, without predefining the number of springs, the connection topology, or the geometric outline of the element. This formal generality indicates that as long as the spring system parameters $\{\theta_i, L_i, k_i\}$ satisfy the aforementioned governing equations, the element can exactly reproduce the specified isotropic elastic constitutive law regardless of its final shape (including concave polygons). It should be noted that the choice of known and unknown parameters is not unique. Once the element geometry and spring topology are prescribed, a natural and convenient choice is to treat $\theta_i$ and $L_i$ as known and the spring stiffnesses $k_i$ as unknowns. This choice is adopted in the following derivation; in principle, other choices are also admissible, provided that the governing equations remain solvable.

### 2.2 Representative Tessellation Elements

When the above framework is applied to specific polygonal elements, the spring layout is determined by the element topology, and the spring stiffnesses are obtained by solving the governing equations from the previous section. For simplicity, the vertices of a polygonal element are treated as element nodes, and no internal nodes are introduced. This section takes representative IELSM elements such as triangles, quadrilaterals, and hexagons as examples to discuss the possibility of extending element shapes to irregular shapes, and finally presents a series of regular polygonal elements for subsequent fracture simulations. Note that the additional volumetric modulus in Eq. (15) is independent of the element shape; therefore, its value will not be reiterated in the subsequent derivations of different IELSM elements.

#### 2.2.1 Equilateral Triangular Element

As shown in FIGURE 2, the natural spring layout for a triangular element is along its three edges. Let the strain energy weights of the springs be the unknowns: $w_j = k_j L_j^2$ ($j$ = 1, 2, 3). However, the isotropic governing equations derived in Section 2.1 consist of five conditions, namely four rotational invariance equations (Eq. (9)) and one elastic parameter matching equation (Eq. (14)). For a general triangular geometry, the first four rotational invariance equations already constitute four independent real constraints (or two independent complex equations in complex form), whereas there are only three unknowns; hence, the number of equations exceeds the number of unknowns. For the system to have a solution, the orientation angles of the triangle must satisfy special symmetry conditions. Through the complex form, the first four equations can be transformed into $\sum w_j e^{i2\theta_j} = 0$ and $\sum w_j e^{i4\theta_j} = 0$, from which the weight ratios are required to be real, thereby forcing the interior angles of the triangle to satisfy a specific relation. It is derived (see Appendix A for details) that this condition is equivalent to the triangle being equilateral. For an equilateral triangle, the high degree of symmetry reduces the rank of the constraints, so that a nonzero solution exists; for any other triangle shape, the homogeneous equations admit only the zero solution and thus cannot satisfy the overall modulus matching requirement. Therefore, a triangular element based on three edge springs cannot be generalized to arbitrary shapes, which fundamentally explains why the classical LSM is limited to equilateral triangular elements.

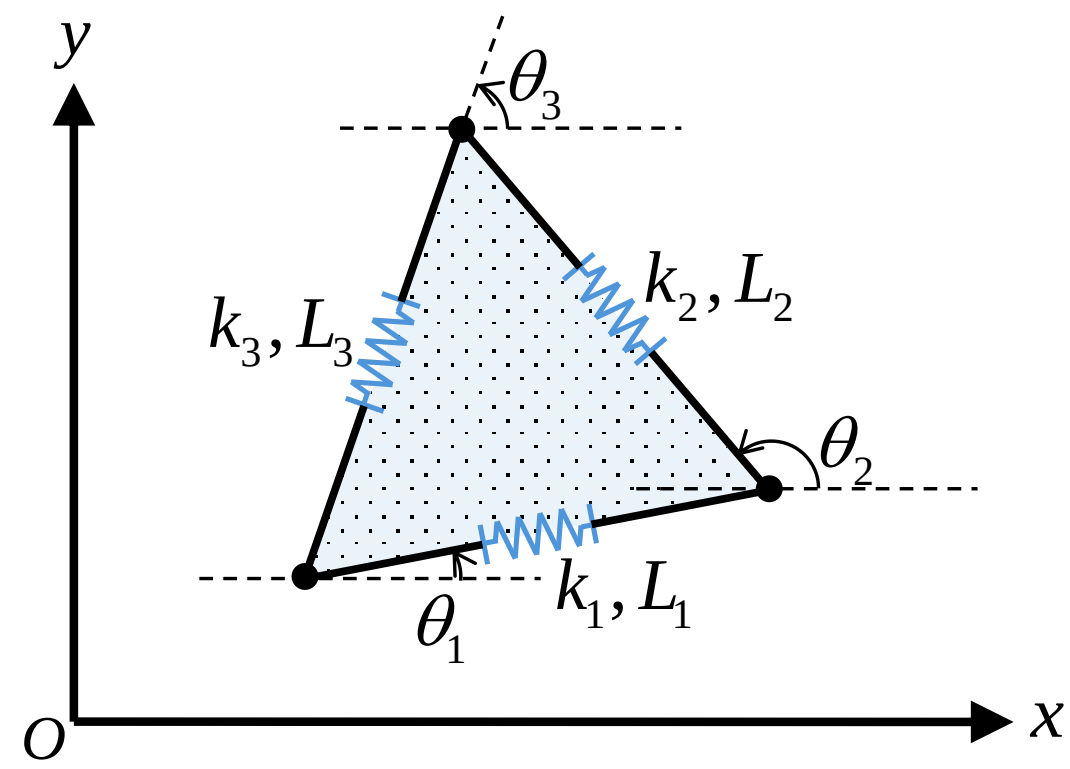

FIGURE 2 Schematic of the triangular IELSM element.

For an equilateral triangular element, assuming $L_i = L$, $\theta_3 = \theta_1 + \pi/3$, and $\theta_2 = \theta_3 + \pi/3$, substituting the geometric parameters into Eqs. (9) and (14) gives:

$$\begin{cases} L^2\left(k_1\cos(2\theta_1)+k_2\cos(2\theta_1+\frac{4}{3}\pi)+k_3\cos(2\theta_1+\frac{2}{3}\pi)\right)=0, \\ L^2\left(k_1\sin(2\theta_1)+k_2\sin(2\theta_1+\frac{4}{3}\pi)+k_3\sin(2\theta_1+\frac{2}{3}\pi)\right)=0, \\ L^2\left(k_1\cos(4\theta_1)+k_2\cos(4\theta_1+\frac{8}{3}\pi)+k_3\cos(4\theta_1+\frac{4}{3}\pi)\right)=0, \\ L^2\left(k_1\sin(4\theta_1)+k_2\sin(4\theta_1+\frac{8}{3}\pi)+k_3\sin(4\theta_1+\frac{4}{3}\pi)\right)=0, \\ L^2\left(k_1+k_2+k_3\right)=\frac{\sqrt{3}L^2}{4}\frac{4E}{1+\nu}. \end{cases} \tag{17}$$

Solving yields the relation between the spring stiffnesses $k_i$ and ($E$, $\nu$) as:

$$k_1=k_2=k_3=\frac{\sqrt{3}E}{3(1+\nu)}. \tag{18}$$

**2.2.2 Rectangular Element**

As shown in FIGURE 3, for an arbitrary quadrilateral with four non-coincident nodes, six springs are arranged along its four edges and two diagonals. Assuming the element geometry ($\theta_i$ and $L_i$) is known and the spring stiffnesses $k_i$ are unknown variables, there are six unknowns, which exceed the five governing equations; the system is therefore underdetermined and has infinitely many solutions. To obtain a unique solution, additional constraints must be introduced. For simplicity, the element shape is restricted here to a rectangle, and the symmetry thereby introduced reduces the rank of the governing equations and simplifies them. To derive a quadrilateral element of general shape, more relaxed constraints must be introduced. The choice of such constraints is crucial: one must identify exactly the missing constraints and ensure that the equation system admits all-positive solutions, which is generally very difficult. Therefore, this paper suggests introducing several internal nodes when deriving a general quadrilateral element, so as to add spring connections and thus enlarge the solution space of the underdetermined system; physical criteria such as minimum total spring stiffness can then be imposed, and the problem can be treated as a quadratic programming problem to obtain relatively uniform positive-stiffness solutions.

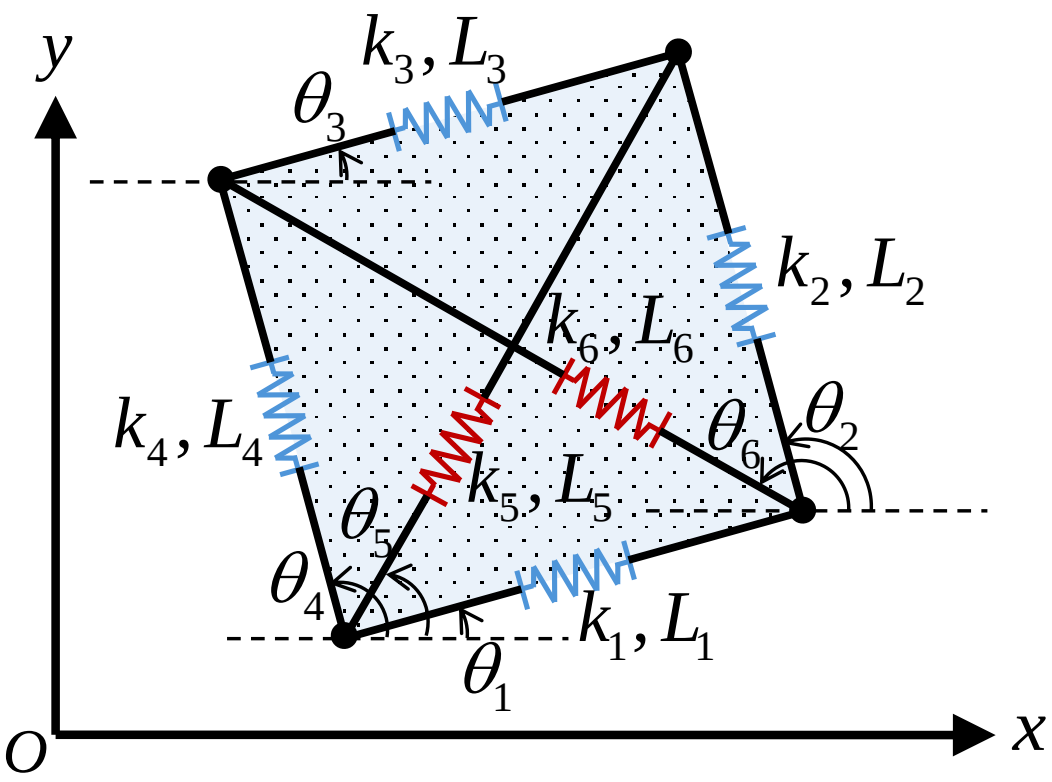

FIGURE 3 Schematic of the rectangular IELSM element.

Assuming the quadrilateral element in FIGURE 3 is rectangular, then $L_1 = L_3$, $L_2 = L_4$, $L_5 = L_6$, $\theta_2 = \theta_1 + \pi/2$, $\theta_3 = \theta_1$, $\theta_4 = \theta_1 + \pi/2$. Let the aspect ratio of the rectangle be $\alpha = L_2/L_1$; then $\theta_5 = \theta_1 + \arctan(\alpha)$ and $\theta_6 = \theta_1 + \arctan(-\alpha)$. Assuming the spring stiffnesses satisfy $k_1 = k_3$, $k_2 = k_4$, $k_5 = k_6$, substituting the above geometric parameters and stiffness symmetry conditions into Eqs. (9) and (14) gives:

$$\begin{cases} 2L^2\left(k_1+k_5-(k_2+k_5)\alpha^2\right)\cos 2\theta_1 = 0, \\ 2L^2\left(k_1+k_5-(k_2+k_5)\alpha^2\right)\sin 2\theta_1 = 0, \\ \dfrac{2L^2}{1+\alpha^2}\left((k_2+k_5)\alpha^4+(k_1+k_2-6k_5)\alpha^2+k_1+k_5\right)\cos 4\theta_1 = 0, \\ \dfrac{2L^2}{1+\alpha^2}\left((k_2+k_5)\alpha^4+(k_1+k_2-6k_5)\alpha^2+k_1+k_5\right)\sin 4\theta_1 = 0, \\ L^2\left(k_1+k_2+k_3+k_4+k_5+k_6\right) = \dfrac{4EL^2}{1+\nu}. \end{cases} \tag{19}$$

Solving yields the relation between the spring stiffnesses $k_i$ and $(E, \nu)$ as:

$$k_1 = k_3 = \frac{E(3\alpha^2-1)}{4(1+\nu)\alpha^2}, k_2 = k_4 = \frac{E(3-\alpha^2)}{4(1+\nu)\alpha^2}, k_5 = k_6 = \frac{E(1+\alpha^2)}{4(1+\nu)\alpha^2}. \tag{20}$$

From Eq. (20), the restriction on the rectangular shape that makes all spring stiffnesses positive can be derived:

$$\frac{\sqrt{3}}{3} < \alpha < \sqrt{3}. \tag{21}$$

This indicates that, for all axial spring stiffnesses to be positive, the rectangular element must not be too slender. Moreover, as $\alpha$ increases, the stiffnesses of the four edge springs ($k_1 - k_4$) become smaller relative to those of the diagonal springs ($k_5$, $k_6$); as $\alpha$ approaches the limit, the edge spring stiffnesses approach zero. In other words, when the rectangular element is too slender, its overall stiffness matrix tends to become singular (rank-deficient), leading to numerical instability. Therefore, this paper considers only the special case of a rectangular element, namely the square element ($\alpha = 1$), whose advantages in computational efficiency and accuracy have been demonstrated in our previous work [13].

### 2.2.3 Regular Hexagonal Element

For an arbitrary hexagonal element, the number of unknown spring stiffnesses is far greater than the number of isotropic governing equations (five), so the system is highly underdetermined. Therefore, this paper restricts the element to the regular hexagonal case and uses symmetry to provide an analytical reference. Owing to its high geometric symmetry, the regular hexagon is another important element shape in lattice spring models. As shown in FIGURE 4, for a regular hexagon, the springs can be divided into three categories according to their connectivity: the first comprises the six edge springs, with length $L$ and stiffness $k_1$; the second comprises springs connecting opposite vertices (i.e., long diagonals), with length $2L$ and stiffness $k_2$; and the third comprises springs connecting vertices separated by one vertex (i.e., short diagonals), with length $\sqrt{3}L$ and stiffness $k_3$. Different combinations of these three spring categories yield the following three representative schemes:

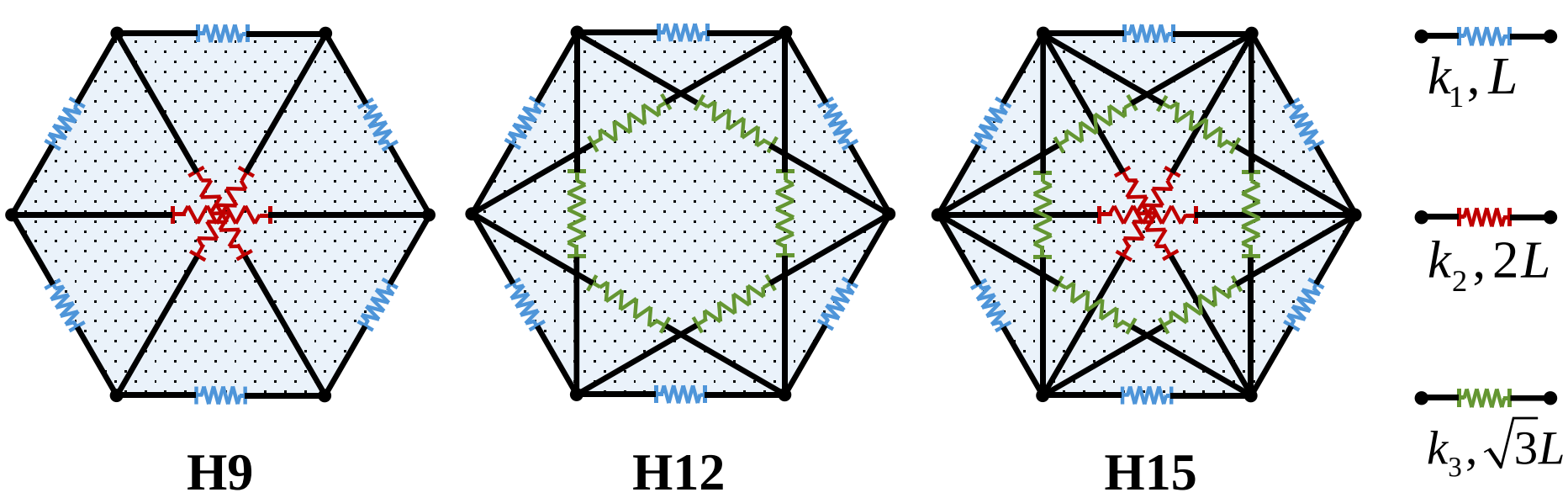

FIGURE 4 Schematic of the regular hexagonal IELSM element.

1. Scheme 1 (H9): contains only the first and second categories, i.e., six edges plus three long diagonals. The number of springs in this arrangement is exactly equal to the number of deformation degrees of freedom, so it is nominally statically determinate; however, the high symmetry of the regular hexagon makes the spring constraints linearly dependent, and an hourglass zero-energy mode actually exists (see the eigenvalue analysis in Section 3.4).

2. Scheme 2 (H12): contains the first and third categories, i.e., six edges plus six short diagonals. This scheme forms a triangular mesh within the element and has three more springs than Scheme 1, corresponding to an internally statically indeterminate arrangement.

3. Scheme 3 (H15): contains all three categories simultaneously, i.e., all edges, all short diagonals, and all long diagonals, totaling 15 springs. It is the union of the first two schemes and provides the richest constraints.

Although the numbers of springs differ, under the highly symmetric geometry of the regular hexagon, all three schemes can be incorporated into a unified algebraic framework. As in the solution of spring stiffnesses for the rectangular element, the spring orientation angles and lengths are fully determined by the regular hexagonal geometry, and the spring stiffnesses of the three schemes can be solved from the governing equations in Section 2.1:

$$\begin{cases} k_1 = k_2 = \dfrac{\sqrt{3}E}{3(1+\nu)}, & \text{H9} \\ k_1 = k_3 = \dfrac{\sqrt{3}E}{4(1+\nu)}, & \text{H12} \\ k_1 = k_2 = k_3 = \dfrac{\sqrt{3}E}{6(1+\nu)}, & \text{H15} \end{cases} \tag{22}$$

## 2.3 Other Element Shapes

In addition to the three representative element shapes in Section 2.2, the theory in Section 2.1 can be readily used to derive elements of other shapes, especially regular polygons. The derivation for arbitrary regular polygonal elements is given in Appendix B. For any regular $N$-gon ($N \geq 5$), as long as the springs are grouped by chord order $\delta$, where $\delta$ denotes the number of edges between the two endpoints of the chord, and the strain-energy weights $w_j$ within the same group are equal, the rotational invariance conditions are automatically satisfied owing to the uniform distribution of orientation angles, without imposing any additional restriction on the spring stiffnesses. The only remaining condition to be satisfied is the overall modulus matching condition. Under the equal-stiffness assumption, an analytical stiffness formula can be obtained, which provides a unified theoretical basis for parameter calibration of regular polygonal IELSM elements. FIGURE 5(a) shows a regular pentagonal element with five nodes and ten springs. Let springs of equal length within the same group have the same stiffness, and assume that the two groups of springs have equal stiffness. According to Eq. (B.13), one set of solutions for the spring stiffnesses in the regular pentagonal element is:

$$k_1 = k_2 = \frac{2E\sin(\frac{2\pi}{5})}{5(1+\nu)}. \quad \text{(Regular pentagon)} \tag{23}$$

However, a regular pentagon cannot be seamlessly tessellated (see FIGURE 5(c)), so it has not appeared in the LSM literature, and the same holds for other regular polygonal elements.

Nevertheless, these regular polygonal elements that cannot be tessellated by a single shape are not entirely inapplicable: seamless tessellations, either periodic or non-periodic, can be achieved by combining elements of different shapes [18]. For example, FIGURE 5(e) shows a schematic of a tessellation mesh formed by regular dodecagonal, regular hexagonal, and square elements. Therefore, the element shapes usable in IELSM can be any combination of shapes that satisfies tessellation. In addition, the present theory can also be applied to concave polygons. As shown in FIGURE 5(b), an L-shaped element has six vertices and nine springs. Assuming its longest edge length is $L$, springs of the same length have the same stiffness, and setting $k_1 = k_2$ or $k_1 = k_3$ gives two sets of solutions for the spring stiffnesses:

$$\text{L-shaped element}\begin{cases} k_1 = k_2 = \dfrac{E}{2(1+\nu)}, k_3 = \dfrac{E}{1+\nu}; \\ k_1 = k_3 = \dfrac{E}{1+\nu}, k_2 = \dfrac{E}{4(1+\nu)}. \end{cases} \tag{24}$$

The construction of this L-shaped element is essentially inherited from the six-spring square element: the orientation combinations of the internal springs are identical, so the L-shaped element can still satisfy the isotropy conditions. On the one hand, this element type provides strong evidence that the IELSM theory supports concave polygons; on the other hand, it can serve as a transition element for mesh refinement of rectangular elements, as shown in FIGURE 5(d). This demonstrates that the generalized IELSM theory developed in this paper has broad applicability and still has room for further extension. However, to focus on the application of IELSM to fracture simulation, the following sections consider only the three element types in Section 2.2 that can be seamlessly tessellated by a single shape.

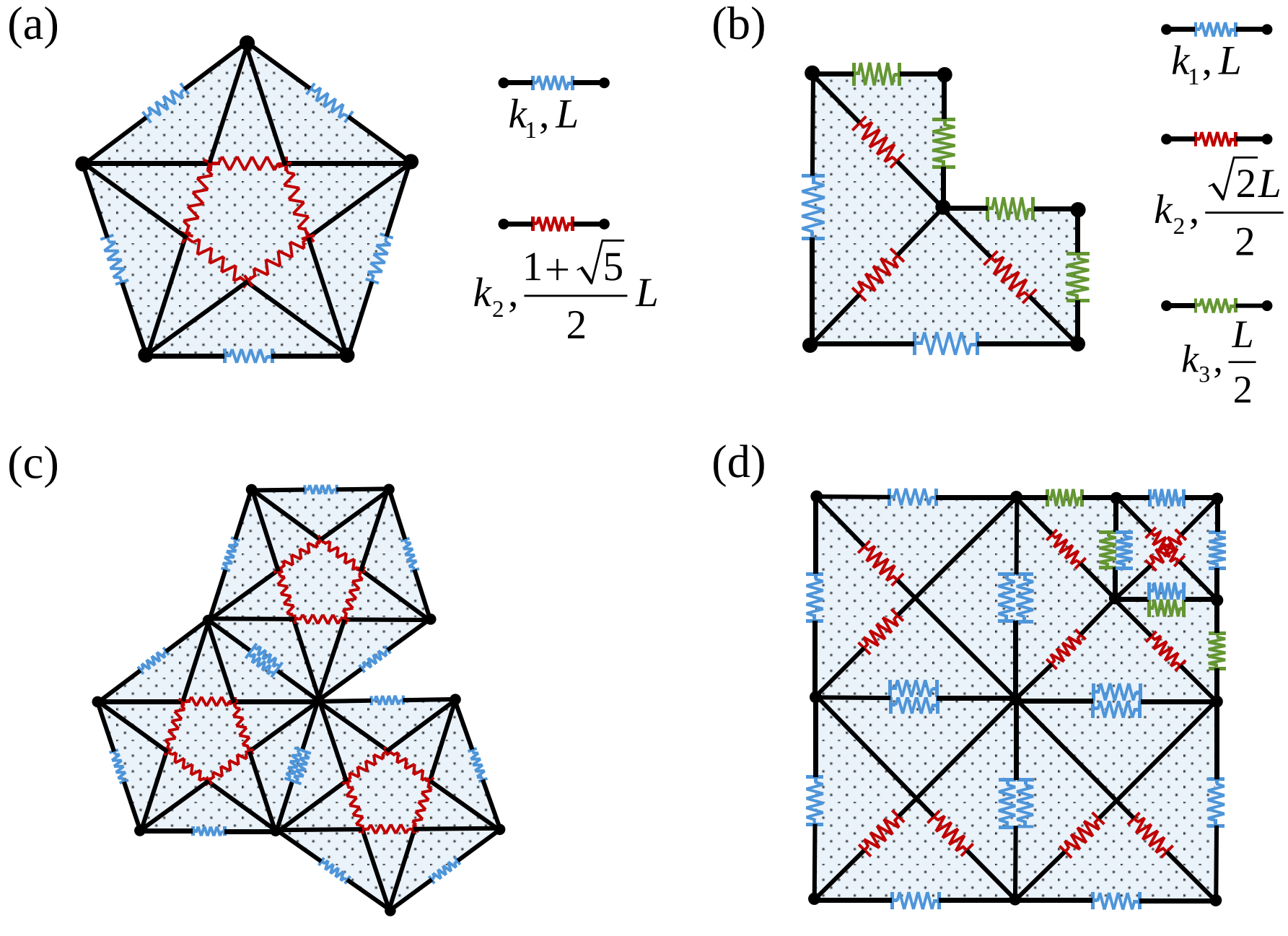

(e)

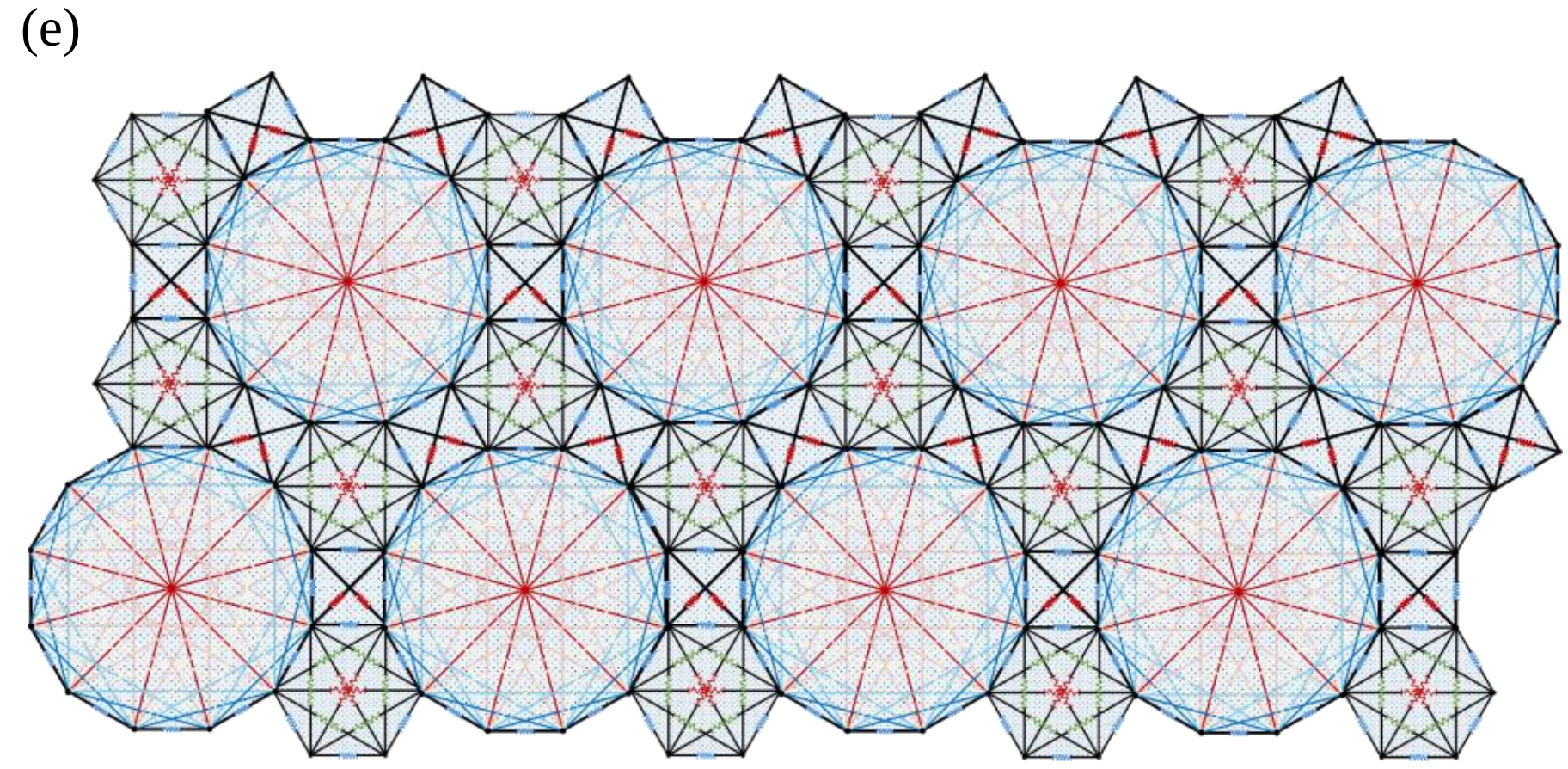

FIGURE 5 Other element shapes: (a) regular pentagonal element, (b) L-shaped element, (c) regular pentagon cannot be tessellated, (d) L-shaped element used as a transition element for square mesh refinement, (e) schematic of a tessellation mesh formed by regular dodecagonal, regular hexagonal, and square elements.

## 3. Numerical Implementation

### 3.1 Derivation of the Element Stiffness Matrix

The element stiffness matrix of the *i*-th spring can be obtained from the element stiffness matrix of a truss element in the finite element method:

$$\mathbf{K}_{n,i}(k_i,\theta_i)=\mathbf{T}^{\mathrm{T}}\mathbf{K}_{n,i}^{e}\mathbf{T} \tag{25}$$

where $\mathbf{K}_{n,i}^{e}$ is the element stiffness matrix of the truss in the local coordinate system, and $\mathbf{T}$ is the two-dimensional coordinate rotation matrix, given respectively by:

$$\mathbf{K}_{n,i}^{e}=\begin{bmatrix} k_i & 0 & -k_i & 0 \\ 0 & 0 & 0 & 0 \\ -k_i & 0 & k_i & 0 \\ 0 & 0 & 0 & 0 \end{bmatrix},$$
$$\mathbf{T}=\begin{bmatrix} \cos\theta & \sin\theta & 0 & 0 \\ -\sin\theta & \cos\theta & 0 & 0 \\ 0 & 0 & \cos\theta & \sin\theta \\ 0 & 0 & -\sin\theta & \cos\theta \end{bmatrix}. \tag{26}$$

The discretization of the additional volumetric constraint follows an approach similar to that in [13]. On this basis, this paper generalizes the additional volumetric constraint element stiffness matrix to arbitrary element shapes. As shown in FIGURE 6, consider an arbitrary polygon with $n$ nodes numbered counterclockwise. The coordinates of the $i$-th node before deformation are $P_i(X_i, Y_i)$, $i$ = 1, 2, …, $n$, and those after deformation are $P_i'(x_i, y_i)$. The displacement of the $i$-th node is $(u_i, v_i)$, with $X_i = x_i - u_i$ and $Y_i = y_i - v_i$. Then the polygon areas before and after deformation are:

$$S = \frac{1}{2}\sum_{i=1}^{n}(X_iY_{i+1} - X_{i+1}Y_i),$$
$$S_d = \frac{1}{2}\sum_{i=1}^{n}(x_iy_{i+1} - x_{i+1}y_i) = \frac{1}{2}\sum_{i=1}^{n}\big((X_i + u_i)(Y_{i+1} + v_{i+1}) - (X_{i+1} + u_{i+1})(Y_i + v_i)\big). \tag{27}$$

where $X_{n+1} = X_1$, $Y_{n+1} = Y_1$, $x_{n+1} = x_1$, $y_{n+1} = y_1$, $u_{n+1} = u_1$, $v_{n+1} = v_1$.

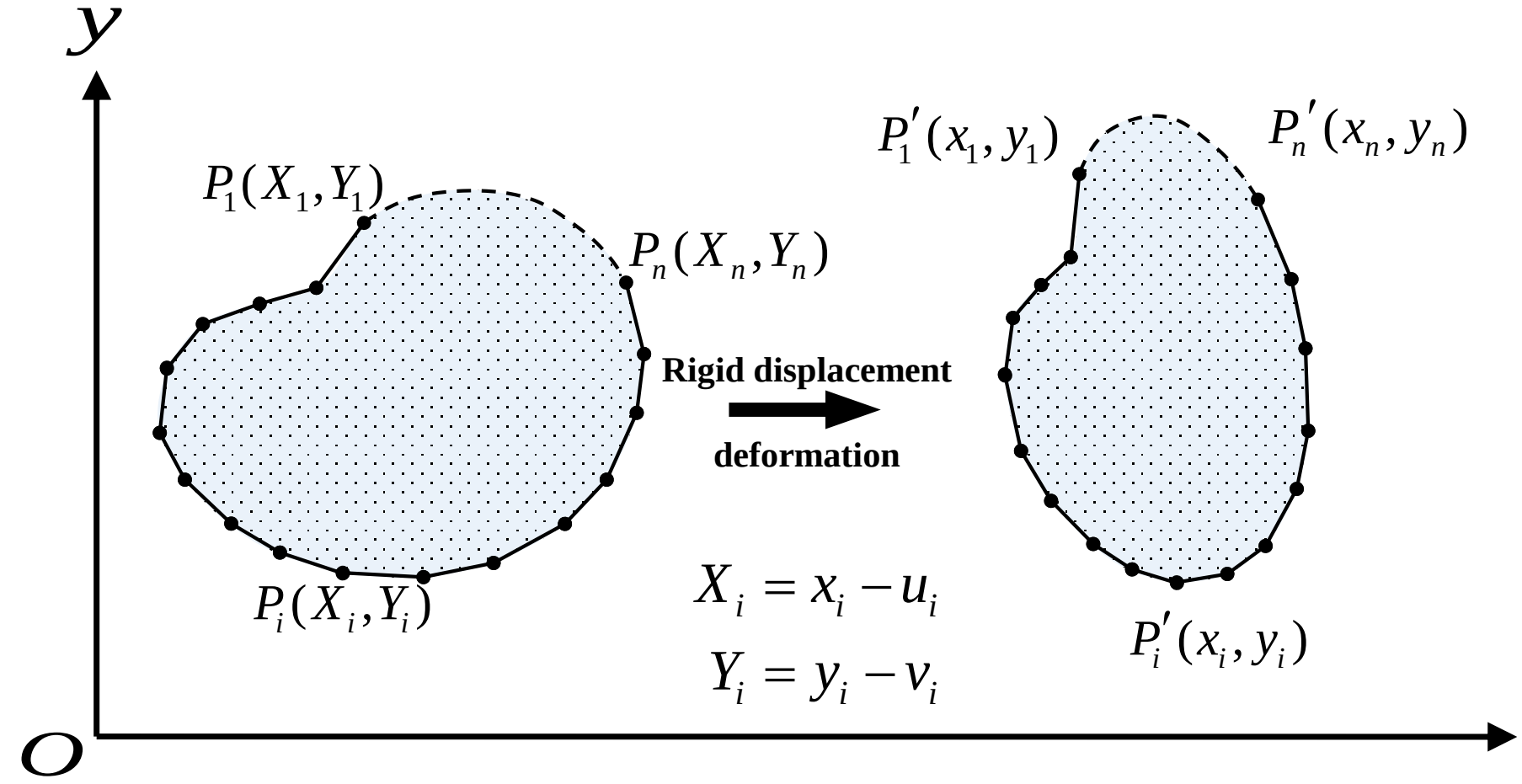

FIGURE 6 Schematic of rigid-body displacement and deformation of a polygon.

The volumetric strain of the polygon is:

$$\varepsilon_v = \frac{S_d - S}{S} = \frac{1}{2S}\sum_{i=1}^{n}\big((X_iv_{i+1} + Y_{i+1}u_i - X_{i+1}v_i - Y_iu_{i+1}) + (u_iv_{i+1} - u_{i+1}v_i)\big). \tag{28}$$

The displacement-related linear term in the volumetric strain is:

$$\begin{aligned}\varepsilon_{v,linear} &= \frac{1}{2S}\sum_{i=1}^{n}\left(X_iv_{i+1} + Y_{i+1}u_i - X_{i+1}v_i - Y_iu_{i+1}\right)\\ &= \frac{1}{2S}\mathbf{cU},\end{aligned} \tag{29}$$

where

$$\mathbf{c} = \begin{bmatrix} Y_2 - Y_n & X_n - X_2 & Y_3 - Y_1 & X_1 - X_3 & \cdots & Y_1 - Y_{n-1} & X_{n-1} - X_1 \end{bmatrix}_{1\times 2n};$$
$$\mathbf{U} = \begin{bmatrix} u_1 & v_1 & u_2 & v_2 & \cdots & u_n & v_n \end{bmatrix}_{2n\times 1}^T. \tag{30}$$

It is pointed out in [13] that under the small-deformation assumption, the nonlinear term in the volumetric strain can be omitted. Then the strain energy of the additional volumetric constraint is:

$$
\begin{aligned}
U_{volume} &\simeq \frac{k_v \varepsilon_{v,linear}{}^2 S}{2} \\
&= \frac{k_v S}{2}\left(\frac{1}{2S}\mathbf{cU}\right)^2 \\
&= \frac{k_v S}{2}\frac{1}{4S^2}(\mathbf{cU})^T(\mathbf{cU}) \\
&= \frac{k_v}{8S}\mathbf{U}^T(\mathbf{c}^T\mathbf{c})\mathbf{U},
\end{aligned}
\tag{31}
$$

this expression leads to the additional volumetric constraint stiffness matrix:

$$
\mathbf{K}_{volume} = \frac{\partial^2 W_{volume}}{\partial \mathbf{U}^2} = \frac{k_v}{4S}\mathbf{c}^T\mathbf{c}. \tag{32}
$$

### 3.2 Hybrid Assembly Strategy for IELSM and Finite Elements

Coupling IELSM with standard finite elements makes it possible to restrict IELSM to the damage-prone region while modeling the rest of the domain with conventional finite elements. The key requirement for such hybrid meshes is $C^0$ continuity along shared edges. ELSM elements are constructed from polygonal corner nodes and adopt a linear displacement field. Along a straight edge, their displacement interpolation is therefore linear and determined only by the two edge nodes. This boundary interpolation is identical to that of CST and Q4 elements. Consequently, IELSM can be directly coupled with CST and Q4 elements through node sharing, without interface elements, overlapping regions, or kinematic constraints. IELSM elements cannot, however, be directly mixed with higher-order finite elements through shared nodes; without special treatment, such a mix would produce non-conforming meshes. The numerical examples in this paper consider IELSM–CST coupling. The compatibility of IELSM and Q4 elements through node sharing follows from the same boundary-interpolation argument, while its numerical verification is left for future work.

Based on the above boundary displacement compatibility, this paper treats IELSM elements as standard superelements with $n$ nodes and $2n$ degrees of freedom, and directly adopts the global stiffness matrix assembly procedure of standard finite elements for hybrid mesh integration, as illustrated in FIGURE 7. Regardless of whether an element comes from IELSM or conventional CST or Q4 elements, its local stiffness matrix is directly added to the global stiffness matrix according to the same degree-of-freedom mapping. When assembling the stiffness matrix in the global coordinate system, the entire solution domain $\Omega$ is first divided into the IELSM element domain $\Omega_{IELSM}$ and the FEM element domain $\Omega_{FEM}$; then, according to the global nodal degree-of-freedom indices corresponding to each type of element, their stiffness matrices are directly added to the global stiffness matrix:

$$\mathbf{K}_{\text{Global}} = \sum_{e\in\Omega_{\text{IELSM}}} \mathbf{K}^{e}_{\text{IELSM}} + \sum_{e\in\Omega_{\text{FEM}}} \mathbf{K}^{e}_{\text{FEM}}, \tag{33}$$

where

$$\mathbf{K}^{e}_{\text{IELSM}} = \sum_{i=1}^{n} \mathbf{K}_{n,i} + \mathbf{K}^{e}_{volume}. \tag{34}$$

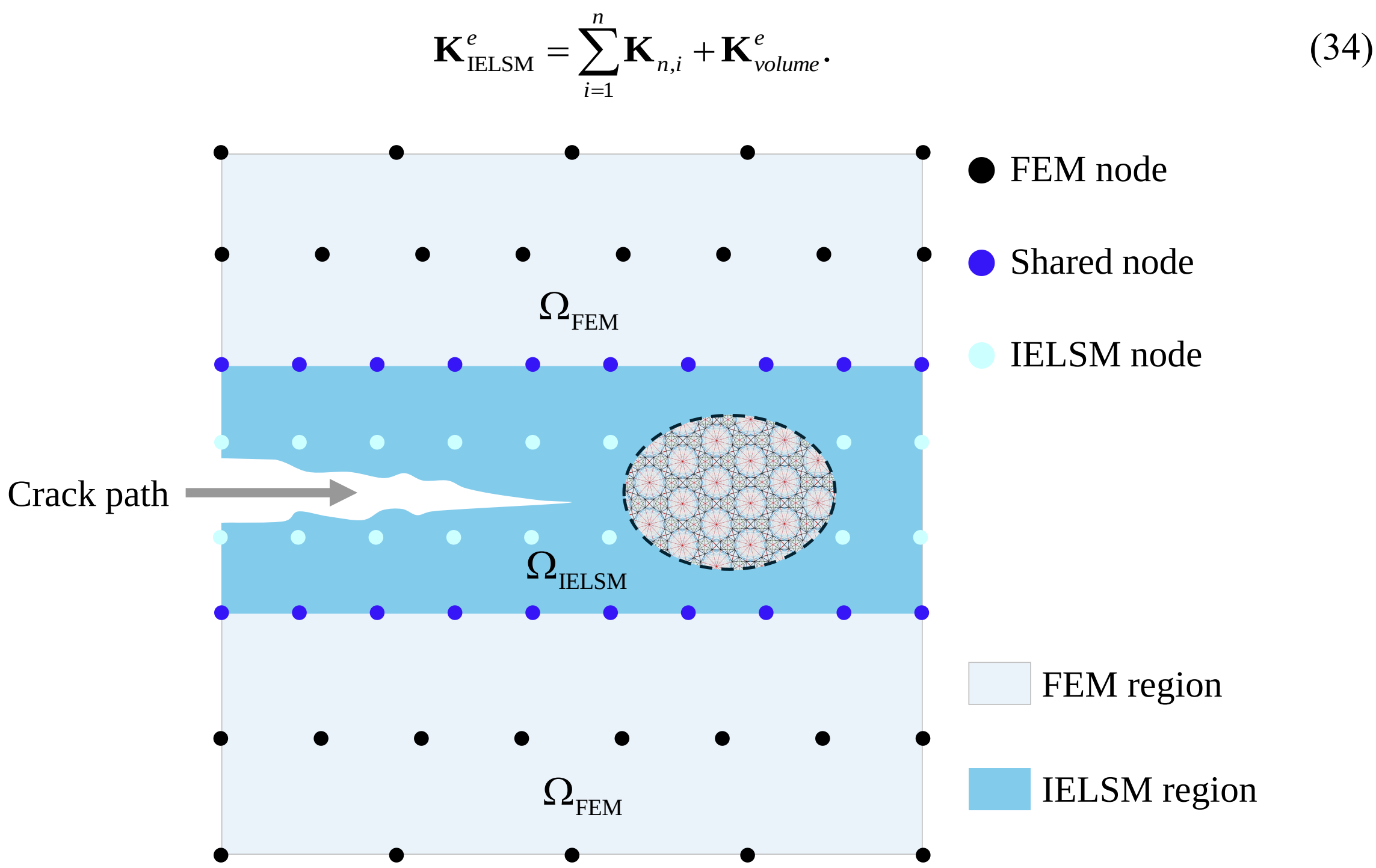

FIGURE 7 Schematic of IELSM-FEM coupling

It should be noted that, in our previous work [13], the additional volumetric constraint was equivalently decomposed into axial, tangential, and rotational springs, which significantly reduces the computational cost of element assembly. Such a decomposition is not adopted in the present study for two reasons. First, the assembly of the global stiffness matrix has only a minor influence on the total computational time of the present simulations. Second, deriving a decomposition strategy applicable to arbitrary regular polygonal elements is relatively cumbersome. Therefore, this decomposition is left for future work.

### 3.3 Stress/Strain Fields for IELSM

Although the derivation of the IELSM element stiffness matrix does not rely on an explicit construction of shape functions, its stress/strain field can still be post-processed using the displacement-strain transformation matrix of standard finite elements. As stated in [13], when reconstructing the displacement field in the region without physical elements within the RVE, the linear displacement field assumption consistent with the axial spring displacement mode is the most physically sound; however, it is also feasible to compute stresses and strains based on a bilinear displacement field assumption, which belongs to the classical stress recovery techniques in finite element theory for improving the continuity of local stress/strain fields. To maintain consistency in post-processing and subsequent discussions, for the equilateral triangular, square, and regular hexagonal elements presented in Section 2.2, this paper uniformly derives their stress/strain fields based on the linear displacement field assumption.

As shown in FIGURE 8, for rectangular and regular hexagonal elements, this paper introduces the geometric centroid of the element as a support node, dividing the rectangular region into four triangles and the regular hexagonal region into six triangles. These subregions use the strain–displacement matrix blocks of CST elements to solve for strains. It should be emphasized that the support node is not a physical node with independent degrees of freedom; its displacement is directly computed as the average of the displacements of the vertices of the corresponding element. Compared with the simple partition of a rectangular element into two triangles along a single diagonal, the radial centroid partition adopted here can effectively avoid the topological anisotropy introduced by the choice of a single diagonal direction. More importantly, based on the averaging assumption for the centroid displacement, the strain tensor constructed by this block reconstruction method naturally satisfies symmetry, while the radial partition avoids interference from artificially biased deformation modes in strain extraction. This block CST linear interpolation method not only preserves the physical continuity of the low-order displacement field at interfaces but also greatly simplifies the strain computation for regular polygons, remaining theoretically consistent with the displacement field assumption of IELSM itself.

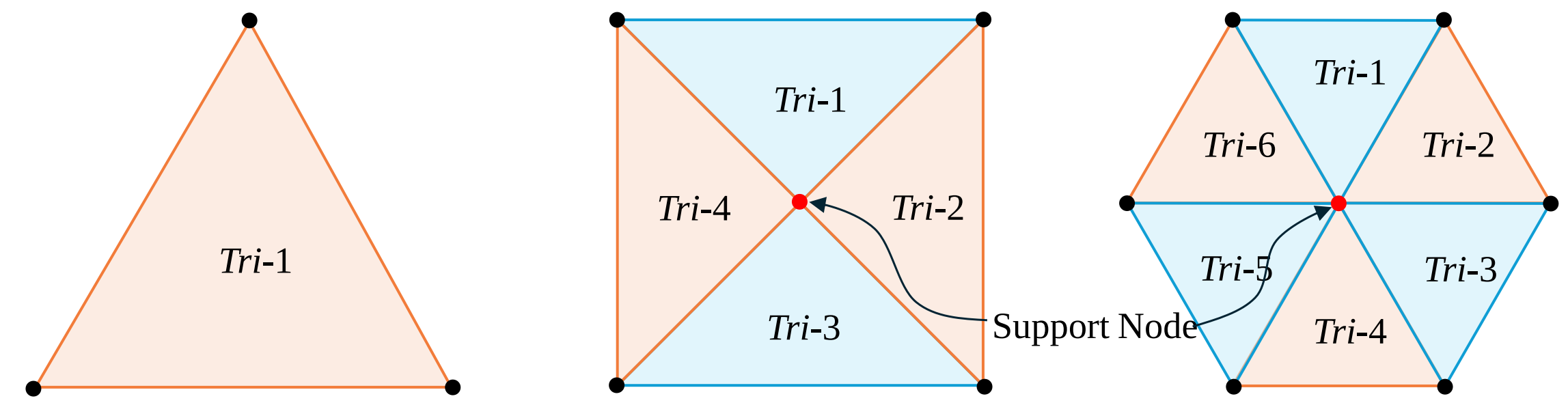

FIGURE 8 Schematic of the radial centroid partition for rectangular and regular hexagonal elements.

After solving the strain components of all sub-triangles within an IELSM element (one for a triangular element, four for a rectangular element, and six for a regular hexagonal element), the nodal strain is taken as the arithmetic average of the strains of all sub-triangles sharing that node.

Based on the physical equations of plane problems in elasticity, the stress field of a representative element can be further obtained through the constitutive relation. However, this post-processing scheme has certain limitations: it has been verified that when the equilateral triangular (TRI), square (SQ), and regular hexagonal (HEX) IELSM elements are mixed in a tessellation, the displacement field accuracy is unaffected, but the stress/strain field accuracy decreases. This phenomenon arises because the block CST recovery relies on the radial centroid partition of each element itself: the strains of the sub-triangles do not share the same interpolation along the common edges of adjacent elements; when adjacent elements have different shapes, the strains recovered on the two sides are inconsistent at the boundary, thereby producing errors (see FIGURE 13). Therefore, when different IELSM elements are mixed, more robust post-processing techniques are needed, and strain smoothing methods such as the smoothed finite element method [46] are promising approaches to solve this problem.

**3.4 Eigenvalue Analysis**

The advantages of the IELSM square element stiffness matrix in terms of numerical stability have been reported in [13]. Compared with CST and Q4 elements in FEM, its modes contain no parasitic shear or hourglass zero-energy modes, which enables it to achieve both computational accuracy and efficiency on structured meshes. In this section, taking the plane strain condition as an example, eigenvalue analyses are performed for the element stiffness matrices of the IELSM equilateral triangular element and the three regular hexagonal elements presented in Section 2.2.

The nodes of each element are numbered counterclockwise. The nodal coordinates of the equilateral triangular element are $P_1(0, 0)$, $P_2(L, 0)$, $P_3(L/2, \sqrt{3}L/2)$; those of the regular hexagonal element are $P_1(0, 0)$, $P_2(L, 0)$, $P_3(3L/2, \sqrt{3}L/2)$, $P_4(L, \sqrt{3}L)$, $P_5(0, \sqrt{3}L)$, $P_6(-L/2, \sqrt{3}L/2)$. Using the additional bulk modulus $k_v$ and axial spring stiffnesses $k_i$ given in Eqs. (15), (18), (20), and (22), and following the assembly scheme in Eq. (34), the element stiffness matrices $\mathbf{K}_{\mathrm{IELSM,\,Tri}}$, $\mathbf{K}_{\mathrm{IELSM,\,H9}}$, $\mathbf{K}_{\mathrm{IELSM,\,H12}}$, and $\mathbf{K}_{\mathrm{IELSM,\,H15}}$ can be obtained for each element type. The detailed forms of the element stiffness matrices are given in Appendix C. For plane stress conditions, the corresponding stiffness matrices can be directly obtained by replacing $E$ and $v$ in the plane-strain expressions of Appendix C with $E(1 + 2v)/(1 + v)^2$ and $v/(1 + v)$, respectively.

To evaluate the numerical stability and deformation characteristics of the four types of element stiffness matrices, it is necessary to solve the corresponding standard eigenvalue problem. The specific solution equation is:

$$\mathbf{K}\phi_i = \lambda_i \phi_i, \tag{35}$$

where **K** represents the stiffness matrix to be analyzed, $\lambda_i$ indicates the $i$-th eigenvalue, and $\phi_i$ represents the corresponding eigenvector.

The modes of the four element stiffness matrices and their corresponding eigenvalues are shown in TABLE 1. All elements share the same volumetric deformation mode and three rigid-body displacement modes (translation and rotation), which indicates that IELSM does not generate spurious elastic forces under rigid-body displacements and also provides a direct verification of rotational invariance. For all four element stiffness matrices, every nonzero eigenvalue is greater than 0 within the Poisson's ratio range of (-1, 0.5), indicating that each stiffness matrix is positive semidefinite. For the IELSM regular hexagonal elements, the eigenvalue distributions of the element stiffness matrices for the three spring arrangements show significant differences: $\mathbf{K}_{\text{IELSM, H9}}$ has an additional zero eigenvalue, which means that the stiffness matrix is rank-deficient and that a non-physical zero-energy hourglass mode exists in the element, so numerical stability cannot be guaranteed; $\mathbf{K}_{\text{IELSM, H12}}$ has no zero-energy mode, but apart from the volumetric deformation mode, its nonzero eigenvalues exhibit two different magnitudes, indicating that the discrete system exhibits a certain mesh-induced stiffness anisotropy in different shear directions; $\mathbf{K}_{\text{IELSM, H15}}$ performs best: owing to the largest number of intra-element springs, all its nonzero modal eigenvalues except the volumetric deformation mode are exactly equal. This indicates that the stiffness matrix of H15 completely eliminates mesh-induced stiffness anisotropy at the modal level, and its behavior is consistent with that of a continuous isotropic medium.

TABLE 1 Comparison of stiffness matrix eigenvalues and modal shapes (plane strain).

| | $\lambda_1$ | $\lambda_2$ | $\lambda_3$ | $\lambda_4$ | $\lambda_5$ | $\lambda_6$ | $\lambda_7$ | $\lambda_8$ | $\lambda_9$ | $\lambda_{10}$ | $\lambda_{11}$ | $\lambda_{12}$ |
|---|---|---|---|---|---|---|---|---|---|---|---|---|
| $\mathbf{K}_{\text{IELSM, Tri}}$ | $\frac{\sqrt{3}E}{2(1-2\nu)(1+\nu)}$ | $\frac{\sqrt{3}E}{2(1+\nu)}$ | $\frac{\sqrt{3}E}{2(1+\nu)}$ | 0 | 0 | 0 | - | - | - | - | - | - |
| | | | | | | | - | - | - | - | - | - |
| $\mathbf{K}_{\text{IELSM, H9}}$ | $\frac{\sqrt{3}E}{2(1-2\nu)(1+\nu)}$ | $\frac{\sqrt{3}E}{1+\nu}$ | $\frac{\sqrt{3}E}{1+\nu}$ | $\frac{\sqrt{3}E}{1+\nu}$ | $\frac{\sqrt{3}E}{2(1+\nu)}$ | $\frac{\sqrt{3}E}{2(1+\nu)}$ | $\frac{\sqrt{3}E}{2(1+\nu)}$ | $\frac{\sqrt{3}E}{2(1+\nu)}$ | 0 | 0 | 0 | 0 |
| $\mathbf{K}_{\text{IELSM, H12}}$ | $\frac{\sqrt{3}E}{2(1-2\nu)(1+\nu)}$ | $\frac{3\sqrt{3}E}{4(1+\nu)}$ | $\frac{3\sqrt{3}E}{4(1+\nu)}$ | $\frac{3\sqrt{3}E}{4(1+\nu)}$ | $\frac{3\sqrt{3}E}{4(1+\nu)}$ | $\frac{3\sqrt{3}E}{4(1+\nu)}$ | $\frac{3\sqrt{3}E}{4(1+\nu)}$ | $\frac{\sqrt{3}E}{4(1+\nu)}$ | $\frac{\sqrt{3}E}{4(1+\nu)}$ | 0 | 0 | 0 |
| $\mathbf{K}_{\text{IELSM, H15}}$ | $\frac{\sqrt{3}E}{2(1-2\nu)(1+\nu)}$ | $\frac{\sqrt{3}E}{2(1+\nu)}$ | $\frac{\sqrt{3}E}{2(1+\nu)}$ | $\frac{\sqrt{3}E}{2(1+\nu)}$ | $\frac{\sqrt{3}E}{2(1+\nu)}$ | $\frac{\sqrt{3}E}{2(1+\nu)}$ | $\frac{\sqrt{3}E}{2(1+\nu)}$ | $\frac{\sqrt{3}E}{2(1+\nu)}$ | $\frac{\sqrt{3}E}{2(1+\nu)}$ | 0 | 0 | 0 |

## 4. IELSM with isotropic damage model

In the present model, an internal damage variable $d$ is recorded at each node, and the damage of springs and additional volumetric constraints is weighted by the damage of all nodes to which they are connected. This section presents the combination of IELSM with an isotropic damage model. First, the constitutive relation of the isotropic damage model is defined as:

$$\boldsymbol{\sigma} = (1-d)\bar{\boldsymbol{\sigma}}, \tag{36}$$

where $d$ is a scalar damage index ranging from 0 to 1, $\boldsymbol{\sigma}$ is the stress tensor, and $\bar{\boldsymbol{\sigma}} = \mathbf{C}_0 : \boldsymbol{\varepsilon}$ is the effective stress tensor in the undamaged state. $\mathbf{C}_0$ is the fourth-order elastic tensor and $\boldsymbol{\varepsilon}$ is the strain tensor. In this paper, there is no difference between computing damage using strain or stress, because the stress field herein is computed from the strain field.

For the $i$-th spring, assuming its connected nodes are $I$ and $J$ with damage values $d_I$ and $d_J$, respectively, the degraded stiffness of the spring is:

$$k_i^d = \left(1 - \frac{d_I + d_J}{2}\right) k_i. \tag{37}$$

For the additional volumetric constraint within the $j$-th IELSM element, assuming it is connected to $N$ nodes, the degraded additional bulk modulus is:

$$k_{v,j}^d = \left(1 - \frac{1}{N}\sum_{m=1}^{N} d_m\right) k_{v,j}. \tag{38}$$

The general form of the damage initiation function is defined as:

$$f\left(\bar{\sigma}_{eq}(\bar{\boldsymbol{\sigma}}), r\right) = \bar{\sigma}_{eq}(\bar{\boldsymbol{\sigma}}) - r, \tag{39}$$

where $\bar{\sigma}_{eq}(\bar{\boldsymbol{\sigma}})$ is a scalar equivalent effective stress determined by the yield surface, and $r$ is the current damage threshold, which represents the maximum value of the equivalent effective stress ever reached in the loading history. The initial value of $r$ is the uniaxial tensile strength $f_t$ of the material. To ensure the irreversibility of damage, $r$ is computed as:

$$r = \max\left(f_t, \max\left(\bar{\sigma}_{eq}\right)\right). \tag{40}$$

The evolution of the internal damage variable $d$ used in this study follows an exponential law, which has been widely used to simulate the post-peak stress-strain response of quasi-brittle materials characterized by gradual softening with a long tail [33,38,47–51], and takes the form:

$$d = 1 - \frac{f_t}{r}\exp\left(2H_d\left(1 - \frac{r}{f_t}\right)\right), \tag{41}$$

where the parameter $H_d$ is given by:

$$H_d = \frac{l^*}{2l_{\mathrm{Irwin}} - l^*}, \tag{42}$$

where $l_{\text{Irwin}}$ is Irwin's characteristic length, equal to $l_{\text{Irwin}}=EG_f/(f_t)^2$, $G_f$ is the fracture energy, $l^*$ is the characteristic length. In the literature [49–51], the characteristic length scale in FEM-type works is estimated based on the area/volume of the element type, whereas the characteristic length scale for IELSM elements is estimated following the scheme in [34].

To account for different loading types, several failure criteria are adopted in this work: (1) the standard Rankine criterion [52], (2) the truncated Rankine criterion [38], (3) the Drucker-Prager criterion [53], and (4) the Ottosen criterion [54] . The yield functions of these failure criteria and the corresponding parameters or variables are given in TABLE 2, where the parameter $K_2$ of the Ottosen criterion is assumed to be 1 [49], and $\zeta$, $\xi$, and $K_1$ are determined from the uniaxial tensile, uniaxial compressive, and shear failure strengths.

TABLE 2 Yield functions of the four failure criteria and their parameters or variables.

| Failure criterion | Yield function | Parameters or variables |
|---|---|---|
| Standard Rankine [52] | $f\left(\bar{\sigma}_{eq}(\bar{\boldsymbol{\sigma}}),r\right)=\langle\bar{\sigma}_1\rangle-r$ | $\langle\bar{\sigma}_i\rangle=\dfrac{\bar{\sigma}_i+\lvert\bar{\sigma}_i\rvert}{2}$ |
| Truncated Rankine [38] | $f\left(\bar{\sigma}_{eq}(\bar{\boldsymbol{\sigma}}),r\right)=\sum_{i=1}^{3}\langle\bar{\sigma}_i\rangle-r$ | |
| Drucker-Prager [53] | $f\left(I_1,J_2,r\right)=\dfrac{f_c-f_t}{2f_c}I_1+\dfrac{f_c+f_t}{2f_c}\sqrt{3J_2}-r$ | $I_1=\text{trace}(\bar{\boldsymbol{\sigma}})$<br>$\bar{\mathbf{s}}=\bar{\boldsymbol{\sigma}}-\text{trace}(\bar{\boldsymbol{\sigma}})\mathbf{I}/3$<br>$J_2=\bar{\mathbf{s}}:\bar{\mathbf{s}}/2$ |
| Ottosen [54] | $f\left(I_1,J_2,\theta,r\right)=\dfrac{\eta\sqrt{J_2}+\zeta I_1+\sqrt{(\eta\sqrt{J_2}+\zeta I_1)^2+4\xi J_2}}{2\rho}-r$ | $\eta=K_1\cos\theta$<br>$\cos3\theta=\dfrac{3\sqrt{3}}{2}\dfrac{J_3}{J_2^{3/2}}$<br>$J_3=\det(\bar{\mathbf{s}})$<br>$\rho=f_c/f_t$ |

It should be noted that defining the damage variable at nodes and then weighting the damage of springs and additional volumetric constraints by the damage of their connected nodes, as done in this section, is not the only choice. Alternative schemes include defining damage on regular polygonal elements or on each spring. Each of the three schemes has its own trade-offs. The node-based definition has the highest sensitivity to localization and the best continuity of the damage field, but updating the damage at each node requires traversing its associated elements and springs, resulting in relatively high computational cost. The element-based definition has the lowest cost, but damage is discontinuous between elements, the smeared crack band width is directly tied to the element size, and mesh sensitivity is stronger. The spring-based definition directly corresponds to the failure of physical components and has a clear physical picture, but when coupled with FEM, damage information must be transferred the discrete spring representation and the continuum finite element representation, resulting in the worst compatibility. A quantitative comparison among the node-based, element-based, and spring-based definitions is left for future work.

## 5. Numerical Validation

This study has proposed various regular polygonal IELSM elements. In this section, the TRI, SQ, and HEX IELSM elements are selected for numerical validation. In view of the eigenvalue analysis in Section 3.4, which showed that H9 has a zero-energy mode and H12 exhibits stiffness anisotropy in shear directions, the HEX element adopts the H15 spring arrangement. This section employs a coupled solution scheme of IELSM elements and FEM elements. Unless otherwise specified, these are collectively referred to as IELSM-TRI/CST, IELSM-SQ/CST and IELSM-HEX/CST. In addition, a standard FEM-CST comparison scheme is set up based on the IELSM-TRI/CST mesh, i.e., the TRI elements of IELSM are converted into CST elements. The mesh generation method for the IELSM-FEM coupling in this section is as follows: the problem domain is first divided into an IELSM element region and an FEM element region; regularly arranged polygonal elements are placed in the IELSM element region, and unstructured CST elements are then filled in the other regions using the built-in meshing module of ABAQUS. The IELSM region is the region where damage is expected to occur.

This work adopts quasi-static incremental analysis with a fixed increment step. The damage solution uses an implicit equilibrium-explicit damage extrapolation (IMPL-EX) scheme: within each load step, the linear tangent stiffness system under a frozen damage field is solved once; damage does not participate in the iteration but is recomputed only once at the end of the load step, after which the stiffness is updated for the next step. This scheme naturally bypasses the ill-conditioned tangent problem caused by damage softening, but at the cost of losing unconditional stability: since damage is explicitly extrapolated, the load increment must be sufficiently small to ensure that the damage history evolves along the true equilibrium path. It should be noted that the loading step sizes for all examples in this section have been verified, and the solved damage fields have converged with sufficient accuracy.

Computational efficiency is an important aspect in evaluating the practical value of IELSM. This paper aims to simulate small-deformation linear elastic fracture problems: element assembly is performed only once at the beginning of the computation to obtain the global stiffness matrix;

thereafter, the damage evolution at local nodes directly reduces the global stiffness matrix, thus avoiding repeated assembly of the global stiffness matrix and improving computational efficiency. Therefore, compared with the cumulative time for solving the static equilibrium equations at each load step, the stiffness matrix assembly time is negligible. For this reason, Sections 5.3 – 5.5 record the CPU solution times of each scheme to compare the computational efficiency of different schemes.

### 5.1 Bending test

In this example, plane stress conditions are assumed, and a square plate with side length $l$ = 10 m and unit thickness is considered. The Young's modulus $E$ and Poisson's ratio $v$ are specified as 10 GPa and 0.49, respectively. The geometry and boundary conditions are shown in FIGURE 9(a): on the bottom edge ($y$ = 0), $v(x, 0) = 0$; at the bottom-left corner (0, 0), $u(0, 0) = 0$; and on the top edge ($y = l$), a linearly distributed normal stress in the x$x$-direction is applied, $\sigma_y(x, l) = \sigma(1 - 2x/l)$.

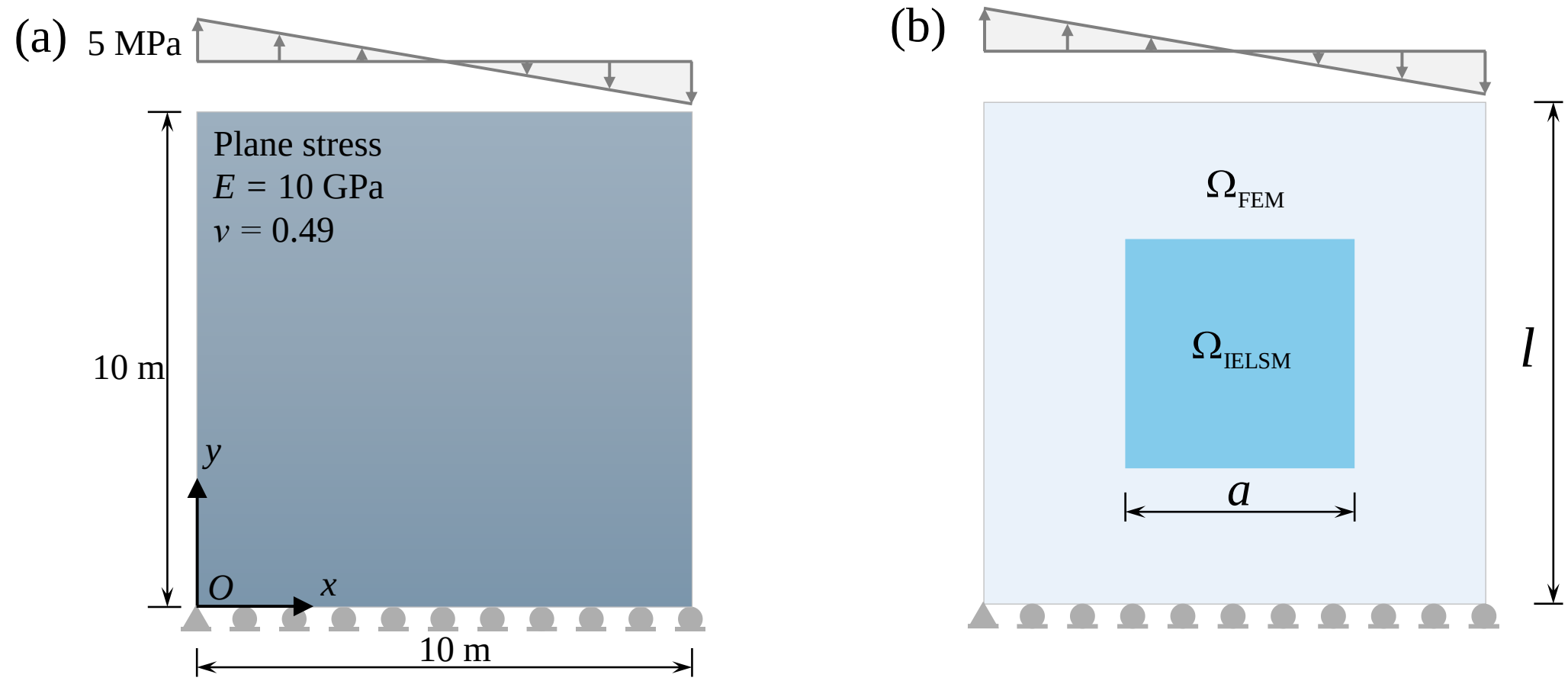

FIGURE 9 Bending test: (a) geometry and boundary conditions; (b) partition of the IELSM region $\Omega_{IELSM}$ and the FEM region $\Omega_{FEM}$.

The above loading induces pure bending deformation in the plate. Solving the stress function with the boundary conditions yields the analytical stress field:

$$\begin{cases} \sigma_{xx} = 0, \\ \sigma_{yy} = \sigma\left(1 - \dfrac{2x}{l}\right), \\ \tau_{xy} = 0. \end{cases} \tag{43}$$

From Hooke's law, the analytical strain field is:

$$
\begin{cases}
\varepsilon_{xx} = -\dfrac{\nu\sigma}{E}\left(1-\dfrac{2x}{l}\right), \\
\varepsilon_{yy} = \dfrac{\sigma}{E}\left(1-\dfrac{2x}{l}\right), \\
\varepsilon_{xy} = 0.
\end{cases}
\tag{44}
$$

Integrating the geometric equations and substituting the boundary conditions yields the analytical displacement field:

$$
\begin{cases}
u = -\dfrac{\sigma}{E}\dfrac{y^2}{l} - \dfrac{\nu\sigma}{E}\left(x-\dfrac{x^2}{l}\right), \\
v = \dfrac{\sigma}{E}\left(1-\dfrac{2x}{l}\right)y.
\end{cases}
\tag{45}
$$

To verify the coupling capability of IELSM and FEM, a coupled IELSM-FEM mesh is designed as shown in FIGURE 9(b), in which the square IELSM element region is located at the center of the problem domain. The parameter $a$ denotes the side length of the IELSM region, and the parameter $l$ denotes the side length of the problem domain. The mesh refinement details are shown in FIGURE 10, where the distributions of FEM elements and IELSM elements are plotted in gray and blue, respectively, and the FEM element size is consistent with the IELSM element size. The meshes in FIGURE 10 consider two variables: the element size $L$ and the IELSM region size $a$, which are used for convergence analysis and for analyzing the influence of the IELSM region proportion, respectively. The relative error indicator is defined as:

$$
\begin{cases}
E_{r,displacement} = \dfrac{\left\|\boldsymbol{u}-\boldsymbol{u}_a\right\|_2}{\left\|\boldsymbol{u}_a\right\|_2}, \\
E_{r,stress} = \dfrac{\left\|\boldsymbol{\sigma}-\boldsymbol{\sigma}_a\right\|_2}{\left\|\boldsymbol{\sigma}_a\right\|_2}, \\
E_{r,strain} = \dfrac{\left\|\boldsymbol{\varepsilon}-\boldsymbol{\varepsilon}_a\right\|_2}{\left\|\boldsymbol{\varepsilon}_a\right\|_2},
\end{cases}
\tag{46}
$$

where $\|\cdot\|_2$ denotes the $L^2$ norm of a vector, $\boldsymbol{u} = [u_1, v_1, u_2, v_2, \dots, u_n, v_n]$ is the vector containing all nodal displacement components, $\boldsymbol{\sigma} = [\sigma_{xx,1}, \sigma_{yy,1}, \sigma_{xy,1}, \sigma_{xx,2}, \sigma_{yy,2}, \sigma_{xy,2}, \dots, \sigma_{xx,n}, \sigma_{yy,n}, \sigma_{xy,n}]$ is the vector containing all nodal stress components, $\boldsymbol{\varepsilon} = [\varepsilon_{xx,1}, \varepsilon_{yy,1}, \varepsilon_{xy,1}, \varepsilon_{xx,2}, \varepsilon_{yy,2}, \varepsilon_{xy,2}, \dots, \varepsilon_{xx,n}, \varepsilon_{yy,n}, \varepsilon_{xy,n}]$ is the vector containing all nodal strain components, and $\boldsymbol{u}_a$, $\boldsymbol{\sigma}_a$ and $\boldsymbol{\varepsilon}_a$ are the vectors containing the analytical solutions obtained from Eqs. (45), (43), and (44), respectively.

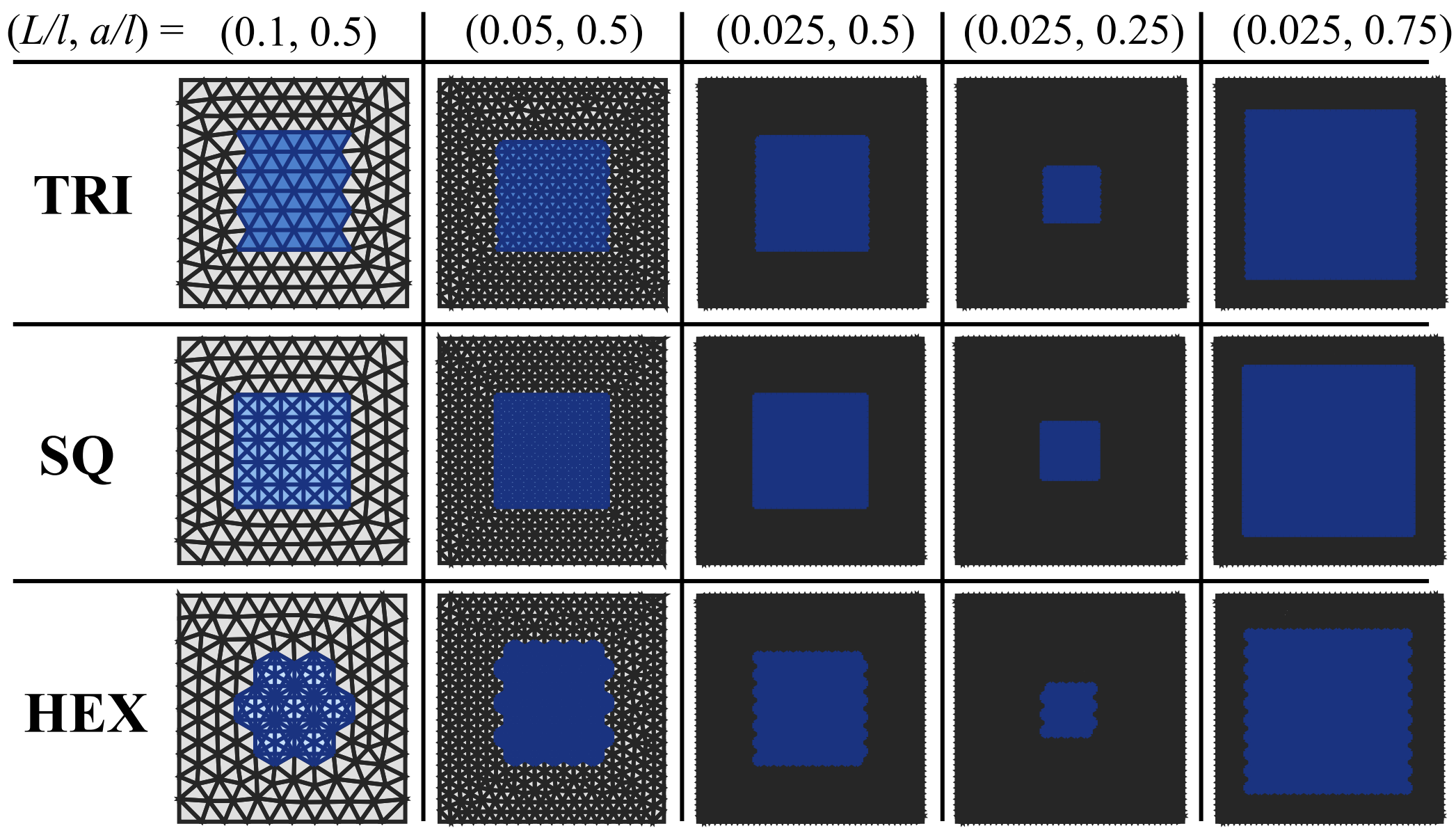

FIGURE 10 Bending test: IELSM-FEM coupled mesh.

FIGURE 11(a) – (c) show the convergence of the displacement, stress, and strain fields relative to the analytical solutions for different mesh sizes at a fixed $a/l$ = 0.5. Here, FEM-Q4, FEM-CST-structured, and FEM-CST-unstructured correspond to results computed using a structured mesh of Q4 elements, a structured mesh of CST elements, and an unstructured mesh of CST elements, respectively; IELSM-SQ corresponds to results computed using full-domain IELSM square elements. As shown in FIGURE 11(a)–(c), at the same mesh density, the displacement field accuracy obtained by the three IELSM-FEM coupling schemes is close, and the errors are close to those of IELSM-SQ and smaller than those of full-domain FEM. The convergence rate $R$ is defined as the slope of the error–element size curve in a log-log plot. The displacement field error convergence rates of the three IELSM-FEM coupling schemes ($R$ = 2.05–2.19) are all higher than those of IELSM-SQ ($R$ = 2.01) and full-domain FEM ($R$ = 1.88 – 2.01); the stress and strain field error convergence rates (stress field $R$ = 1.32 – 1.38, strain field $R$ = 1.29 – 1.35) are close to those of IELSM-SQ (stress field $R$ = 1.38, strain field $R$ = 1.39) and lower than those of full-domain FEM (stress field $R$ = 1.39 – 1.56, strain field $R$ = 1.45 – 1.58). The lower convergence rates of the stress and strain fields in the IELSM coupling schemes compared with full-domain FEM are an expected result of IELSM using linear displacement field recovery (with accuracy comparable to CST), consistent with its inherent linear displacement field assumption. FIGURE 12 shows the $L_2$ error norms for different $a/l$ at a fixed $L/l$ = 0.025, where $a/l$ = 0 means that the IELSM region degenerates to empty and the entire domain is the FEM region, whose result is equivalent to that of FEM-CST-unstructured; $a/l$ = 1 corresponds to the full-domain IELSM result, for which only the IELSM-SQ solution is given because the other

element shapes cannot completely tessellate the square plate and can only provide approximate boundaries. The results in FIGURE 12 indicate that the relative error of the displacement field computed using IELSM-FEM coupling is independent of the IELSM region proportion; the relative error results of the stress and strain fields show that the relative errors of the IELSM-TRI/CST and IELSM-SQ/CST schemes are independent of $a/l$, while the relative error of the IELSM-HEX/CST scheme decreases as $a/l$ decreases. FIGURE 13 further shows contour plots of the absolute error of $\gamma_{xy}$ computed by different discretization schemes at a fixed $L/l = 0.025$. In FIGURE 13(a) – (d), the error distribution in the IELSM region is close to zero, whereas the error distribution in the FEM region is non-uniform; regarding the error distribution at the coupling boundaries, (a) IELSM-TRI/CST performs best, with the coupling boundary inheriting the small internal error; (b) IELSM-SQ/CST exhibits obvious error bands at the upper and lower coupling boundaries; (c) IELSM-HEX/CST, in addition to error bands at the upper and lower boundaries, shows a high-frequency sawtooth distribution of errors at the left and right boundaries; and (d) IELSM-SQ has errors concentrated on the boundaries where boundary conditions are applied. Overall, the displacement field accuracy computed by IELSM is not degraded by coupling with FEM; the stress and strain field errors are more pronounced at the coupling boundaries, which is related to the post-processing scheme.

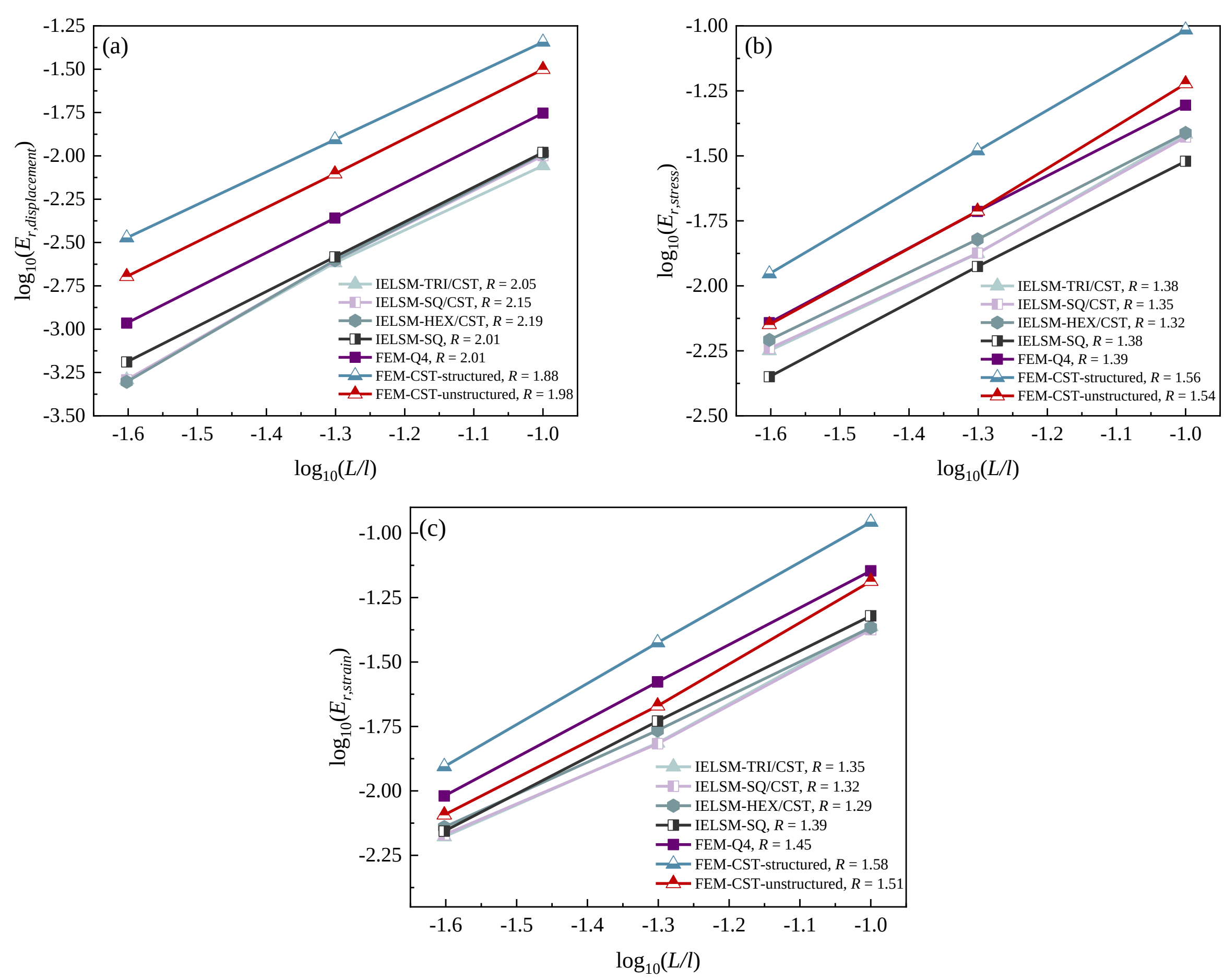

FIGURE 11 Bending test: Convergence analyses. ($a/l$ = 0.5)

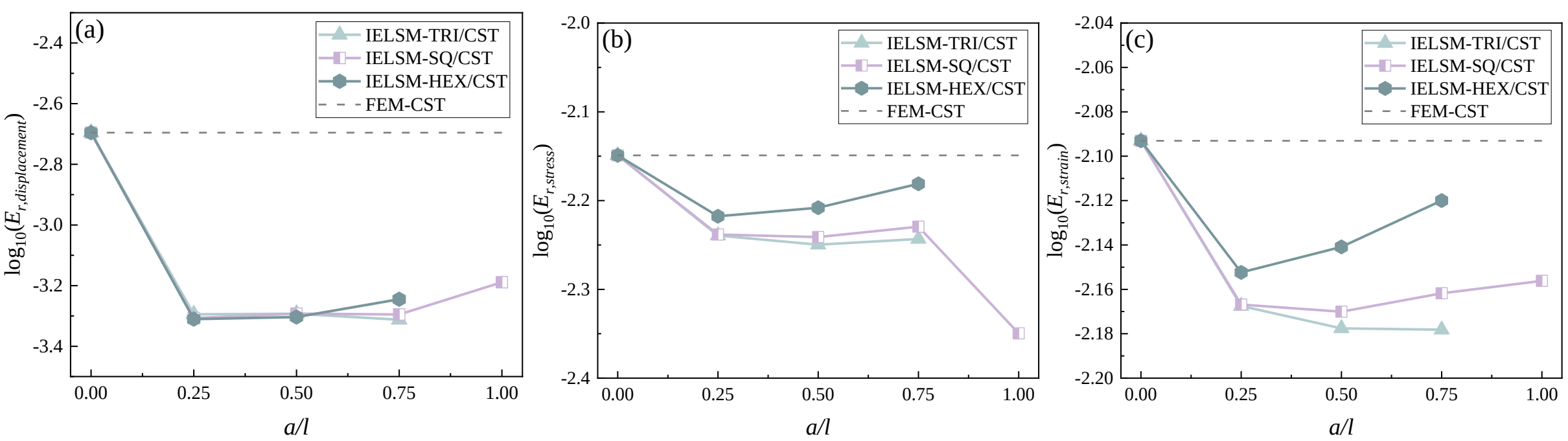

FIGURE 12 Bending tests: Relationship between the relative error of IELSM-FEM coupling computation and $a/l$. ($L/l$ = 0.025)

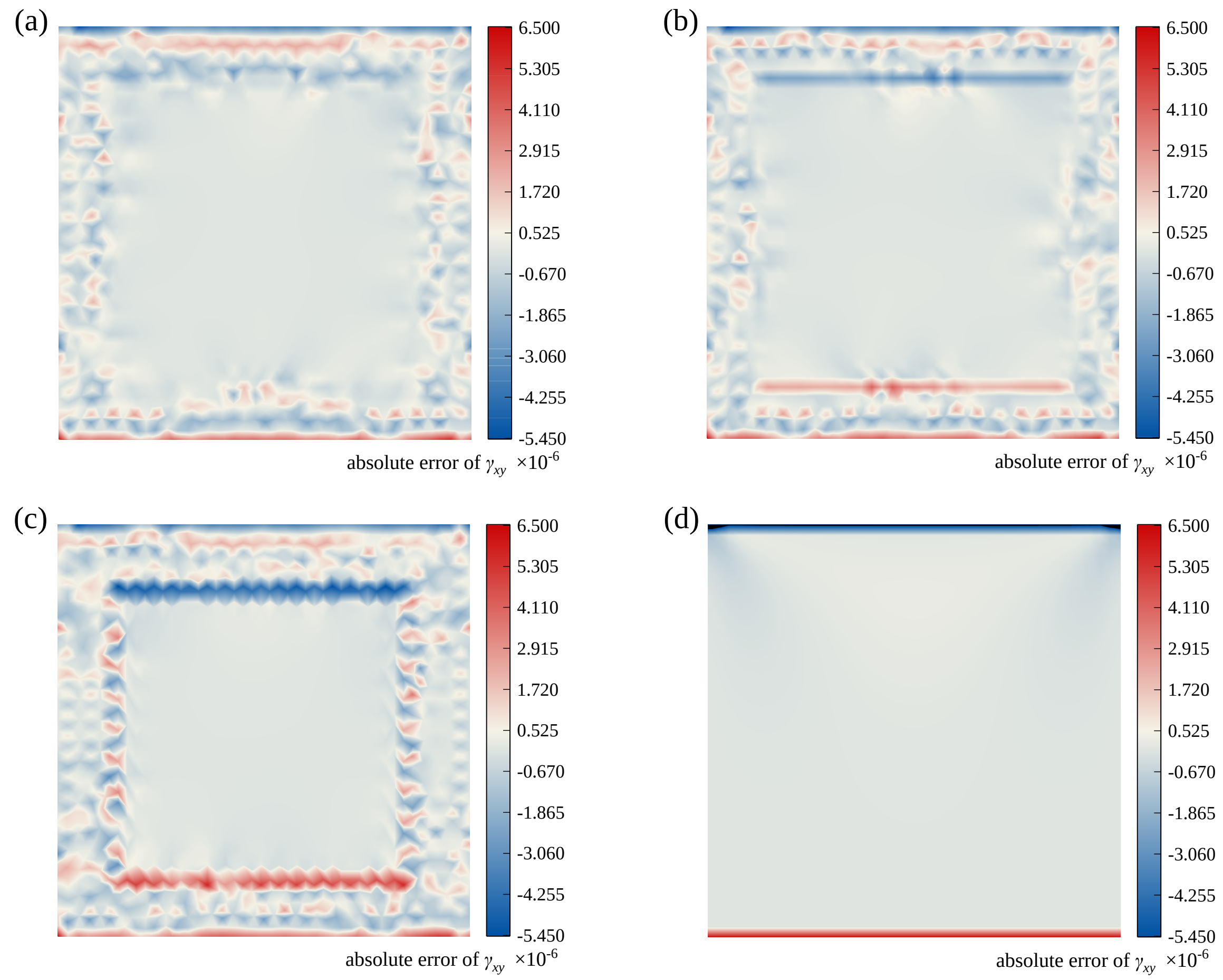

FIGURE 13 Bending test: Absolute error distributions of $\gamma_{xy}$ computed by different discretization schemes: (a) IELSM-TRI/CST ($L/l$ = 0.025, $a/l$ = 0.75); (b) IELSM-SQ/CST ($L/l$ = 0.025, $a/l$ = 0.75); (c) IELSM-HEX/CST ($L/l$ = 0.025, $a/l$ = 0.75); (d) IELSM-SQ ($L/l$ = 0.025, $a/l$ = 1).

## 5.2 Symmetric three-point bending test

This example simulates a three-point bending test of a beam under symmetric boundary conditions, with experimental data taken from Grégoire et al. [55]. Grégoire's tests are frequently used to calibrate the characteristic length of numerical methods coupled with local isotropic damage models and to verify whether the computational model can reproduce the size effect [34,48,56]. FIGURE 14(a) shows the geometry and boundary conditions of this example. The span-to-depth ratio of the beam is kept constant at 2.5, and the beam depth $D_n$ follows the specifications in the experiments, i.e., $D_0$ = 50mm, $D_1$ = 100mm, $D_2$ = 200mm, $D_3$ = 400mm. This example also considers three notch-to-depth ratios: 0.2 (fifth-notched), 0.5 (half-notched), and 0.0 (unnotched). The notch width is 2 mm, consistent with the dimensions recorded in the experiments [55]. All test beams have a thickness of 50 mm, and plane stress conditions are assumed. The material properties used in this example are taken from the literature [48], as shown in TABLE 3. FIGURE 14(b) – (d) show the partition of the IELSM region and the FEM region for the models with the three notch sizes (excluding the meshes used for calibrating the element characteristic length). Only the Standard Rankine failure criterion is used in this example. All tests are subjected to vertical displacement control at specified points, and the corresponding load reactions are recorded. The loading protocol is uniformly set as follows: a step size of 0.0001 mm before a total applied displacement of 0.15 mm, and a step size of 0.0005 mm after a total displacement of 0.15 mm.

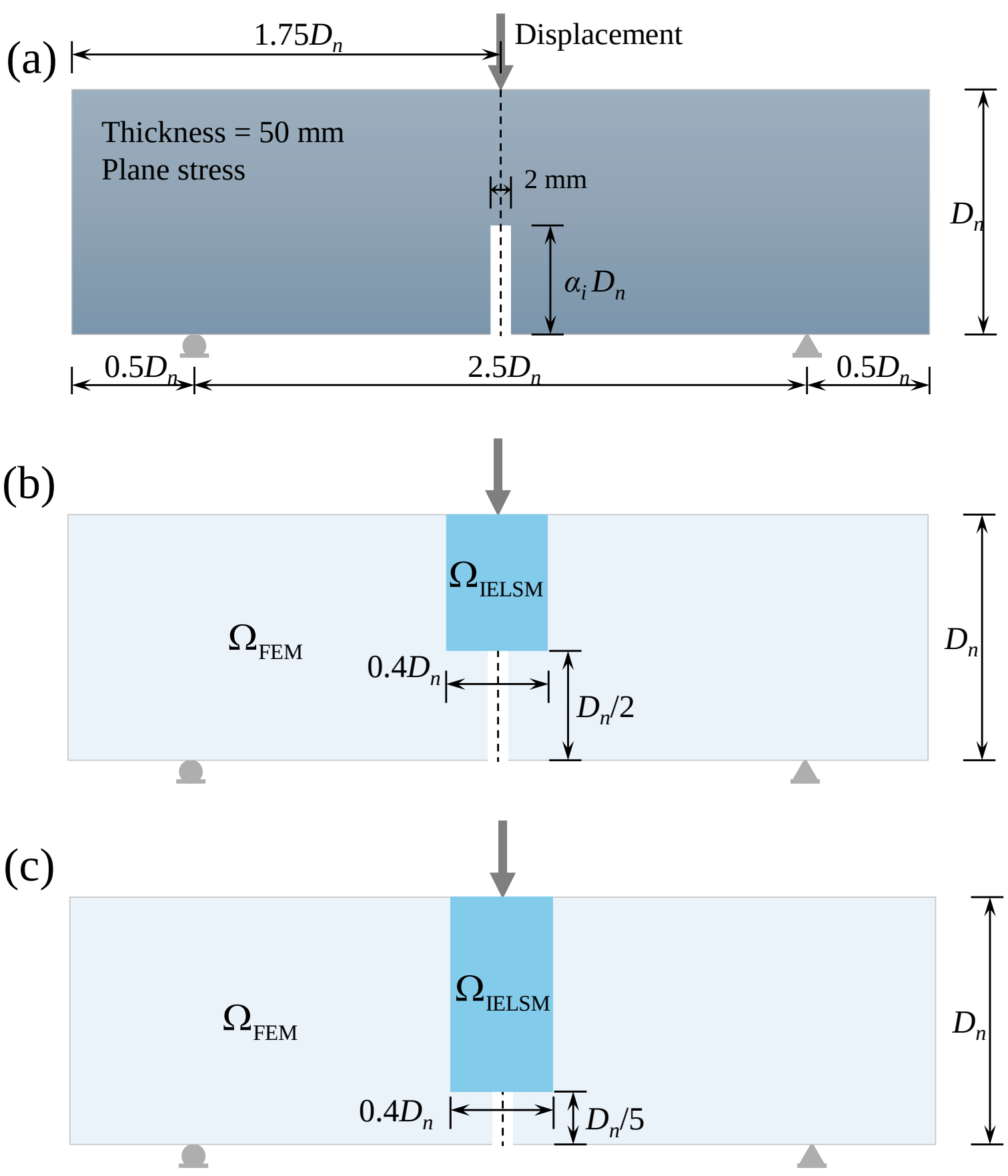

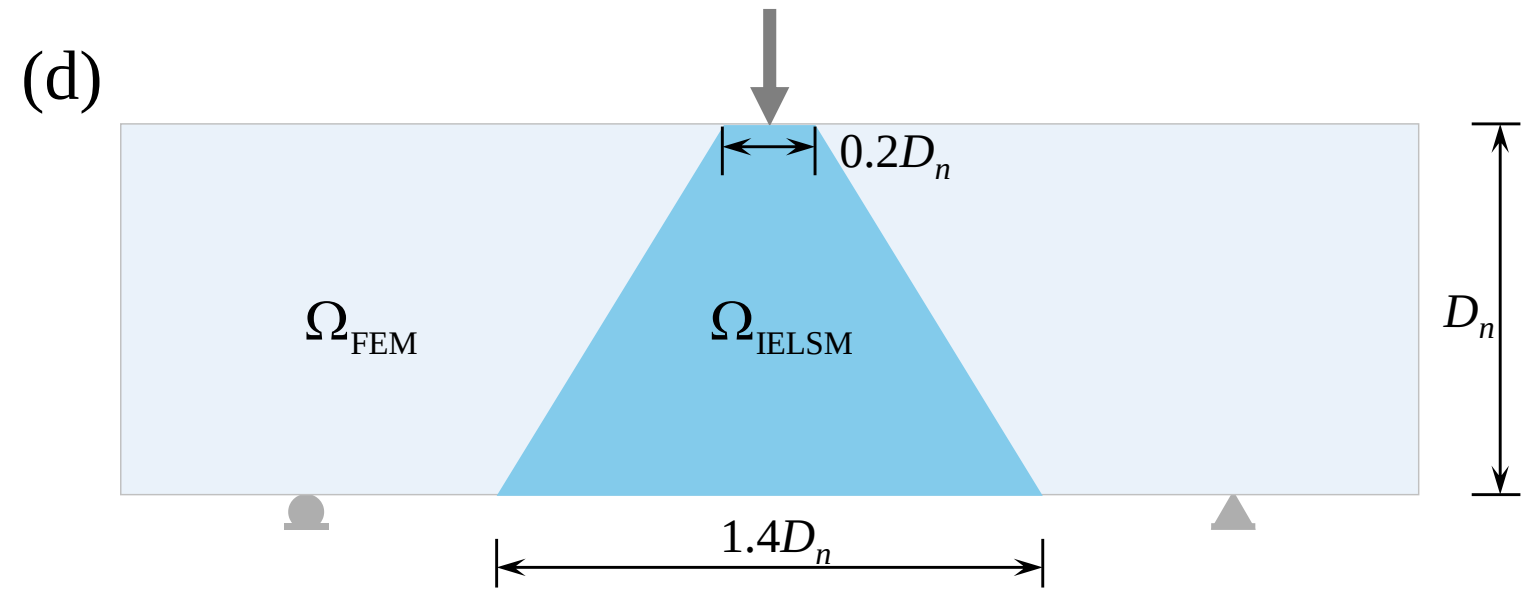

FIGURE 14 Grégoire's beam: (a) geometry and boundary conditions; (b) – (d) partition of the IELSM region $\Omega_{IELSM}$ and the FEM region $\Omega_{FEM}$ for the notched models in (b) and (c) and the unnotched model in (d).

TABLE 3 Grégoire tests: Material properties.

| Symbol | Value |
|---|---|
| $E$ | 37 GPa |
| $v$ | 0.2 |
| $f_t$ | 3.5 MPa |
| $G_f$ | 90 J/m$^2$ |

FIGURE 15, FIGURE 16, and FIGURE 17 show the characteristic length calibration details for the IELSM TRI, SQ, and HEX elements, respectively. In FIGURE 15, the IELSM-TRI/CST mesh is also used to calibrate the characteristic length of the FEM CST element. The characteristic length calibration in this paper follows the following idea: based on the $D_0$ = 50mm (fifth-notched) specimen, the prefabricated notch is changed by a very small size (one element length) and the computation is repeated. The area difference between the two force-displacement curves is the energy required to release one element length of crack propagation in that specimen, from which the apparent fracture energy can be determined. This process is repeated for different mesh sizes, and the characteristic length that makes the apparent fracture energy computed with different mesh sizes closest to the preset fracture energy is selected as the calibrated value, so as to minimize the sensitivity of the characteristic length to the mesh size. The results are given in TABLE 4. In addition, TABLE 5 shows the peak forces computed using various element types with different sizes when fine-tuning the characteristic length. The variation in peak force caused by fine-tuning the characteristic length is about 0.5% for TRI, about 1.2% – 1.3% for SQ, about 5.6% for HEX, and about 1.5% for CST; the variation in peak force caused by mesh size changes is no more than 0.6% for all schemes. Except for the HEX element, which is relatively sensitive to fine-tuning of the characteristic length (about 5.6%), the peak force variation of the other

combinations is within 1.5%; all schemes are insensitive to mesh size changes (≤0.6%), indicating that the calibrated characteristic length parameters are stable.

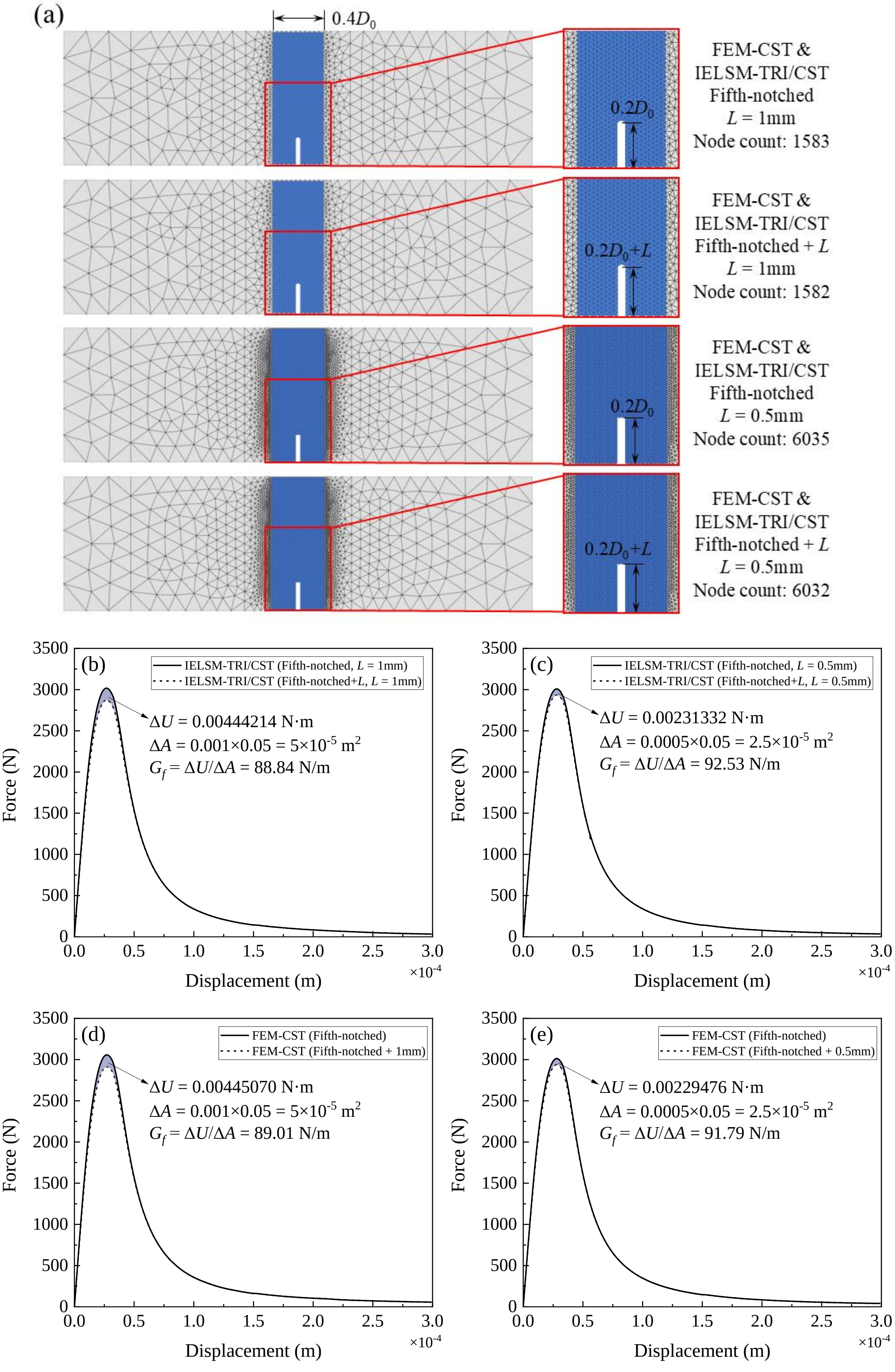

FIGURE 15 Grégoire's beam: Characteristic length calibration for the IELSM TRI element.

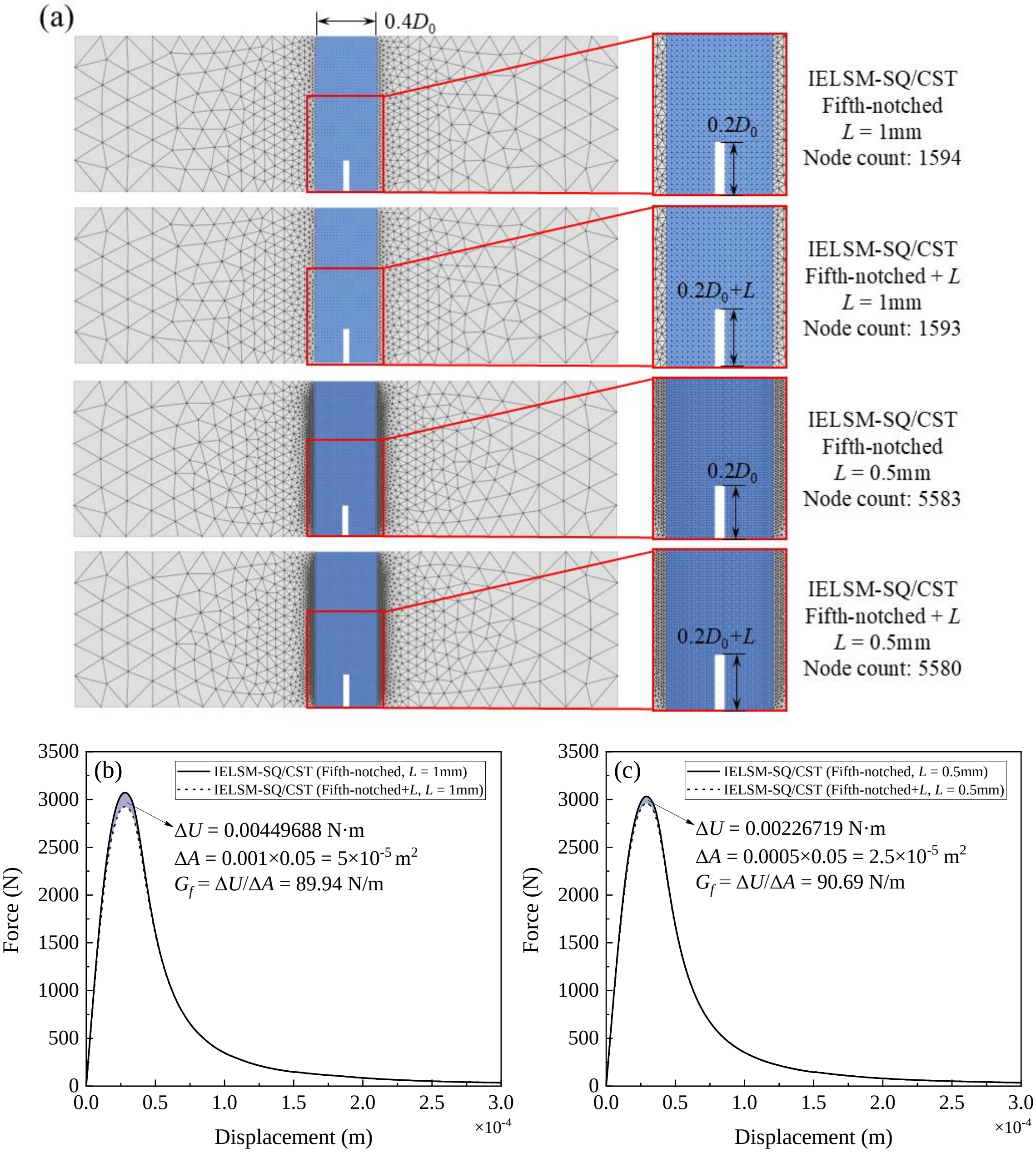

FIGURE 16 Grégoire's beam: Characteristic length calibration for the IELSM SQ element.

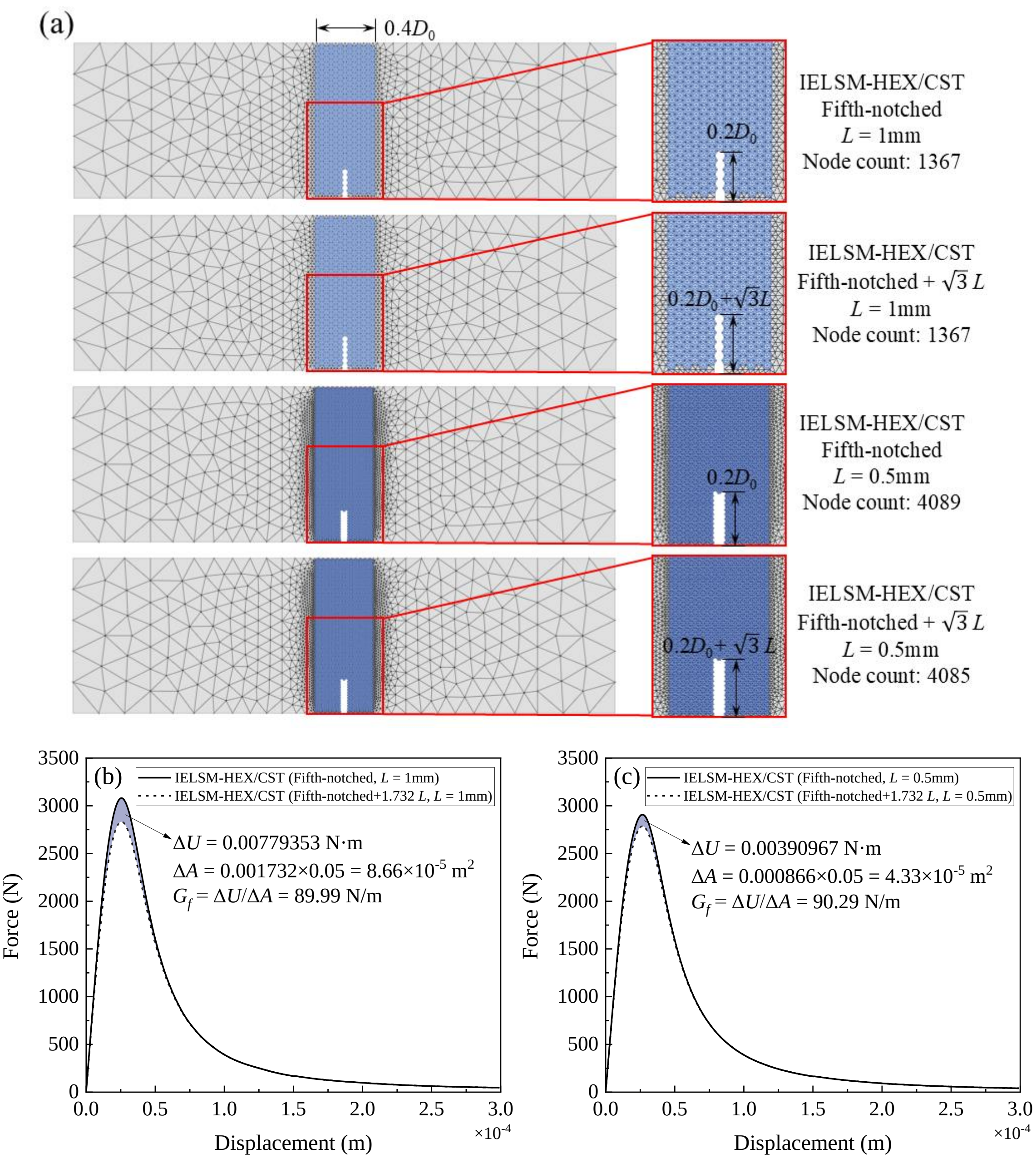

FIGURE 17 Grégoire's beam: Characteristic length calibration for the IELSM HEX element.

TABLE 4 Characteristic lengths used in this study.

| Characteristic length | Value |
|---|---|
| $l^*_{\mathrm{TRI}}$ | 3.0 |
| $l^*_{\mathrm{SQ}}$ | 3.4 |
| $l^*_{\mathrm{HEX}}$ | 8.4 |
| $l^*_{\mathrm{CST}}$ | 3.0 |

TABLE 5 Grégoire's beam: Peak forces for different mesh sizes when fine-tuning the characteristic length.

| Element | $l^*$ | $L$ = 1 mm | $L$ = 0.5 mm | Variation due to $L$ |
|---|---|---|---|---|
| TRI | 3.0 | 3017 | 3008 | 0.3 % |
| | 3.1 | 3001 | 2992 | 0.3 % |
| Variation due to $l^*$ | | 0.5 % | 0.5 % | |
| SQ | 3.4 | 3073 | 3032 | 1.3 % |
| | 3.5 | 3056 | 3018 | 1.2 % |
| Variation due to $l^*$ | | 0.6 % | 0.5 % | |
| HEX | 8.4 | 3080 | 2908 | 5.6 % |
| | 8.5 | 3073 | 2901 | 5.6 % |
| Variation due to $l^*$ | | 0.2 % | 0.2 % | |
| CST | 3.0 | 3057 | 3012 | 1.5 % |
| | 3.1 | 3040 | 2996 | 1.4 % |
| Variation due to $l^*$ | | 0.6 % | 0.5 % | |

FIGURE 18 shows a comparison of the load-crack mouth opening displacement (CMOD) curves among IELSM, experimental [55], and mixed-FEM (referring to the mixed finite element results of Barbat et al. [48], the same below) results for beams with four depths and three notch sizes. For notched beams, the CMOD is taken as the relative displacement of the two nodes at the prefabricated notch at the beam bottom; for unnotched beams, the CMOD is measured at two points on the beam bottom separated by a distance $D_n$, each located at $0.5D_n$ from the beam center. Considering the computation time, the IELSM element size for the unnotched beam with depth 400 mm is $L$ = 2 mm, and $L$ = 1 mm in all other cases. For both notched and unnotched beams, the results of IELSM-TRI/CST, IELSM-SQ/CST, and FEM-CST are very close, and their predicted peak loads are slightly higher than those of IELSM-HEX/CST, with the difference becoming larger for beams with greater depth. The peak loads predicted by IELSM-HEX/CST and mixed FEM are relatively close and agree well with the experimental results. For notched beams, the curves predicted by IELSM are lower than those of mixed FEM in the descending branch; for unnotched beams, the predicted curves are higher than those of mixed FEM in the descending branch. For the unnotched beams with $D_2$ = 200 mm and $D_3$ = 400 mm, the descending branches predicted by IELSM are approximately straight lines. Under displacement-controlled loading, this indicates that after the crack enters steady-state propagation, the structural response evolves along an elastic unloading path: as the crack propagates, the overall stiffness of the beam continuously decreases, and the load decays approximately linearly with crack length; this stage is dominated by elastic strain energy release. Unlike the snap-back phenomenon reported for large-size beams in similar studies using explicit time integration schemes [34], the quasi-static displacement-controlled IMPL-EX solution in this paper does not exhibit snap-back, because the displacement-controlled scheme itself does not trace the snap-back branch and the solution evolves along the elastic unloading line. The two phenomena originate from the solution control method rather than from differences in the model's prediction of fracture energy release. The overall comparison shown in FIGURE 18 indicates that, for both notched and unnotched beams, the IELSM results agree well with the experiments. The present model can capture the structural response of geometrically similar beams of different sizes.

FIGURE 19 shows the damage contours at the end of loading for test cases with different notch sizes at $D_2$ = 200 mm. Elements in which the damage values of all nodes are greater than 0.99 are regarded as failed (this setting is used hereinafter), and they are directly removed to create voids (only in visualization), thereby producing crack propagation paths. For all notched beams, the computed crack path is a straight line from the notch to the top midpoint, and the resulting damage band is centered on the crack path with uniform width; the results predicted by all schemes are consistent. For unnotched beams, damage first accumulates and spreads over a wider region at the beam bottom, then localizes and propagates toward the top. However, for IELSM-TRI/CST, IELSM-SQ/CST, and FEM-CST, local damage accumulation occurs in the middle stage of crack propagation (the post-peak softening stage), finally producing a gourd-shaped damage pattern (with two regions of damage spreading, at the beam bottom and in the post-peak softening stage), whereas the damage evolution of IELSM-HEX/CST is consistent with the description provided in the literature [34,48] using isotropic damage models. To explain this phenomenon, the damage contours at various stages are given in FIGURE 20 for the unnotched beam with $D_2$ =200 mm. Compared with the other schemes, the damage predicted by IELSM-HEX/CST is more concentrated at the center of the beam bottom before reaching the peak load, so energy is released in a concentrated manner into the central crack, without spurious damage spreading; for the other three schemes, before the peak load, damage is instead concentrated around the middle of the beam bottom, and the pre-peak damage diffusion delays localization, which explains their higher peak loads. In the softening stage, in addition to the main through-crack, these three schemes also produce two secondary cracks on the left and right; the secondary cracks may divert the energy required for main crack propagation. The interaction between the main and secondary cracks slows down main crack propagation, which in turn produces further damage diffusion, forming a gourd-shaped damage contour.

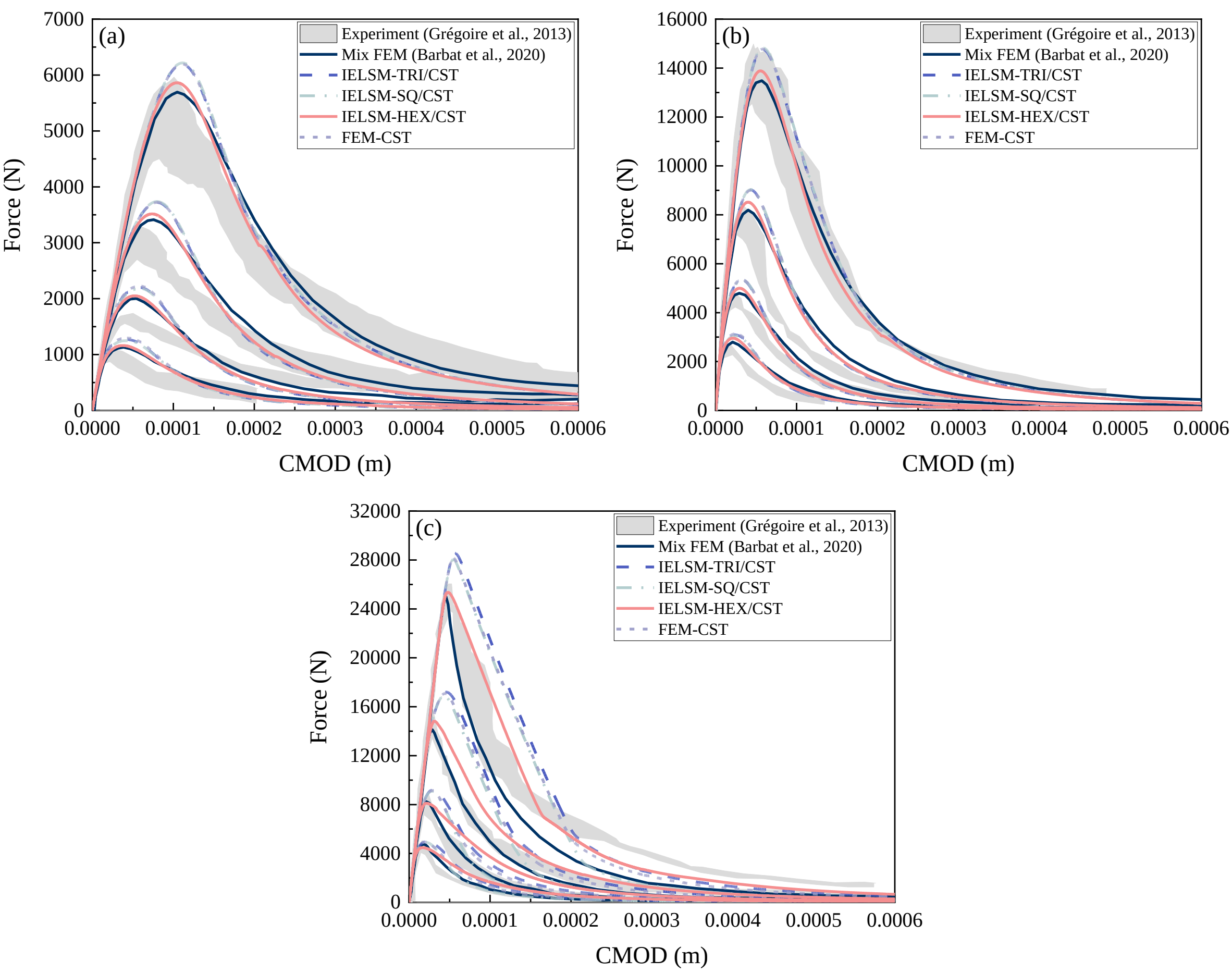

FIGURE 18 Grégoire's beam: Comparison of force-CMOD curves among IELSM, mixed FEM, and experimental results for three-point bending beams with different depths and notch configurations: (a) Half-notched, (b)Fifth-notched, (c)Unnotched.

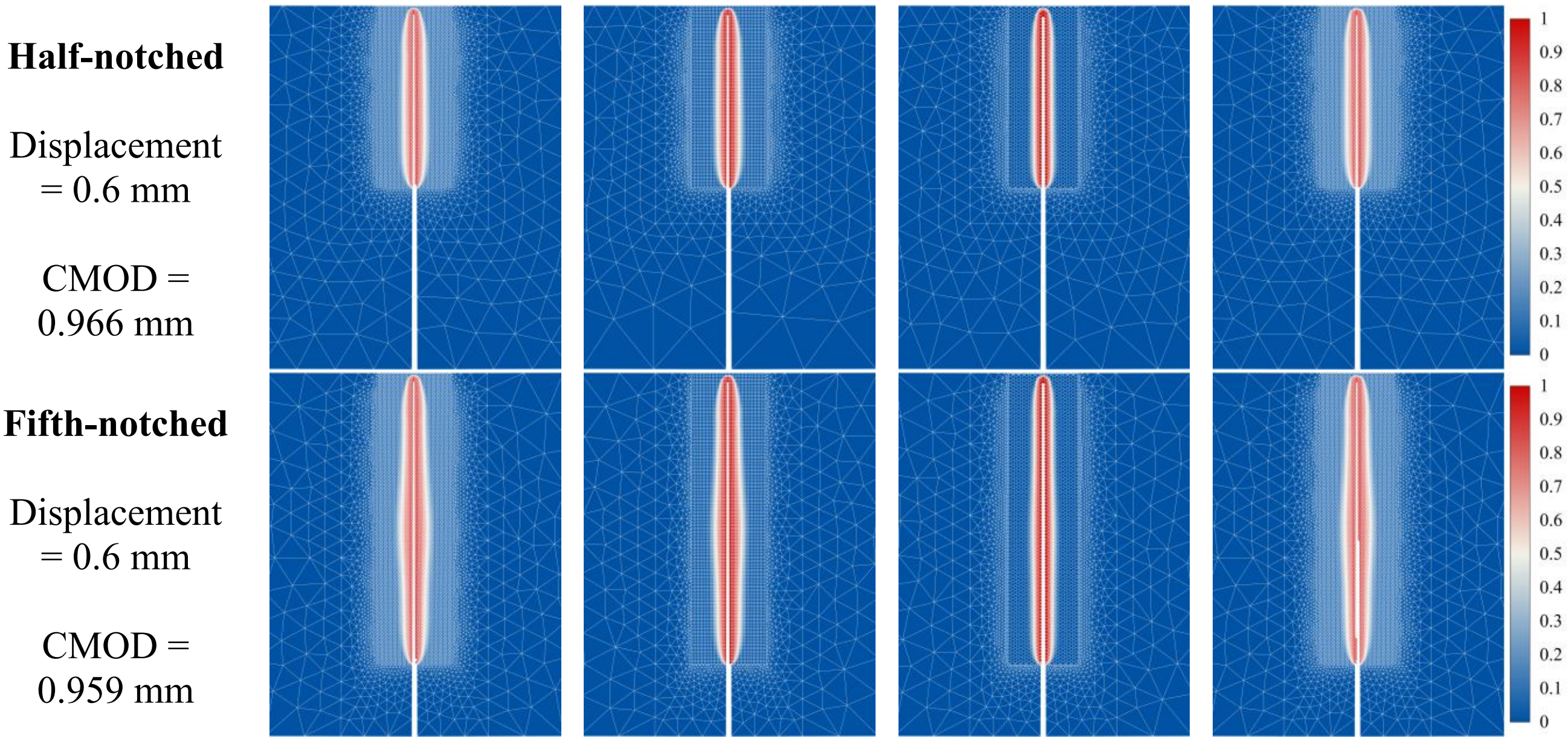

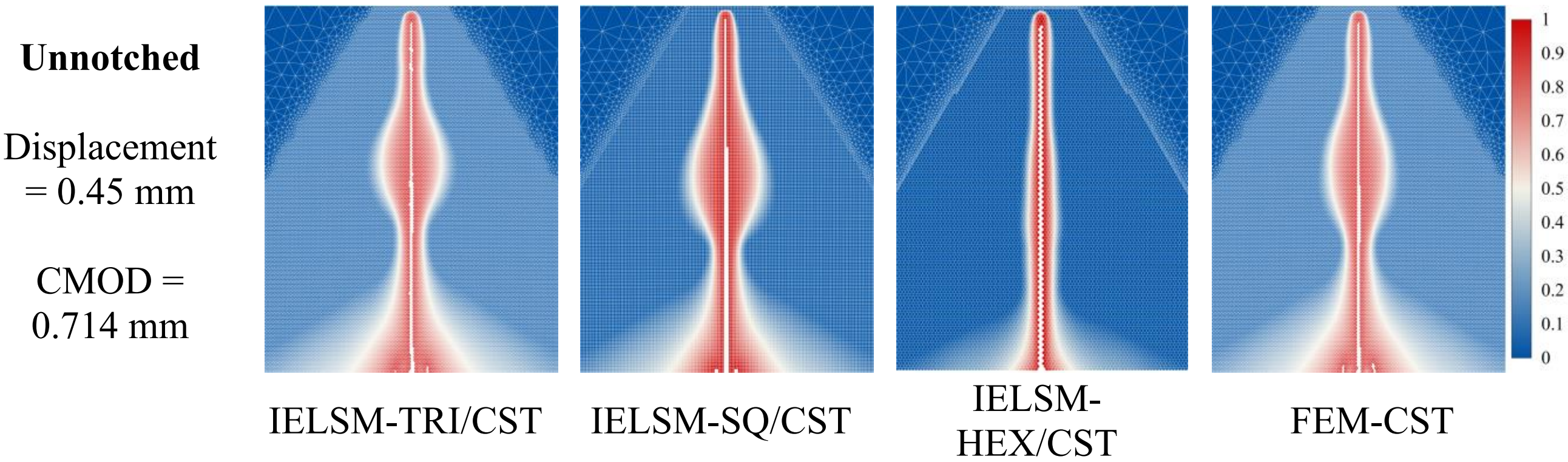

FIGURE 19 Grégoire's beam: Damage contours at the end of loading for test cases with different notch sizes. ($D_2$ = 200 mm)

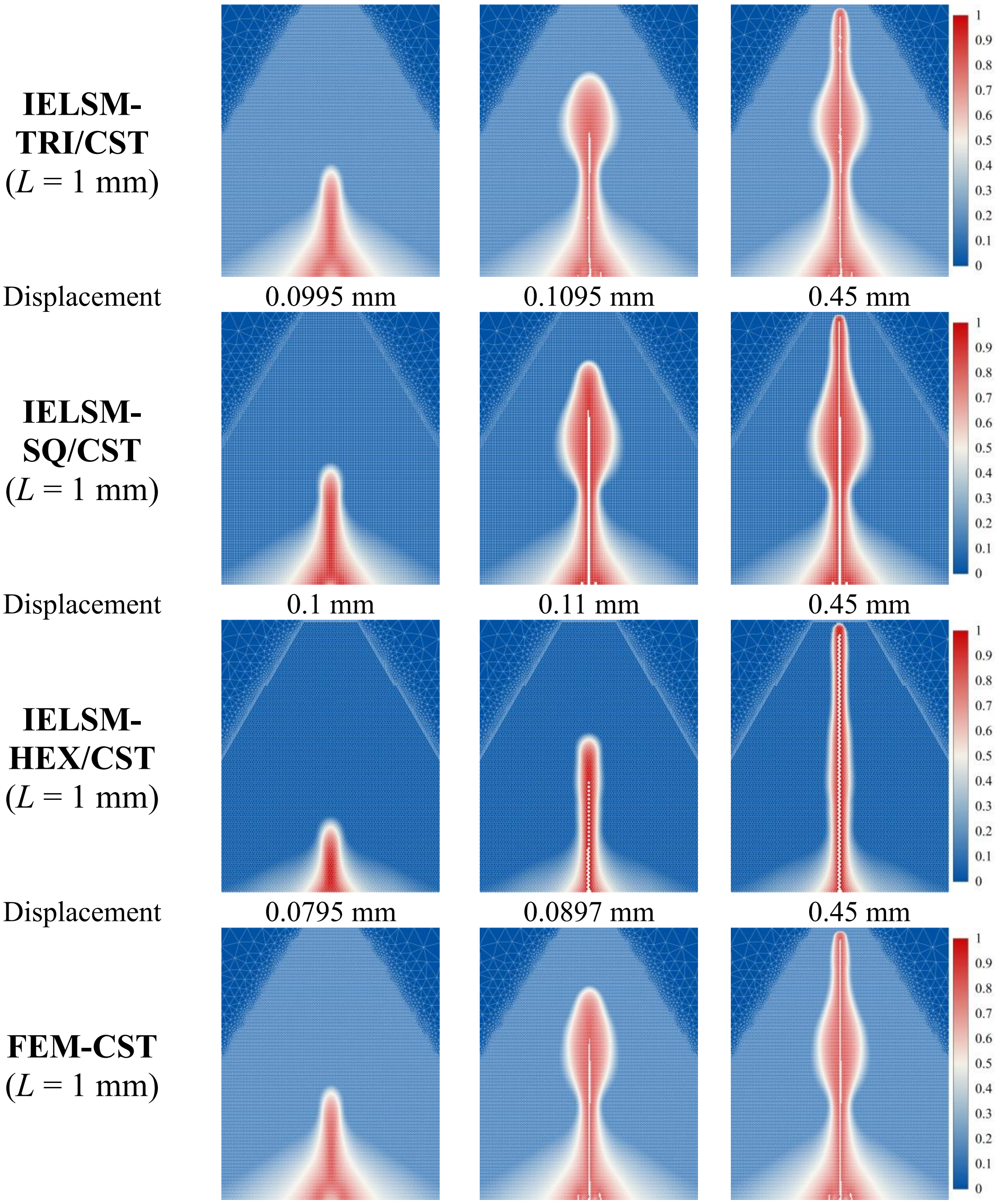

**IELSM-TRI/CST** ($L$ = 1 mm)

Displacement 0.0995 mm 0.1095 mm 0.45 mm

**IELSM-SQ/CST** ($L$ = 1 mm)

Displacement 0.1 mm 0.11 mm 0.45 mm

**IELSM-HEX/CST** ($L$ = 1 mm)

Displacement 0.0795 mm 0.0897 mm 0.45 mm

**FEM-CST** ($L$ = 1 mm)

| Displacement | 0.0995 mm (peak force) | 0.1095 mm (softening) | 0.45 mm (failure) |
|---|---|---|---|

FIGURE 20 Grégoire's beam: Damage contours at various stages predicted by each element type at the same element size $L$ = 1 mm. ($D_2$ = 200 mm, Unnotched)

### 5.3 Asymmetrical three-point bending test

This example simulates a three-point bending test of a beam under asymmetric boundary conditions, with experimental data taken from Gálvez et al. [57]. This benchmark example is widely used to evaluate the ability of numerical models to capture asymmetric crack propagation in quasi-brittle materials. FIGURE 21(a) shows the geometry and boundary conditions of this example. The beam thickness is 50 mm, plane stress conditions are considered, and the prefabricated notch has a width of 2 mm and a depth of 75 mm. FIGURE 21(b) shows the partition of the IELSM region and the FEM region in this example. The beam is subjected to vertical displacement control at specified points, and the corresponding load reactions are recorded. The loading protocol is uniformly set as 0.0001 mm per step until a total displacement of 0.2 mm is reached. The relative displacement of the two nodes at the prefabricated notch at the beam bottom is recorded as the CMOD. The material properties used in this example are taken from Mixture 1 in the experiments of Gálvez et al. [57] and the specific values are given in TABLE 6.

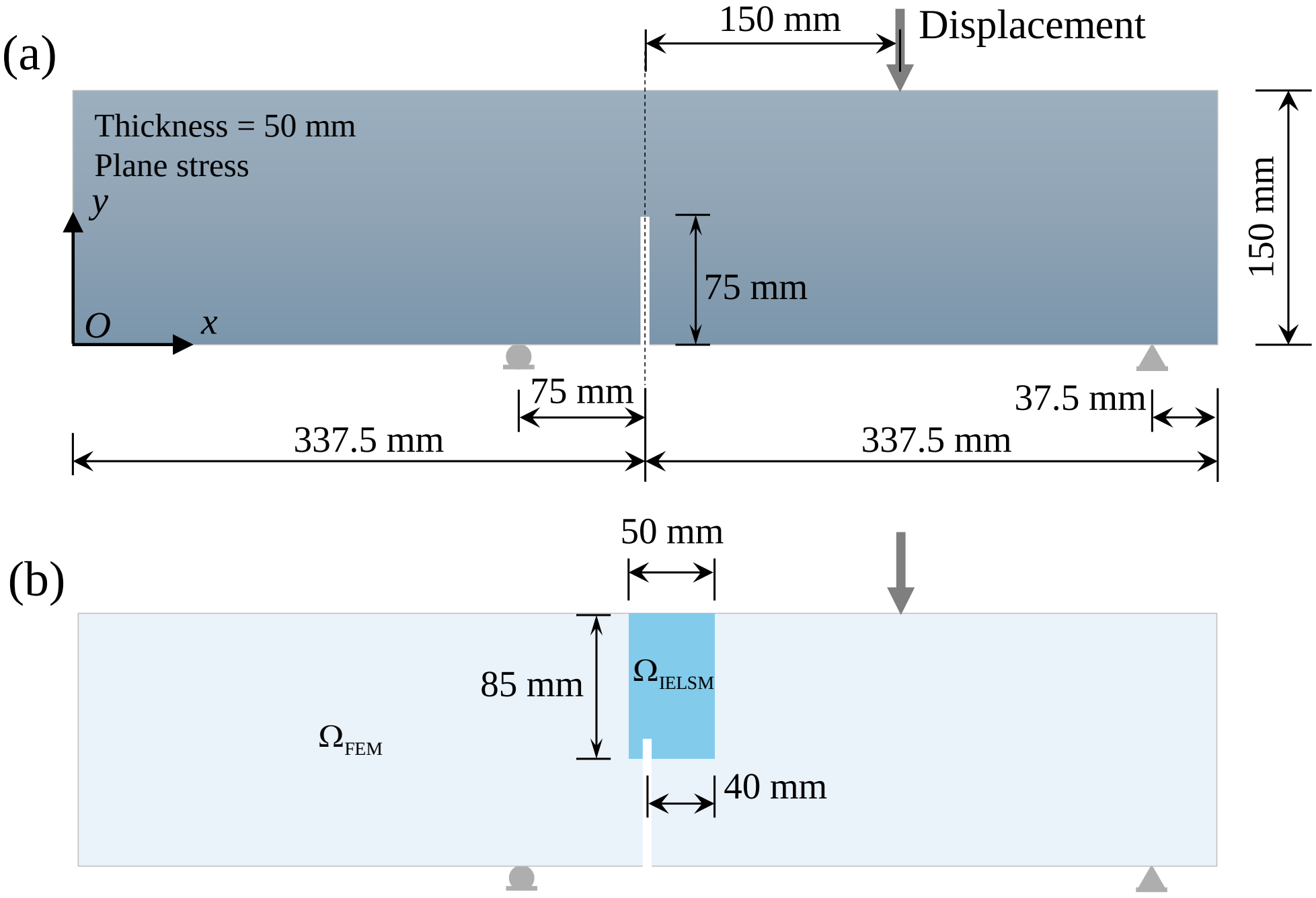

FIGURE 21 Gálvez's beam: (a) geometry and boundary conditions; (b) partition of the IELSM region $\Omega_{IELSM}$ and the FEM region $\Omega_{FEM}$.

TABLE 6 Gálvez's beam: Material properties.

| Symbol | Value |
|---|---|
| $E$ | 38 GPa |
| $v$ | 0.2 |
| $f_t$ | 3.0 MPa |
| $f_c$ | 54 MPa |
| $G_f$ | 69 J/m$^2$ |

As shown in FIGURE 22, the Force-CMOD curves computed using four element types with different damage criteria are compared. FIGURE 22(a) uses the Rankine damage criterion, and the results are larger than the experimental range; FIGURE 22(b) uses the truncated Rankine damage criterion, and the results are smaller than the experimental range; FIGURE 22(c) uses the Drucker-Prager damage criterion, and the results fall within the experimental range. FIGURE 23 presents the Force-CMOD curves computed using different element sizes and element types based on the Drucker-Prager damage criterion. It also compares the IELSM curves with the experimental range [57] and with the results of the isotropic damage model with smooth Rankine-type definition of equivalent strain (IDM-R) reported by Jirásek and Grassl [58], where Structured-Q4 and Unstructured-CST denote the structured Q4 mesh and the unstructured CST mesh, respectively. The results show that the curve predicted by IELSM-HEX/CST agrees with the experimental results and is close to that of Unstructured-CST.. The curves predicted by IELSM-TRI/CST and FEM-CST are highly consistent, but the reaction force in the softening stage decreases faster than the experimental results. The peak force predicted by IELSM-SQ/CST is slightly higher than the experimental range, and its Force-CMOD curve lies entirely above the experimental results. In addition, FIGURE 23 shows that when the element size is refined, the structural responses almost overlap, indicating that the proposed IELSM with an isotropic damage model is insensitive to mesh size.

FIGURE 24 shows the damage contours at the end of loading computed using different element sizes and element types. All schemes provide clear crack paths and damage contours, and a comparison of the results for the same element type with different mesh sizes shows that the damage extent is consistent. The crack contours in FIGURE 24 are further extracted in FIGURE 25 and compared with the experimental results [57] and those of Jirásek and Grassl. [58]. It is observed that the crack contours predicted for the same element type with different element sizes are highly consistent. It should be noted that the meshes in the damage-prone region are structured in the present simulations. Nevertheless, the crack paths predicted by IELSM-TRI/CST, IELSM-HEX/CST, and FEM-CST do not exhibit evident mesh-direction dependence and agree well with the experimental results [57] and with the Unstructured-CST results [58]. In contrast, the crack path predicted by IELSM-SQ/CST is close to that obtained with Structured-Q4 and shows a mesh-direction bias. This can be attributed to the fact that the springs in the square IELSM element are aligned with the element edges and diagonals, which coincide with the structured mesh directions; when the crack tends to propagate obliquely, damage tends to follow these preferential directions. In comparison, the triangular and hexagonal IELSM elements have more diverse spring orientations, which allows them to better accommodate arbitrary crack paths and reduces the mesh-direction bias. FIGURE 26 shows the damage contours at various stages predicted by the four element types on a relatively fine mesh. It can be seen that, compared with the other schemes, IELSM-HEX/CST has the smallest damage diffusion range, the most concentrated damage, and the earliest localization at each stage; under the same displacement load, its crack propagates the longest, and the corresponding load level is the lowest, which is also more consistent with the experimental results.

FIGURE 27 shows the CPU time versus nodes curves for each scheme in this example. The lower the curve, the higher the computational efficiency. The results show that IELSM-TRI/CST and FEM-CST have the highest and relatively similar computational efficiencies, IELSM-HEX/CST has an intermediate computational efficiency, and IELSM-SQ/CST has relatively poor computational efficiency.

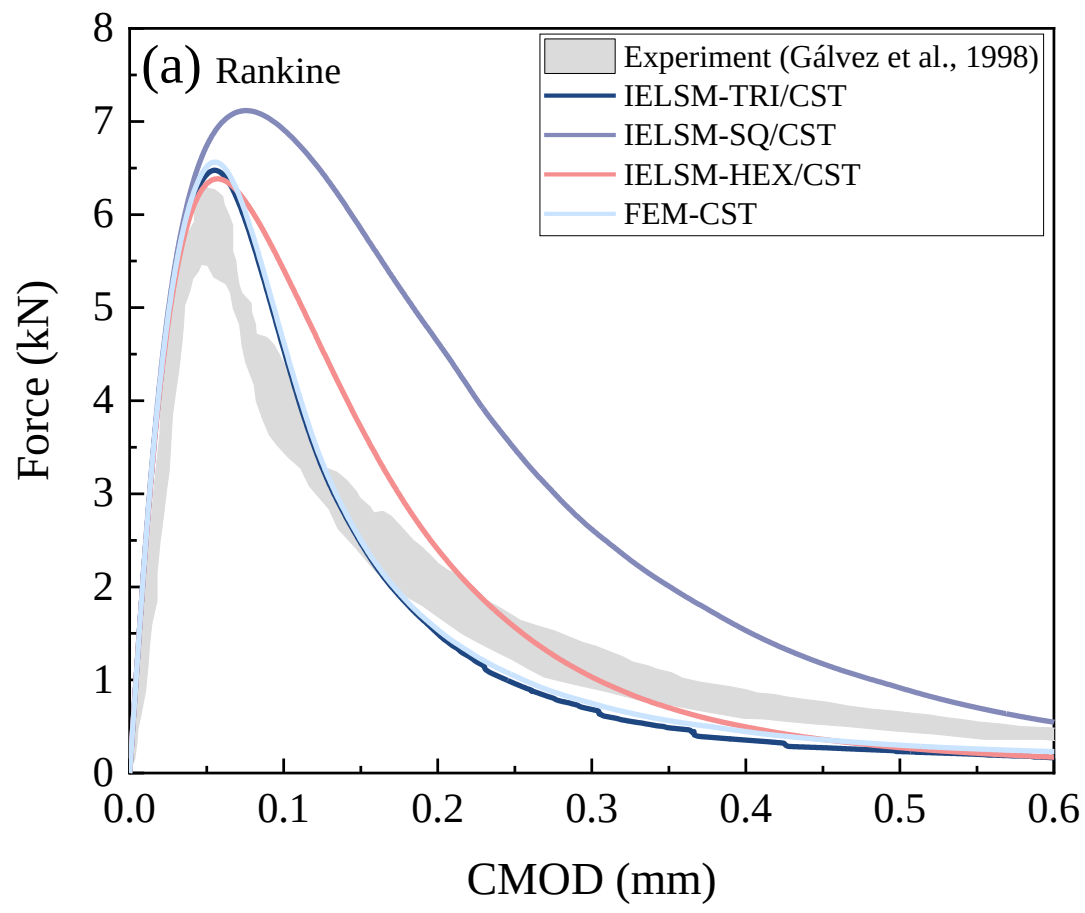

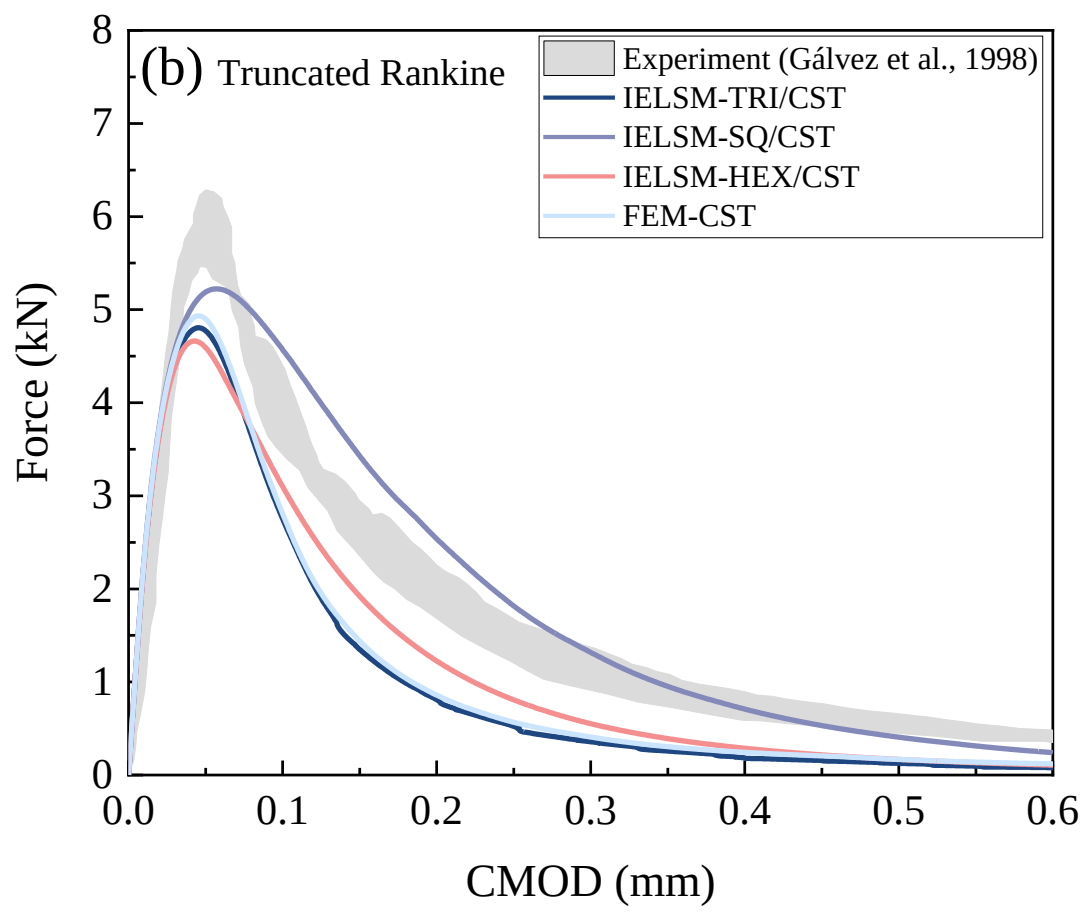

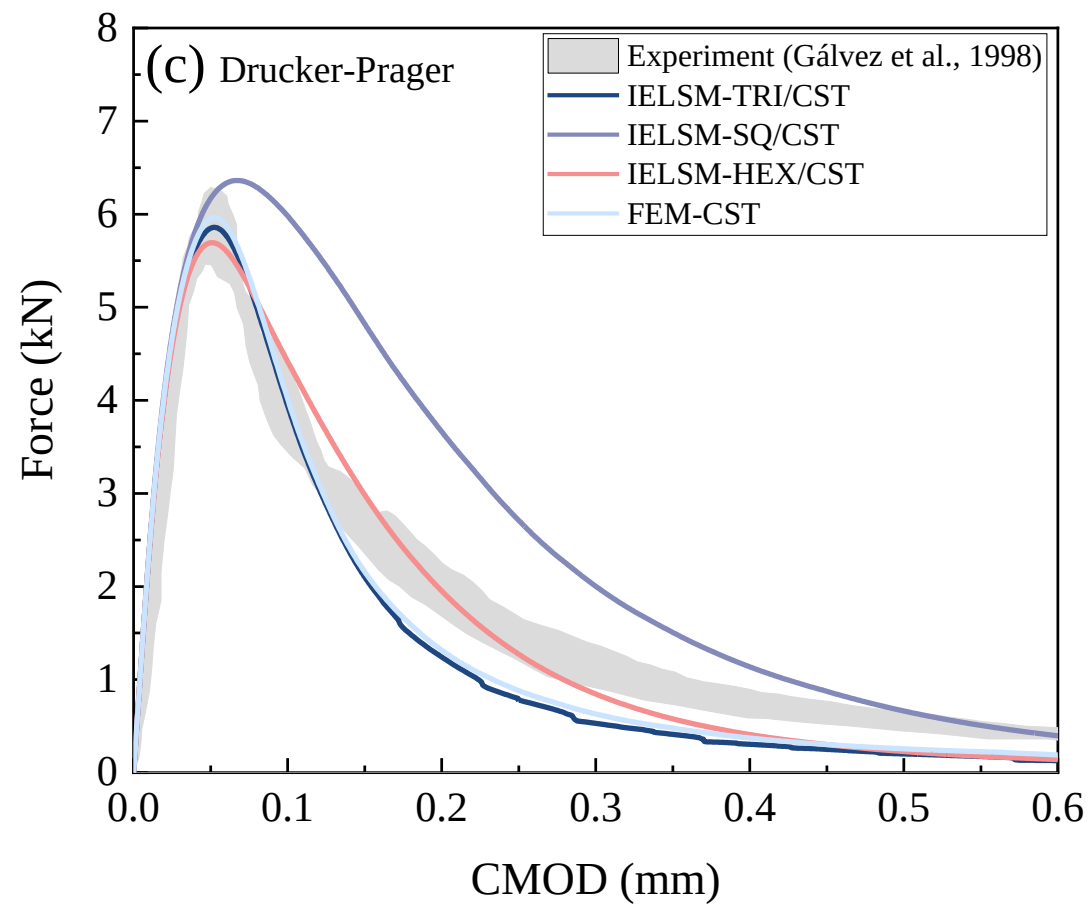

FIGURE 22 Gálvez's beam: Force-CMOD curves computed using four element types under different damage criteria: (a) Rankine, (b) Truncated Rankine, (c) Drucker-Prager.

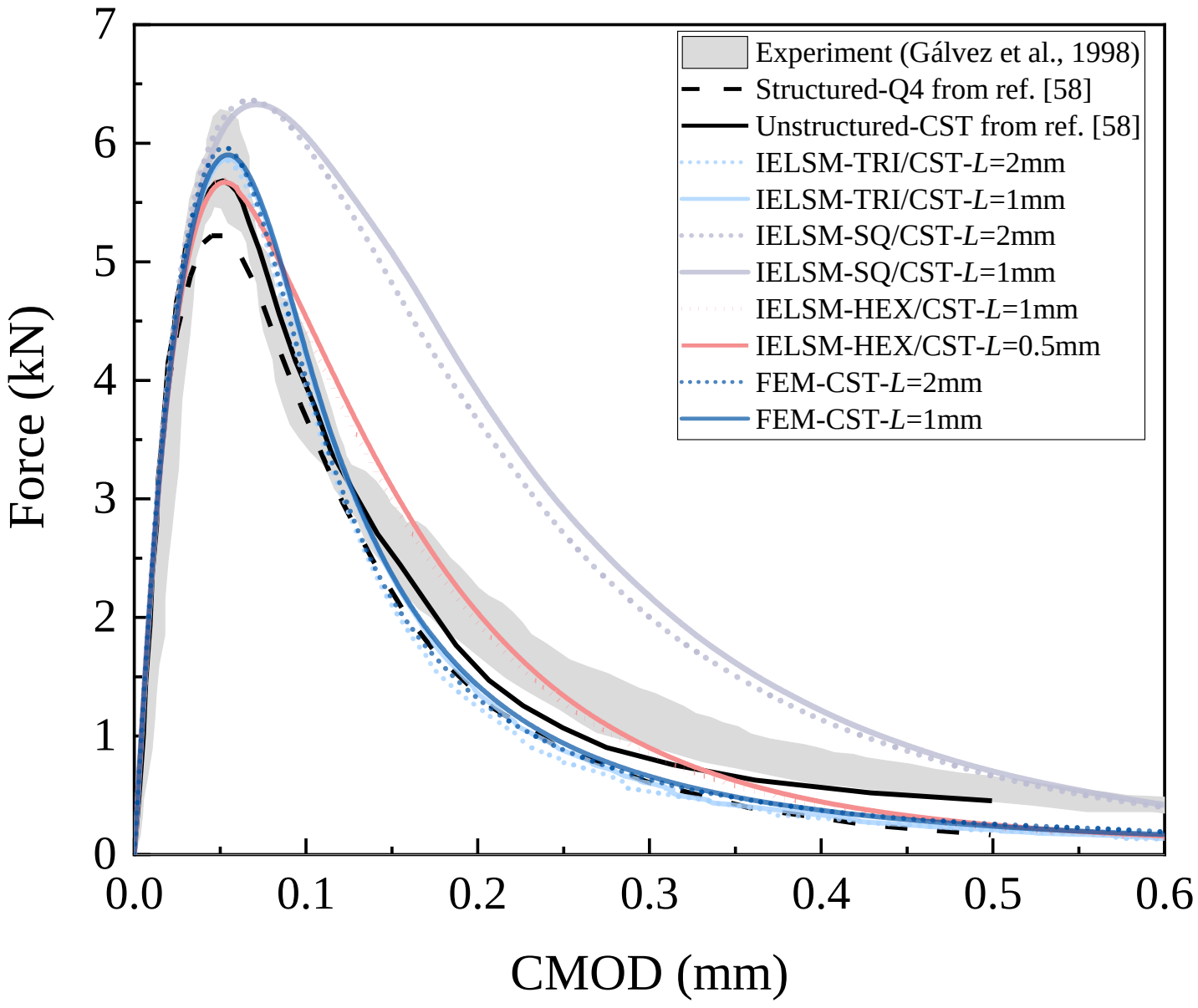

FIGURE 23 Gálvez's beam: Force-CMOD curves computed using different element sizes and element types.

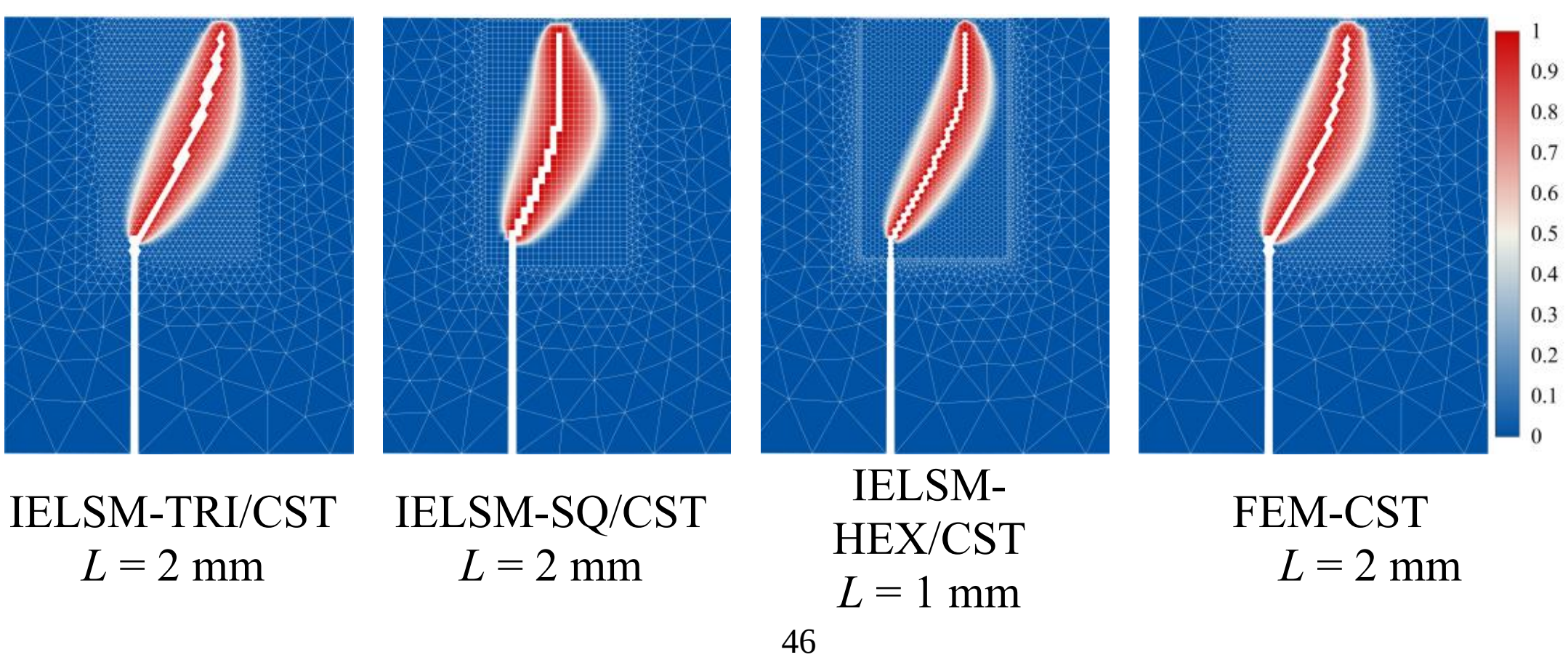

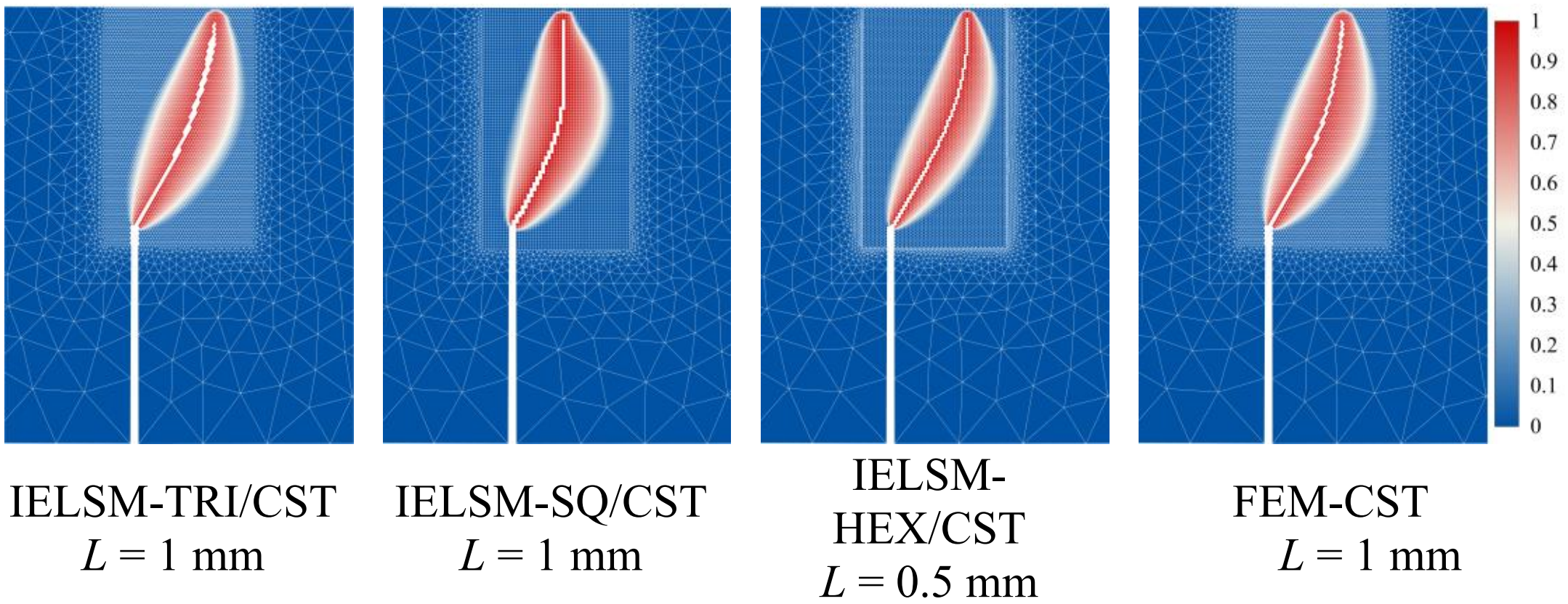

FIGURE 24 Gálvez's beam: Damage contours at the end of loading computed using different element sizes and element types. (Displacement = 0.2 mm)

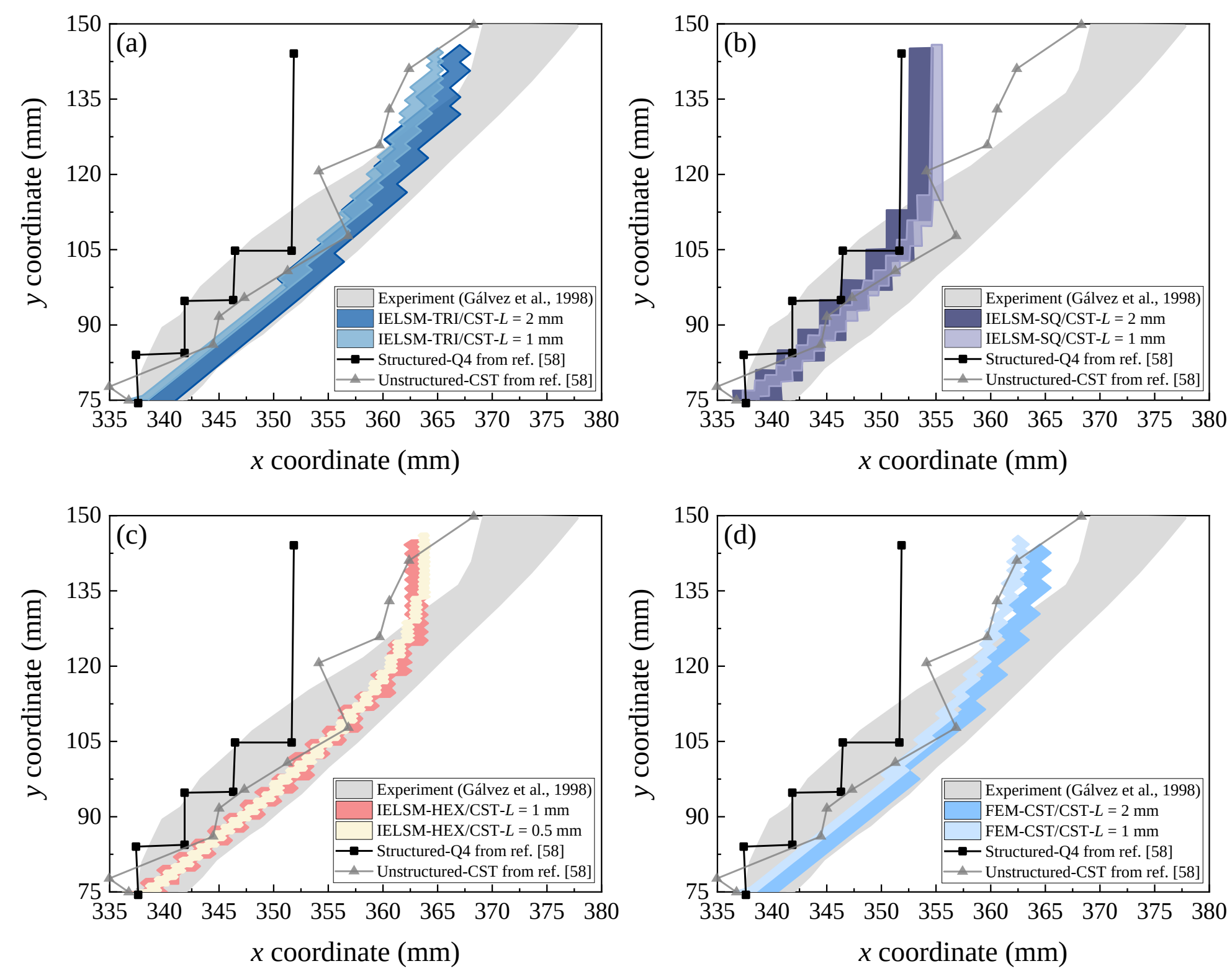

FIGURE 25 Gálvez's beam: Crack contours predicted using different element types and element sizes: (a) IELSM-TRI/CST, (b) IELSM-SQ/CST, (c) IELSM-HEX/CST, (d) FEM-CST.

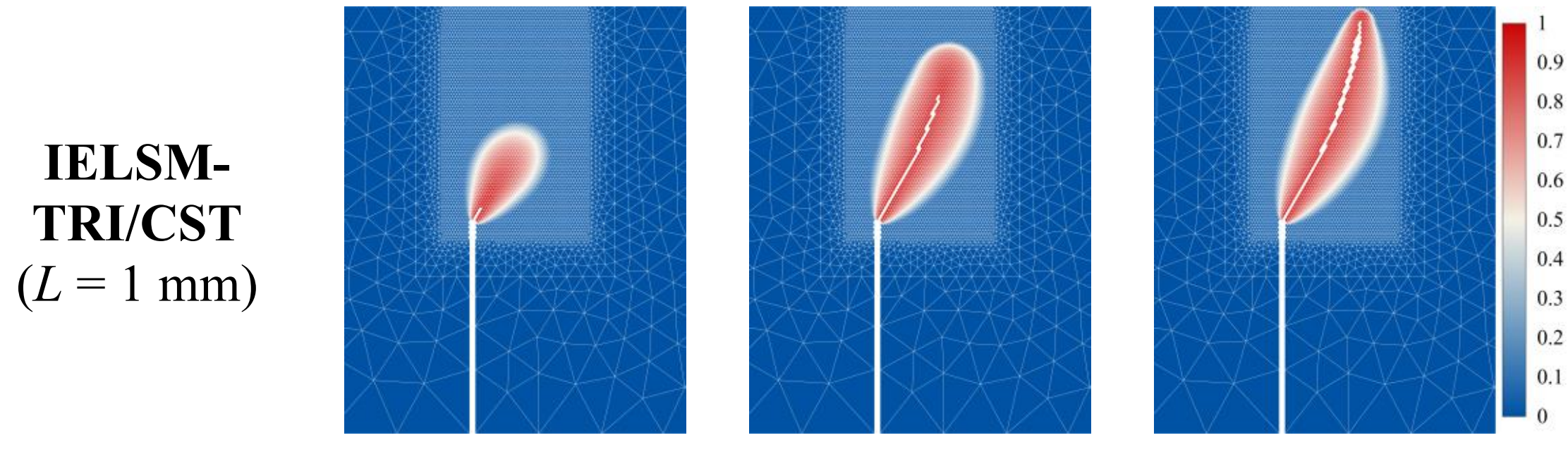

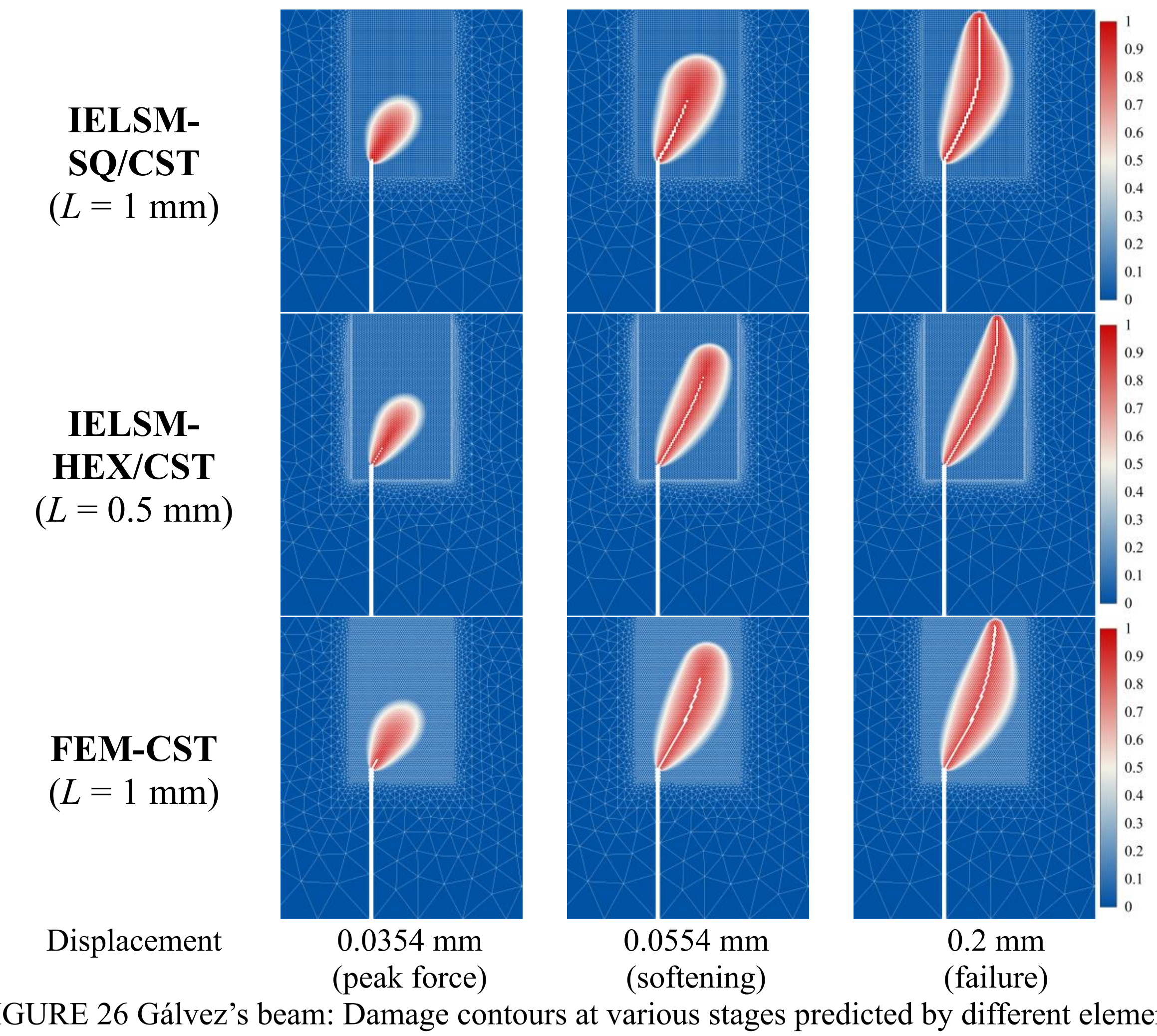

FIGURE 26 Gálvez’s beam: Damage contours at various stages predicted by different element types on a relatively fine mesh.

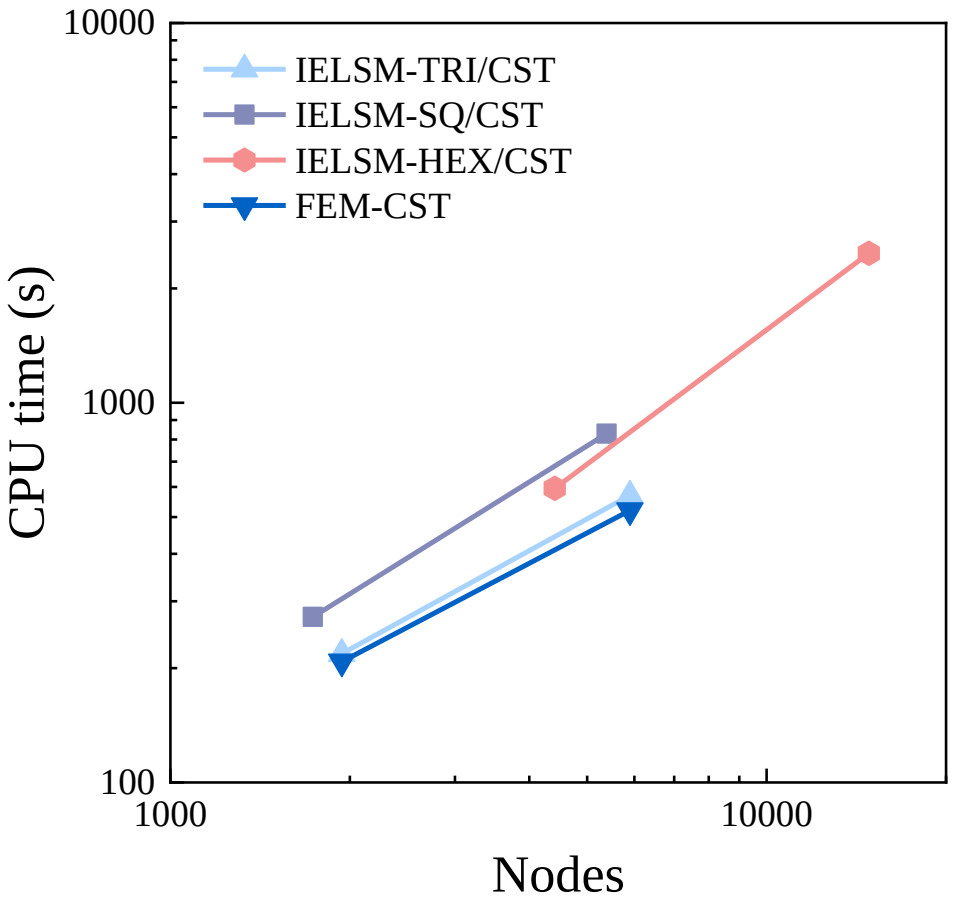

FIGURE 27 Gálvez’s beam: CPU time versus nodes curves for each scheme.

### 5.4 L-shaped panel

This example simulates mixed-mode fracture of an L-shaped concrete panel, with experimental data taken from Winkler et al. [59]. The geometric details and boundary conditions of the L-shaped concrete panel are shown in FIGURE 28(a). The panel thickness is 100 mm, and plane stress conditions are assumed. FIGURE 28(b) shows the partition of the IELSM region and the FEM region in this example. The L-shaped concrete panel is subjected to vertical displacement control at specified points, and the corresponding load reactions are recorded. The loading protocol is 0.0001 mm per step before a total displacement of 0.4 mm, and 0.0005 mm per step after 0.4 mm until a total displacement of 0.8 mm is reached. The material property values recorded in the experiments are not adopted in this example, because numerous studies have shown that using the material properties reported in the experiments usually leads to a structural response in the linear elastic region that is inconsistent with the experimental results [34,60]. Therefore, this example adopts the material properties calibrated by Wu et al. [61], and the specific values are given in TABLE 7, where $f_c/f_t$ = 10 is consistent with Low et al. [34].

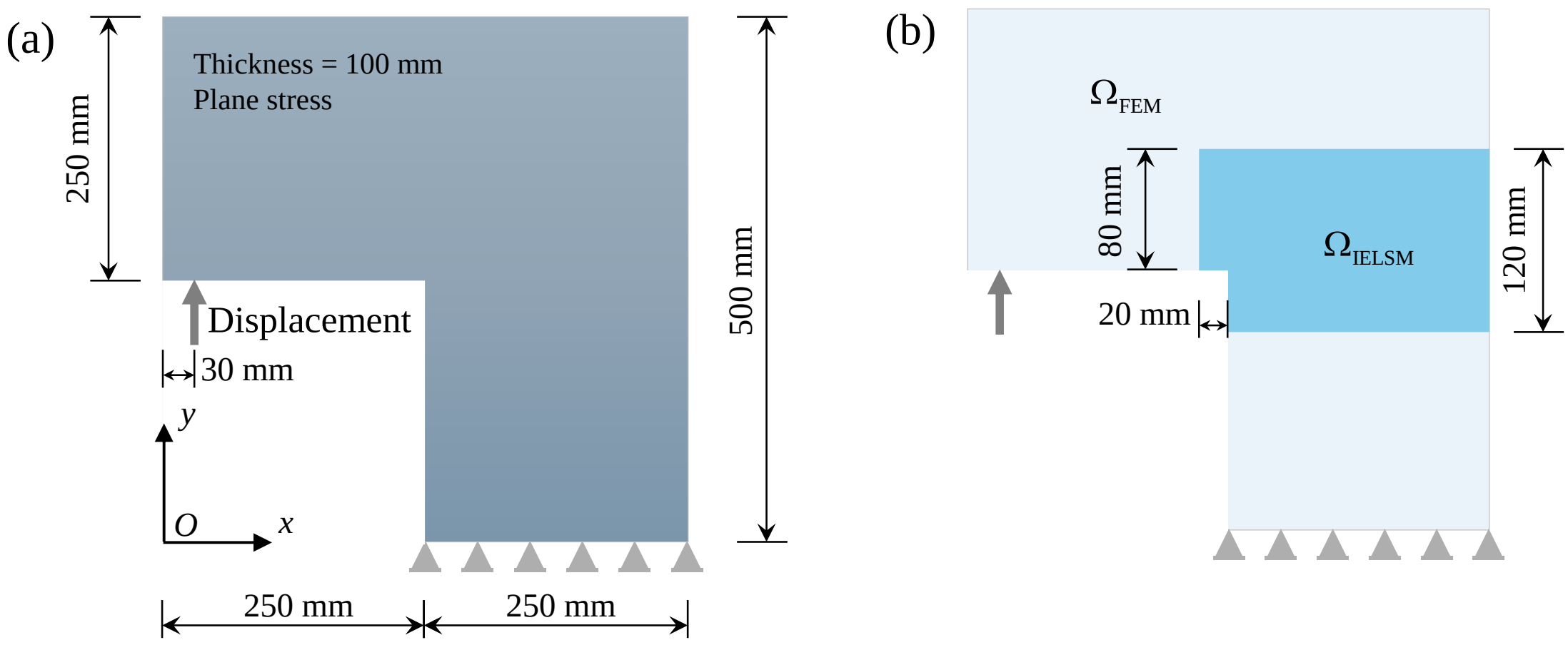

FIGURE 28 L-shaped panel: (a) geometry and boundary conditions; (b) partition of the IELSM region $\Omega_{IELSM}$ and the FEM region $\Omega_{FEM}$.

TABLE 7 L-shaped panel: Material properties.

| Symbol | Value |
|---|---|
| $E$ | 20 GPa |
| $v$ | 0.18 |
| $f_t$ | 2.5 MPa |
| $f_c$ | 25 MPa |
| $G_f$ | 130 J/m$^2$ |

As shown in FIGURE 29, the Force-Displacement curves computed using four element types under different damage criteria are compared. FIGURE 29(c) uses the Drucker-Prager damage criterion, and its results best match the range given by the experimental results. Then, based on the Drucker-Prager damage criterion, the coupling performance of IELSM and FEM is verified. As shown in FIGURE 30, the Force-Displacement curves predicted by IELSM-SQ/CST and IELSM-SQ with different element sizes are compared, where IELSM-SQ refers to the results computed using full-domain IELSM SQ elements. It is observed that the two are basically consistent. The cause of the deviation can be seen in FIGURE 32: the IELSM-SQ result exhibits some damage at the fixed end of the L-shaped concrete panel, whereas the IELSM-FEM coupling scheme does not account for this part of the damage. Based on the Drucker-Prager damage criterion, FIGURE 31 compares the Force-Displacement curves computed using different element sizes and element types with the experimental range. The IELSM-HEX/CST results agree well with the experimental range; the IELSM-TRI/CST results essentially coincide with the FEM-CST results, and the predicted peak force is relatively high; the peak force predicted by IELSM-SQ/CST is somewhat larger. The curves given by all schemes agree with the experimental results in the post-peak softening stage.

FIGURE 32 shows the damage contours at the end of loading computed considering different element types and element sizes. It can be seen that IELSM-HEX/CST has the smallest damage band width, which is also the reason for its smallest predicted peak force. Furthermore, FIGURE 33 extracts the crack contours predicted by different element types and element sizes and compares them with the experimental crack range, as well as with Huang et al. [62] and Low et al. [34]. The results show that, as the element size is refined, the crack contours predicted by different element types exhibit good convergence. Comparing the crack results of IELSM-SQ/CST and IELSM-SQ, it can be seen that the two are completely consistent, indicating that coupling does not affect crack propagation. Regarding the accuracy of crack propagation, the crack range of IELSM-TRI/CST agrees with the experiments and is better than the results of Huang et al. [62] and Low et al. [34]; the results of IELSM-HEX/CST and FEM-CST are close to those of Huang et al. [62] and Low et al. [34]; and the crack of IELSM-SQ/CST deviates from the experimental results and those of other numerical schemes. As mentioned above, the crack predicted by IELSM-SQ/CST has a certain mesh dependence and propagates along the mesh direction. Overall, IELSM-HEX/CST can provide Force-Displacement curves and crack paths that agree well with the experiments.

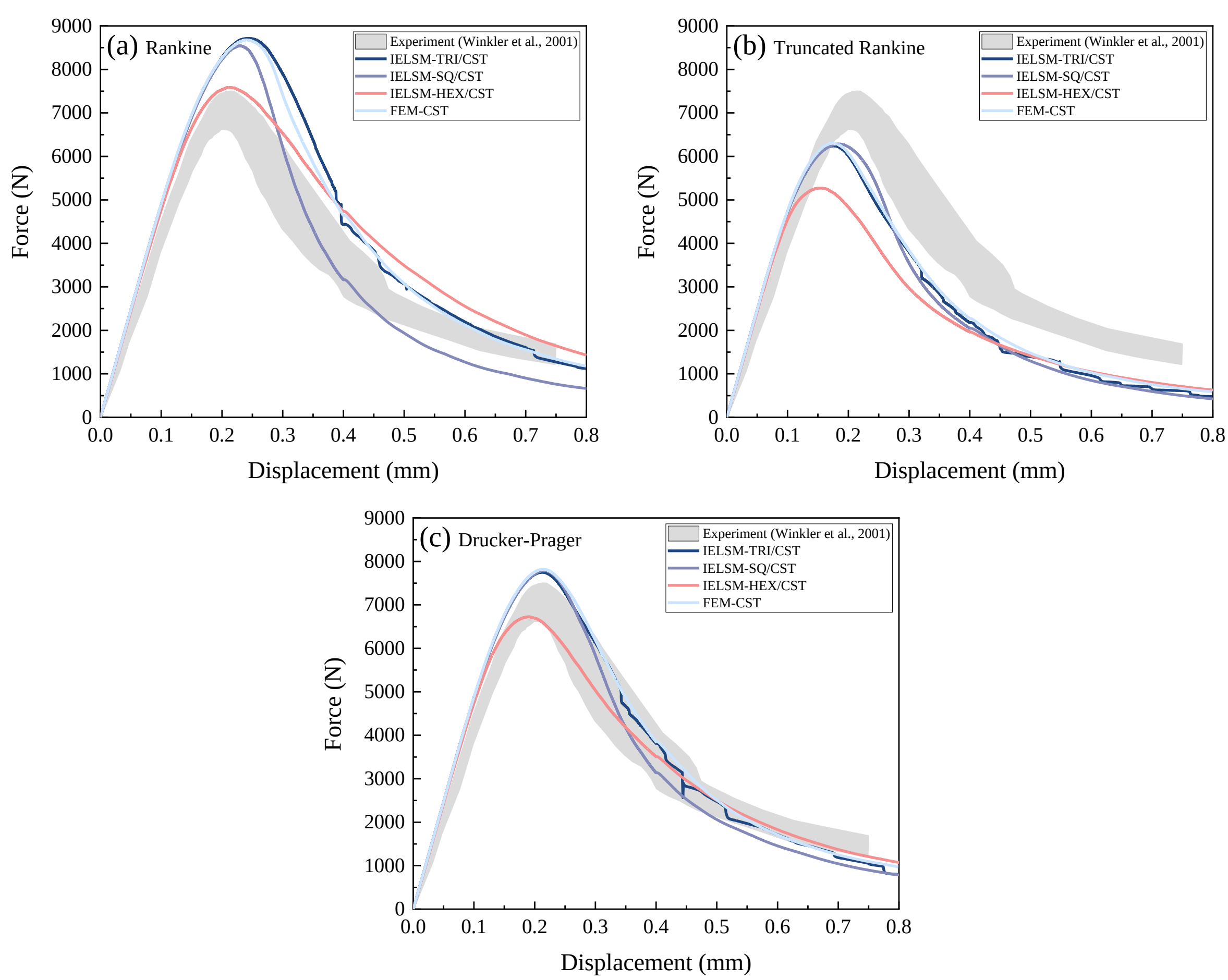

FIGURE 29 L-shaped panel: Force-Displacement curves computed using four element types under different damage criteria: (a)Rankine, (b)Truncated Rankine, (c)Drucker-Prager.

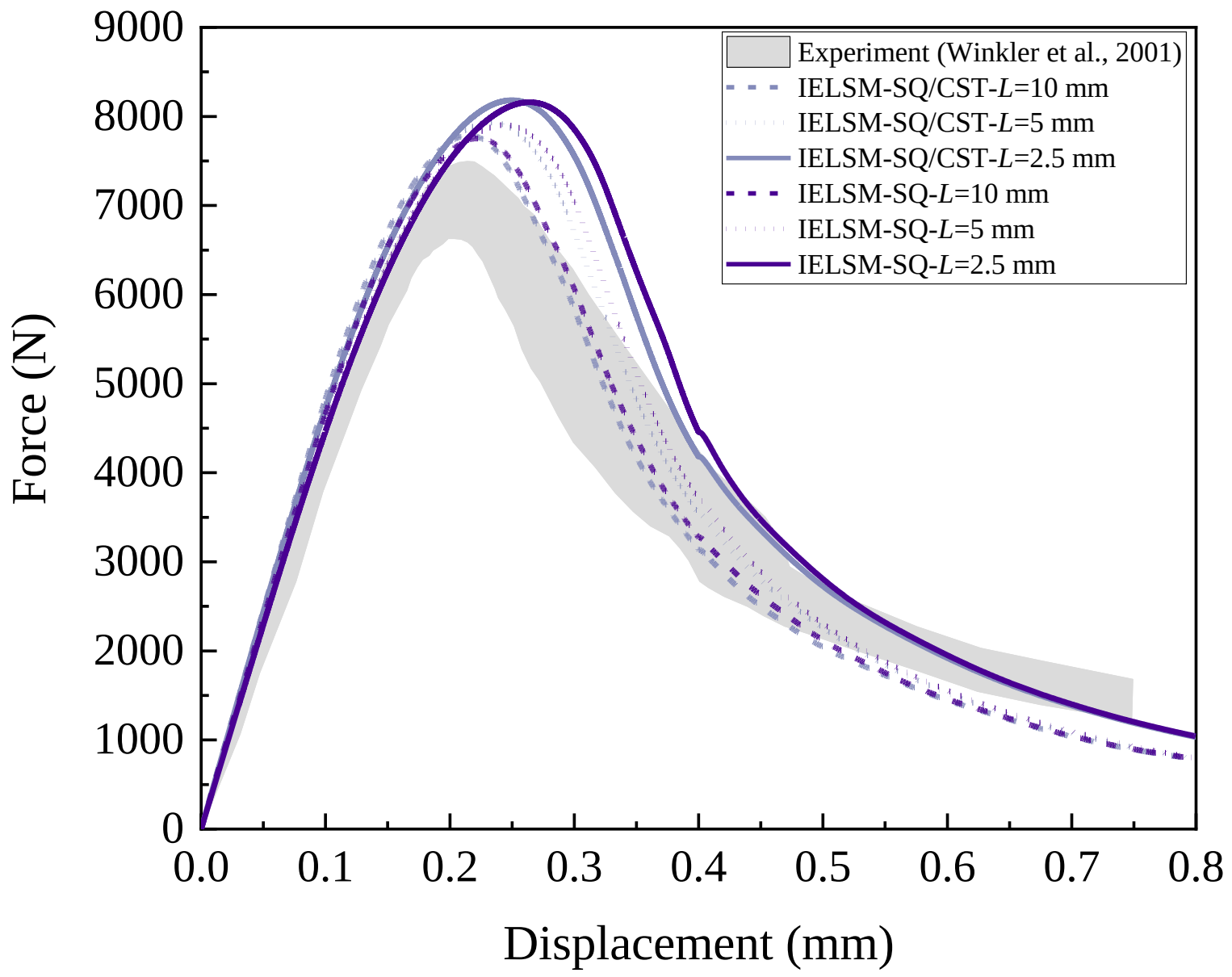

FIGURE 30 L-shaped panel: Comparison of Force-Displacement curves predicted by IELSM-SQ/CST and IELSM-SQ with different element sizes.

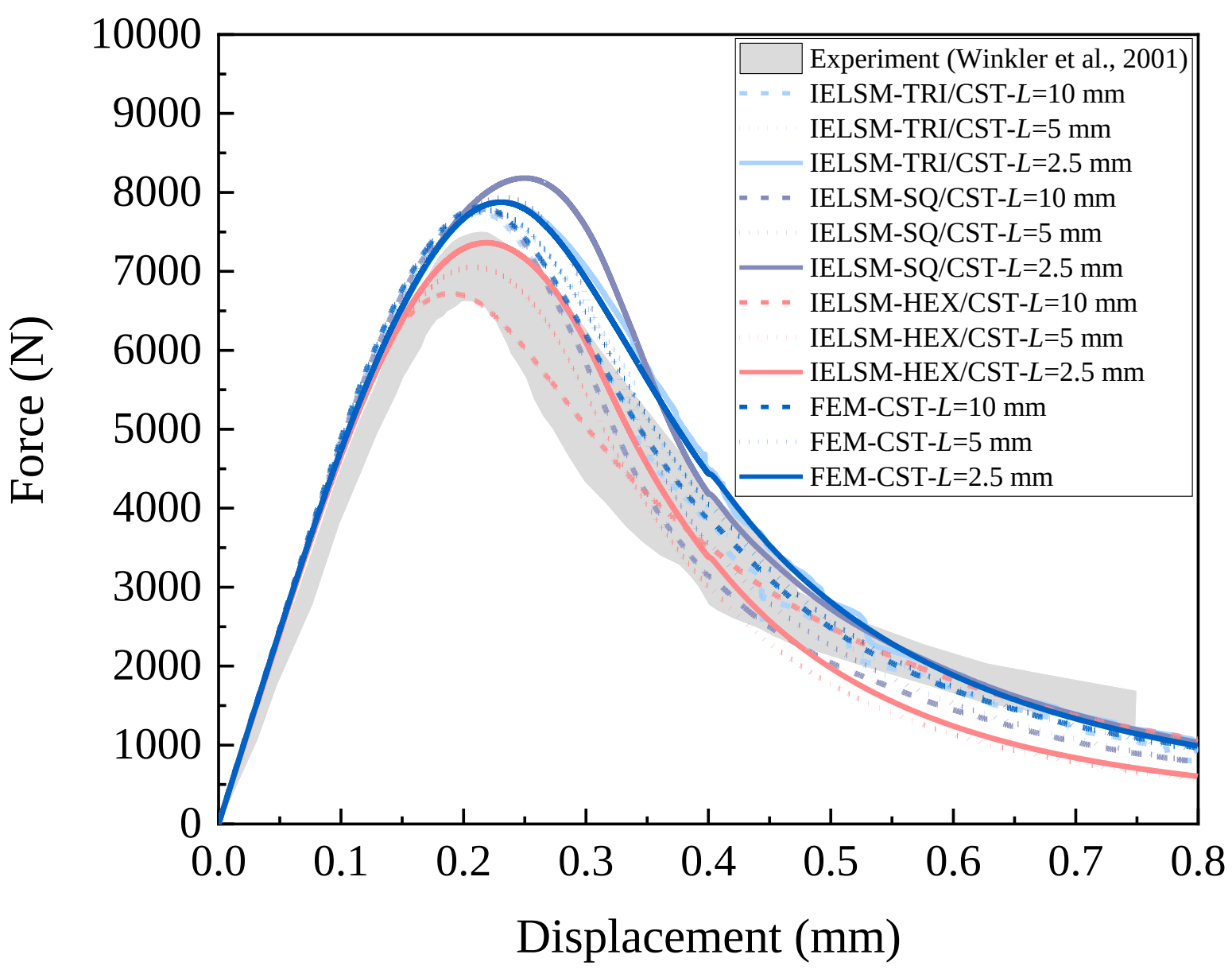

FIGURE 31 L-shaped panel: Force-Displacement curves computed using different element sizes and element types.

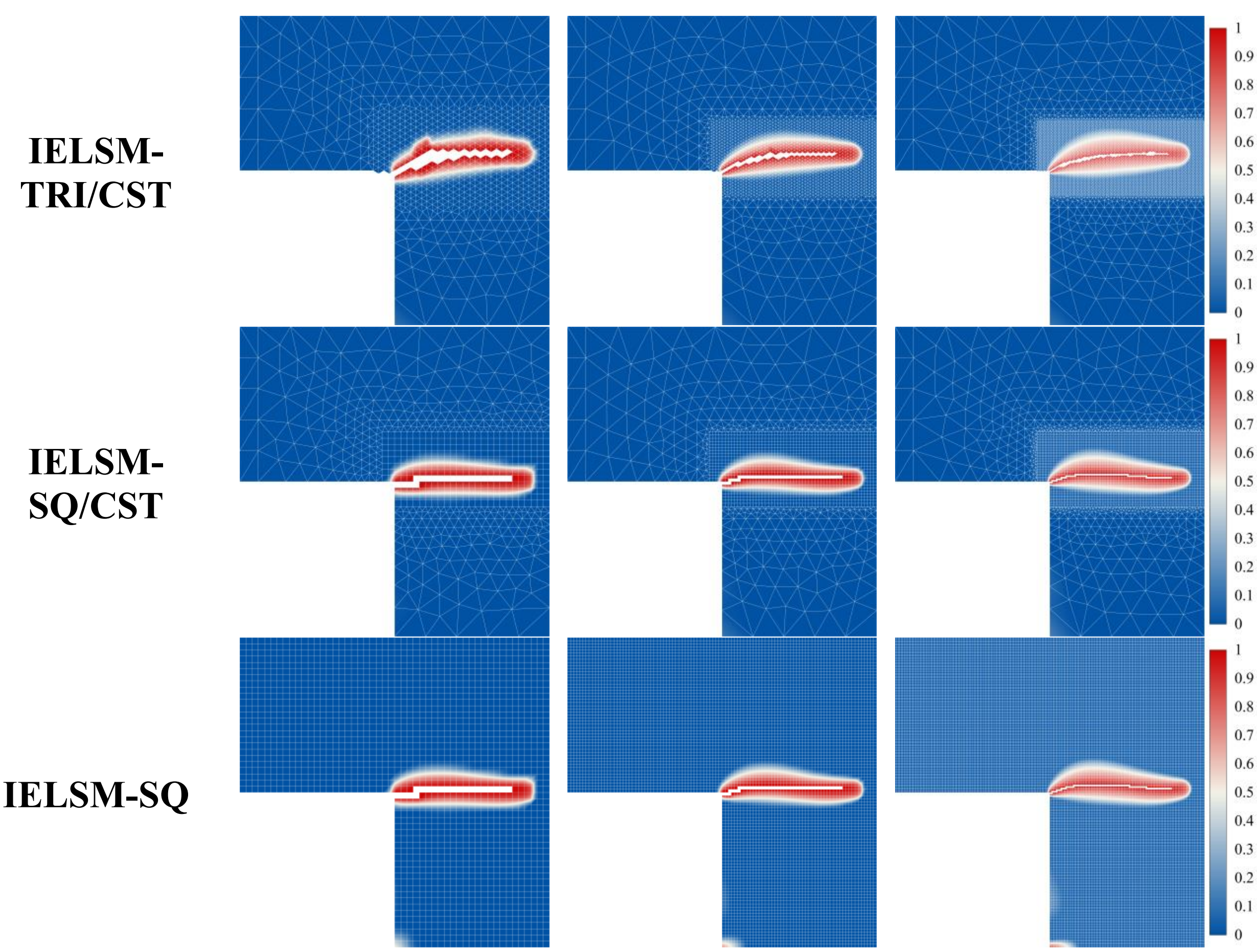

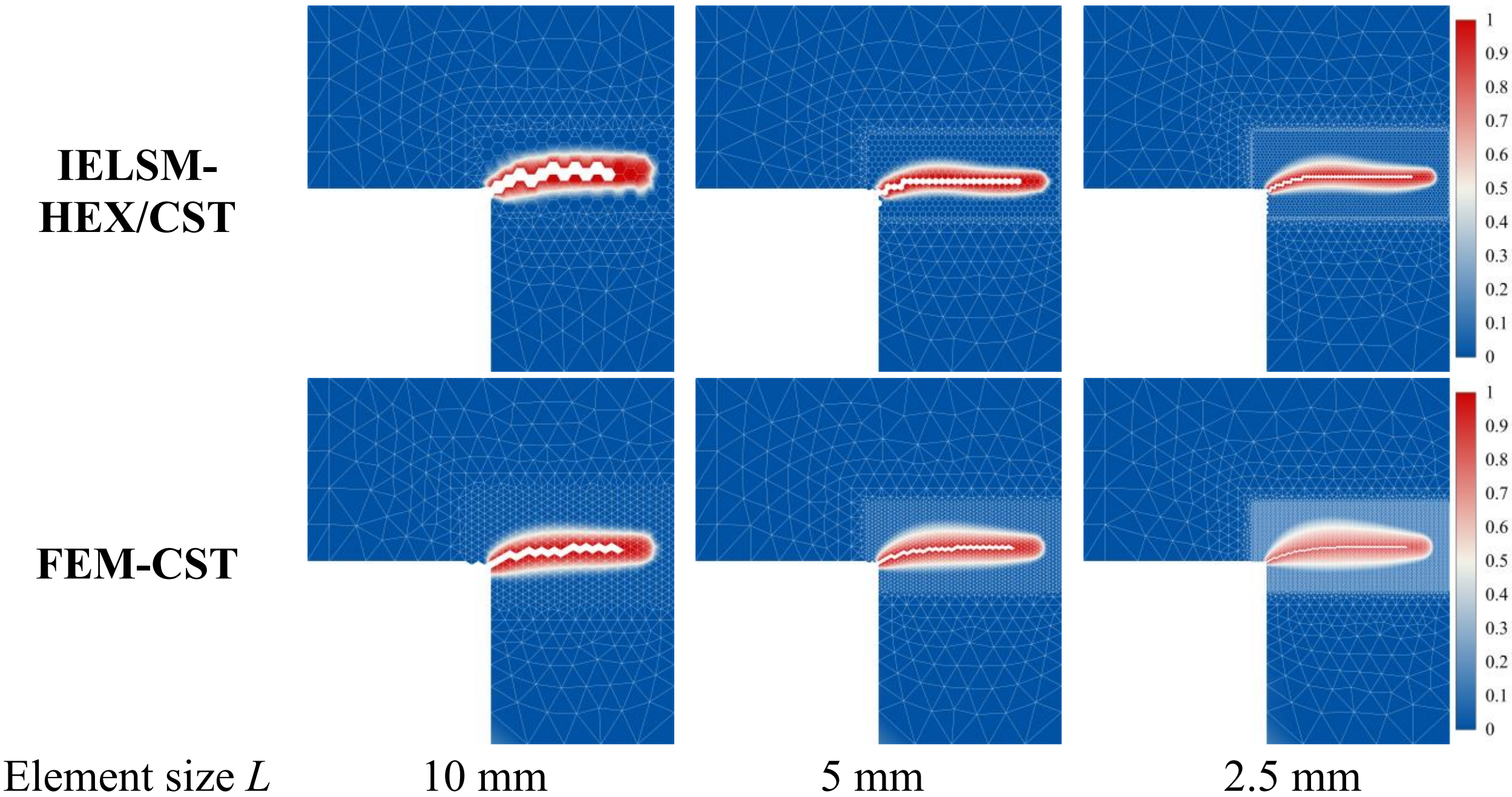

FIGURE 32 L-shaped panel: Damage contours at the end of loading computed using different element sizes and element types. (Displacement = 0.8mm)

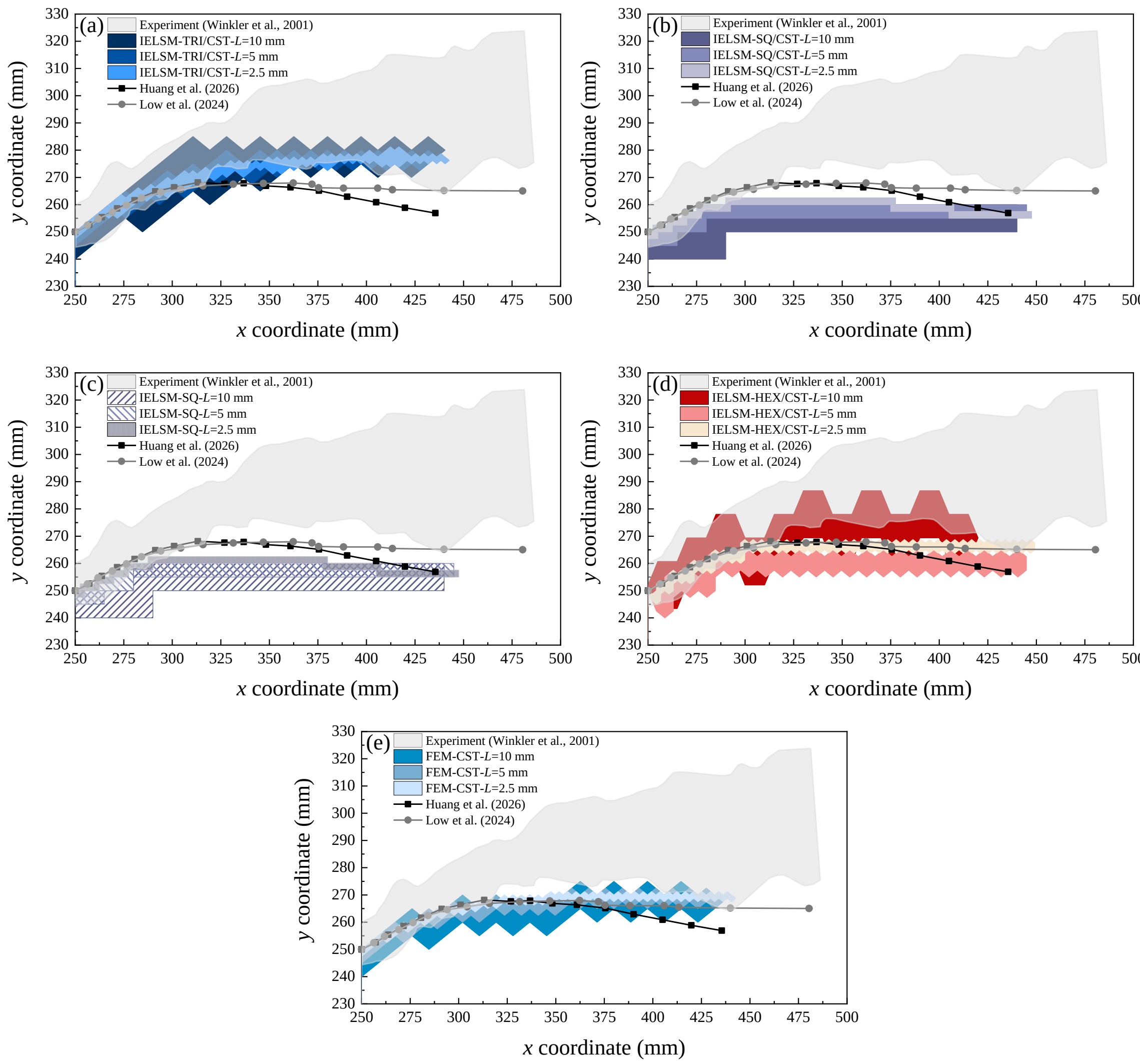

(a)
Experiment (Winkler et al., 2001)
IELSM-TRI/CST-L=10 mm
IELSM-TRI/CST-L=5 mm
IELSM-TRI/CST-L=2.5 mm
Huang et al. (2026)
Low et al. (2024)
y coordinate (mm)
x coordinate (mm)

(b)
Experiment (Winkler et al., 2001)
IELSM-SQ/CST-L=10 mm
IELSM-SQ/CST-L=5 mm
IELSM-SQ/CST-L=2.5 mm
Huang et al. (2026)
Low et al. (2024)
y coordinate (mm)
x coordinate (mm)

(c)
Experiment (Winkler et al., 2001)
IELSM-SQ-L=10 mm
IELSM-SQ-L=5 mm
IELSM-SQ-L=2.5 mm
Huang et al. (2026)
Low et al. (2024)
y coordinate (mm)
x coordinate (mm)

(d)
Experiment (Winkler et al., 2001)
IELSM-HEX/CST-L=10 mm
IELSM-HEX/CST-L=5 mm
IELSM-HEX/CST-L=2.5 mm
Huang et al. (2026)
Low et al. (2024)
y coordinate (mm)
x coordinate (mm)

(e)
Experiment (Winkler et al., 2001)
FEM-CST-L=10 mm
FEM-CST-L=5 mm
FEM-CST-L=2.5 mm
Huang et al. (2026)
Low et al. (2024)
y coordinate (mm)
x coordinate (mm)

FIGURE 33 L-shaped panel: Crack contours predicted using different element types and element sizes: (a) IELSM-TRI/CST, (b) IELSM-SQ/CST, (c) IELSM-SQ, (d) IELSM-HEX/CST, (e) FEM-CST.

As shown in TABLE 8, the computational efficiencies of the coupled and uncoupled schemes in this example are compared, where the IELSM region accounts for 16.9% of the total geometric domain area. It can be seen from the table that the IELSM-FEM coupling scheme significantly reduces the computational scale. Compared with the fully tessellated IELSM-SQ, the coupled scheme IELSM-SQ/CST reduces the number of nodes by 64.3% – 80.9% and the CPU time by 59.1% – 85.2%, and the reduction in computational cost increases with mesh refinement. FIGURE 34 shows the CPU time versus nodes curves for each scheme in this example, and the computational efficiency follows the same trend as in FIGURE 27.

TABLE 8 L-shaped panel: Comparison of computational efficiency between coupled and uncoupled schemes.

| Element size $L$ (mm) | Nodes | | Reduction | CPU time (s) | | Reduction |
|---|---|---|---|---|---|---|
| | SQ/CST | Full-domain SQ | | SQ/CST | Full-domain SQ | |
| 10 | 706 | 1976 | 64.3% | 340 | 832 | 59.1% |
| 5 | 1681 | 7701 | 78.2% | 651 | 2941 | 77.9% |
| 2.5 | 5793 | 30401 | 80.9% | 1992 | 13493 | 85.2% |

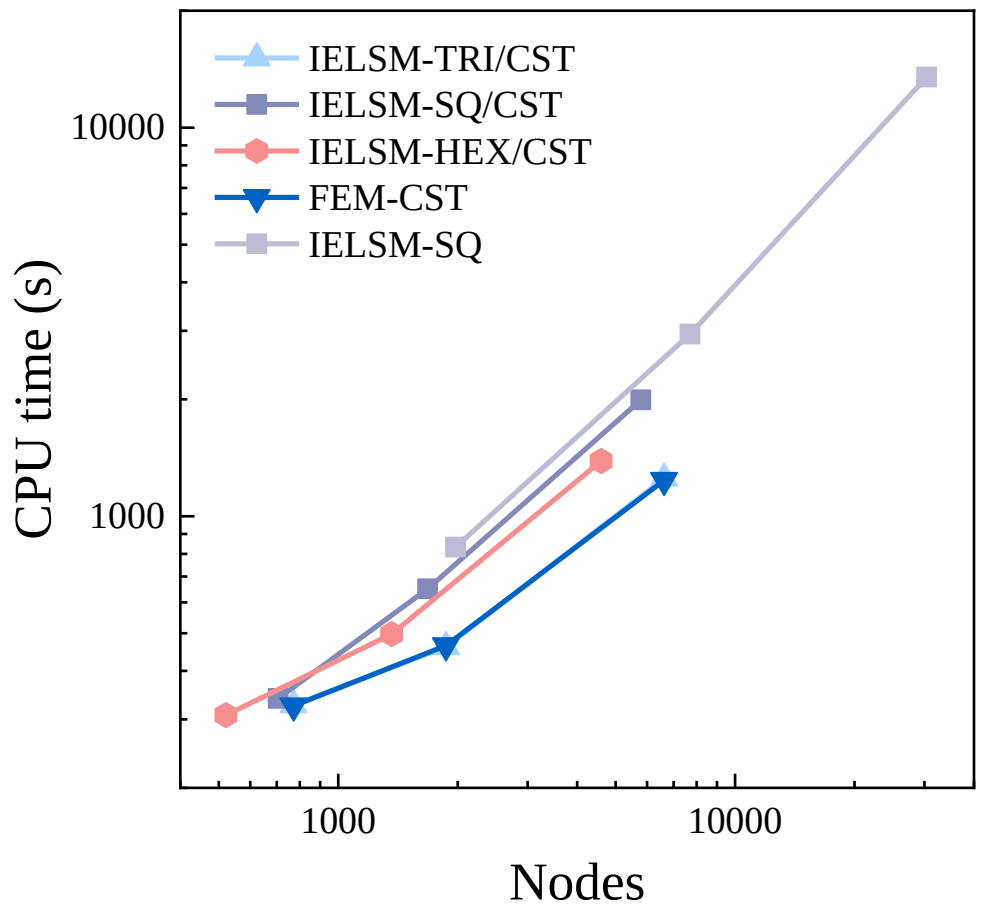

FIGURE 34 L-shaped panel: CPU time versus nodes curves for each scheme.

## 5.5 Double-Edged Notched specimen

This example simulates mixed-mode fracture of a double-edge-notched (DEN) specimen, with experimental data taken from Nooru-Mohamed et al. [63]. This work considers loading path "2" in Nooru-Mohamed et al. [63]. The geometric details and boundary conditions of the DEN specimen are shown in FIGURE 35(a), and FIGURE 35(b) shows the partition of the IELSM region and the FEM region in this example. Only the largest specimen in the experiments, DEN 200 (200 mm × 200 mm × 50 mm), is considered here, as it is widely used for benchmark validation of mixed-mode fracture in various numerical methods. The specimen thickness is 50 mm, and plane stress conditions are considered. An $x$-direction displacement is applied to the upper-left edge of the specimen, and a $y$-direction displacement of the same magnitude is applied to the top edge, so as to create a mixed loading mode; this corresponds to the boundary condition $\delta/\delta_s = 1.0$ in Nooru-Mohamed et al. [63]. The loading protocol is 0.00001 mm per step before a total displacement of 0.06 mm, and 0.00005 mm per step after 0.06 mm until a total displacement of 0.2 mm is reached. The normal/vertical deformation $U_n$ of the specimen is computed as the average relative vertical displacement between points M and M′ and between points N and N′ in FIGURE 35(a), and the shear/lateral deformation $U_s$ is given by the lateral displacement at point S. The material properties used in this example are taken from the experiments [63] and Low et al. [34], and the specific values are given in TABLE 9. Numerous studies on isotropic damage models have shown that the Ottosen failure criterion can capture the response of this benchmark test well [34,49]. This section adopts the Ottosen failure criterion for computation and comparison.

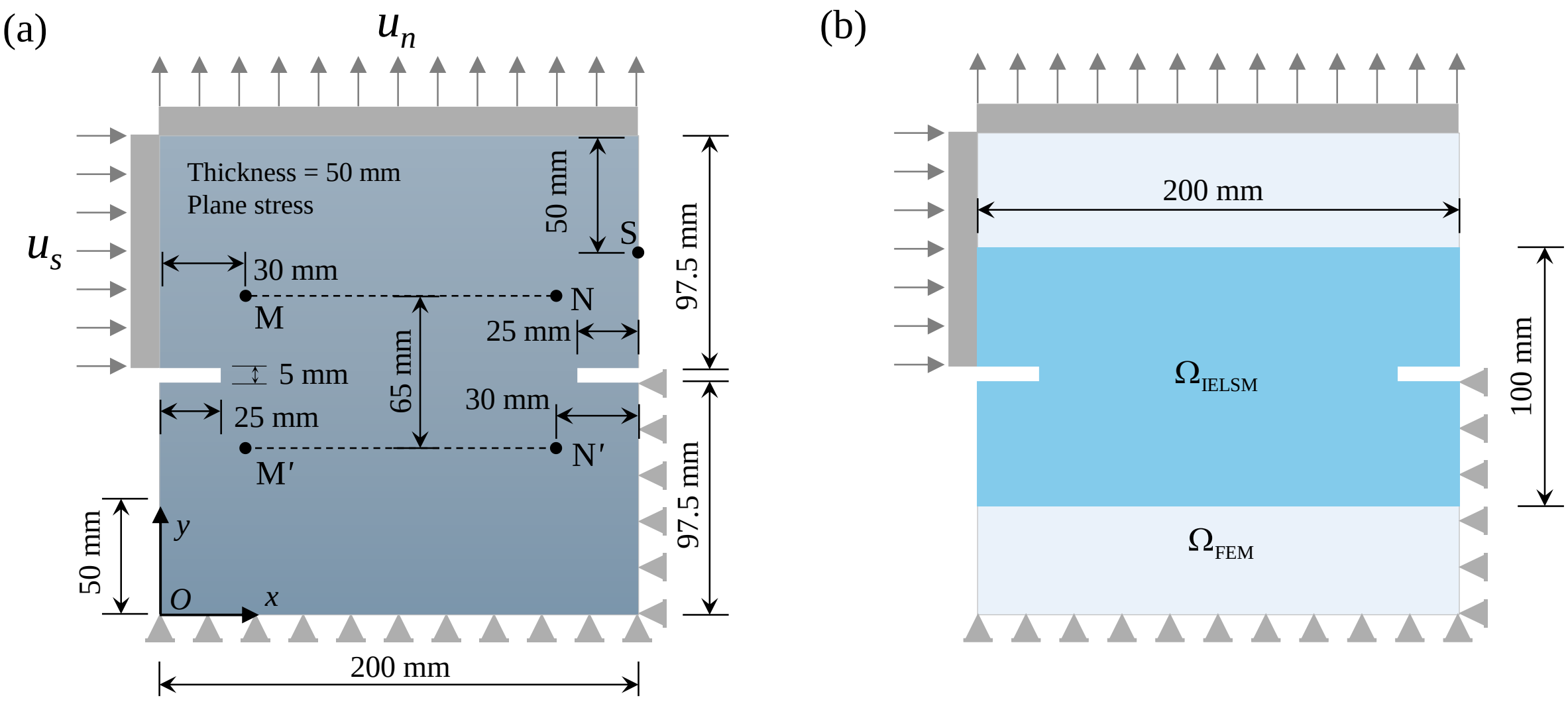

FIGURE 35 DEN: (a) geometry and boundary conditions; (b) partition of the IELSM region $\Omega_{IELSM}$ and the FEM region $\Omega_{FEM}$.

TABLE 9 DEN: Material properties.

| Symbol | Value |
|---|---|
| $E$ | 30 GPa |
| $v$ | 0.2 |
| $f_t$ | 3.0 MPa |
| $f_c/f_t$ | 14.4 |
| $f_s/f_t$ | 1.4 |
| $G_f$ | 110 J/m$^2$ |

FIGURE 36 presents the Force-$U_n$ curves computed using different element types and element sizes. In terms of computational accuracy, IELSM-TRI/CST and IELSM-HEX/CST are both superior to the results of Low et al. [34], and their peak loads fall within the experimental range, whereas IELSM-SQ/CST and FEM-CST are close to the results of Low et al. [34] before the softening stage but have relatively higher load levels in the final failure stage; the Force-$U_n$ curves obtained by each scheme with different mesh sizes show good consistency. FIGURE 37 presents the Force-$U_s$ curves computed using different element types and element sizes. The curve obtained by IELSM-HEX/CST best matches the experimental range, and the Force-$U_s$ curves obtained by each scheme with different mesh sizes show good convergence. FIGURE 38 shows the damage contours at the end of loading computed using different element types and element sizes, and FIGURE 39 extracts the cracks therein for comparison. The results show that the crack paths predicted by all schemes are generally consistent with the experimental observations. Finally, combining FIGURE 36(b), FIGURE 37(b), FIGURE 38, and FIGURE 39(b) and (c), it can be seen that the curves obtained by IELSM-SQ/CST and IELSM-SQ highly coincide, and the damage contours and crack propagation paths are completely consistent, indicating that the coupling of IELSM and FEM has no effect on the accuracy of the method.

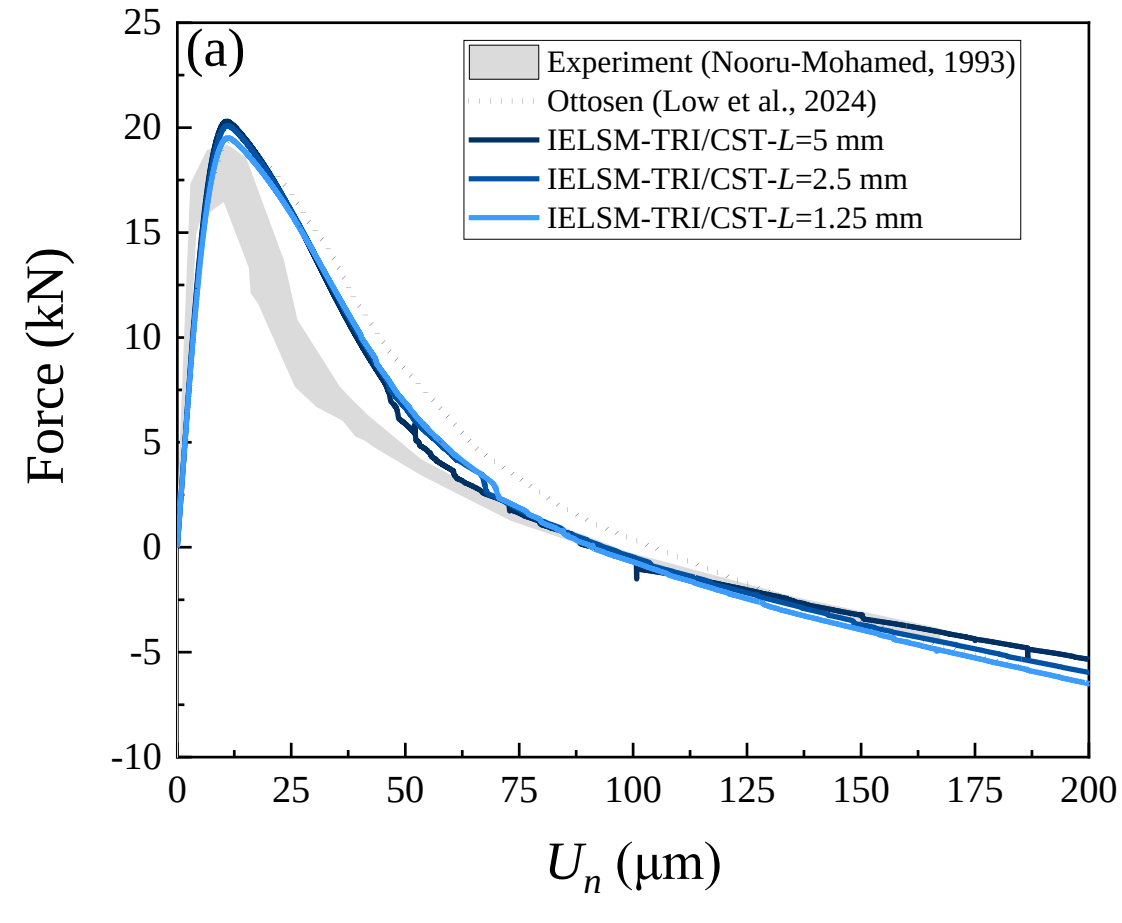

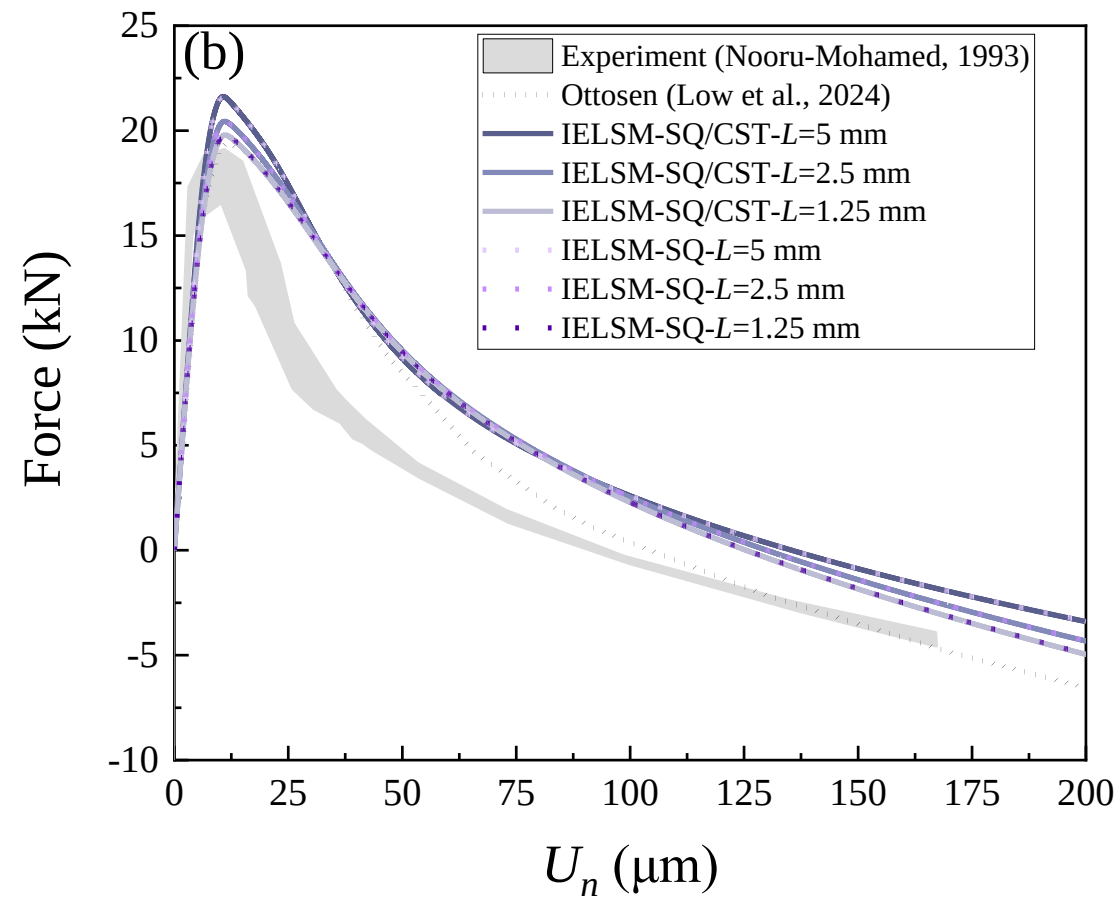

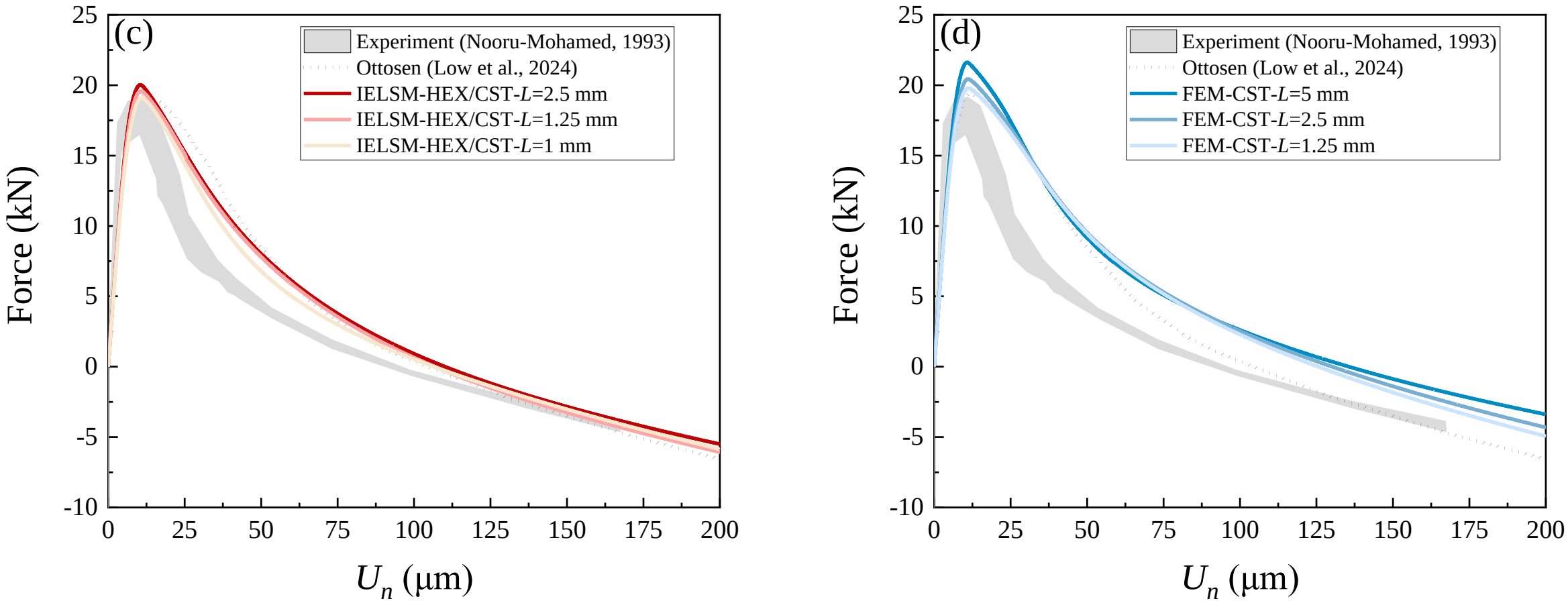

FIGURE 36 DEN: Force-$U_n$ curves computed using different element types and element sizes: (a) IELSM-TRI/CST, (b) IELSM-SQ/CST and IELSM-SQ, (c) IELSM-HEX/CST, (d) FEM-CST.

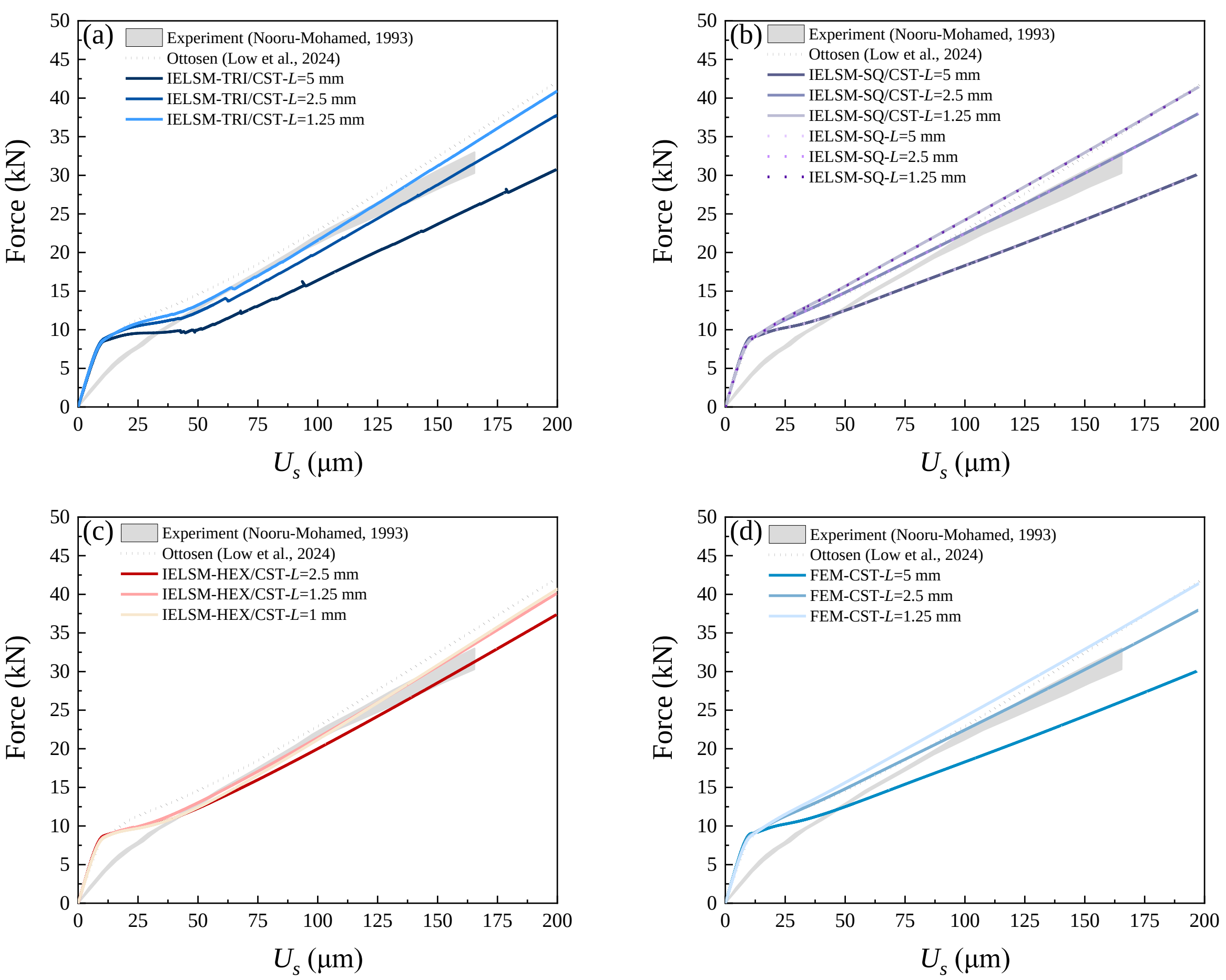

FIGURE 37 DEN: Force-$U_s$ curves computed using different element types and element sizes: (a) IELSM-TRI/CST, (b) IELSM-SQ/CST and IELSM-SQ, (c) IELSM-HEX/CST, (d) FEM-CST.

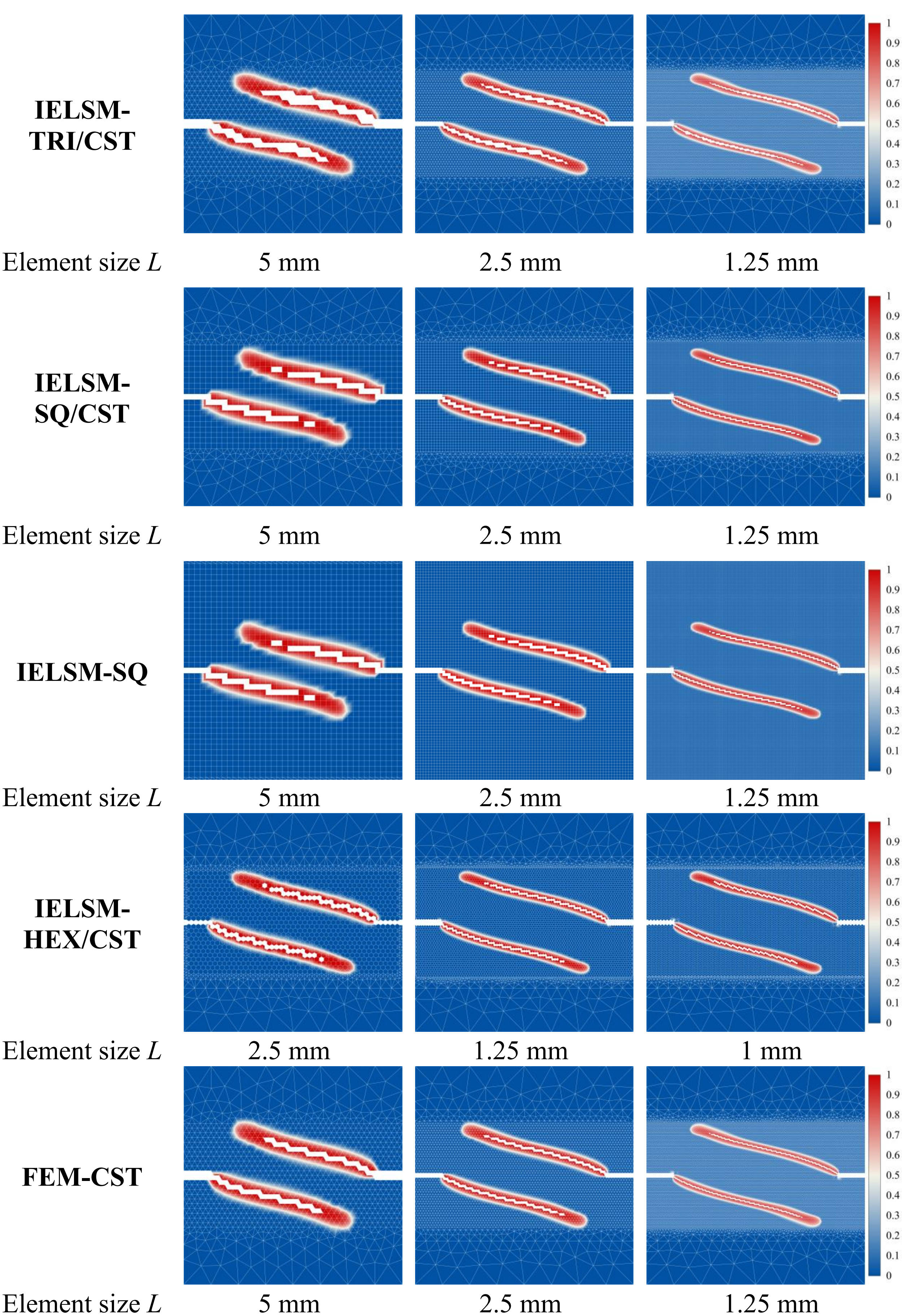

FIGURE 38 DEN: Damage contours at the end of loading computed using different element types and element sizes.

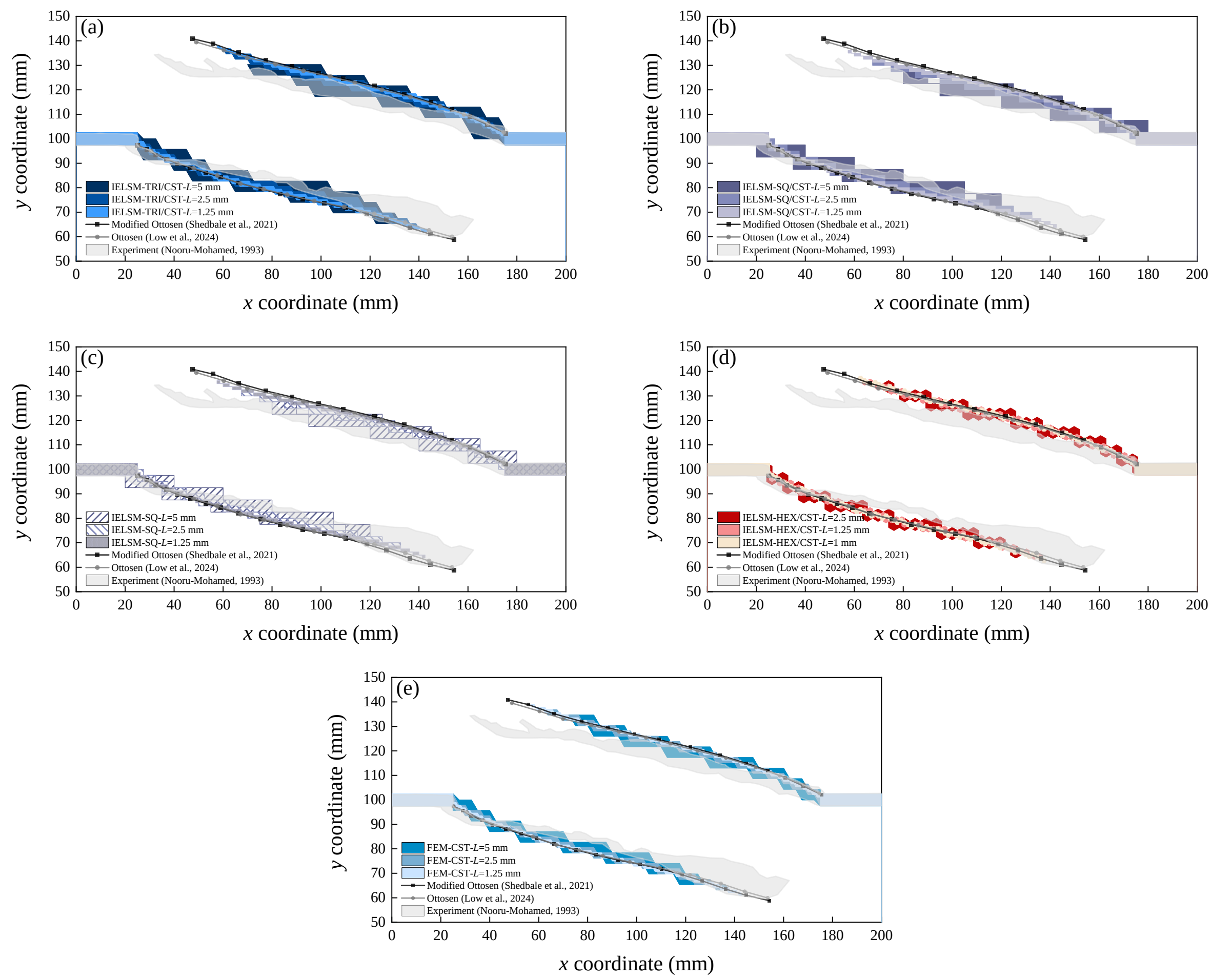
FIGURE 39 DEN: Crack contours computed using different element types and element sizes: (a) IELSM-TRI/CST, (B) IELSM-SQ/CST, (c)IELSM-SQ, (d)IELSM-HEX/CST, (e)FEM-CST.

As shown in TABLE 10, the computational performance of the coupled and uncoupled schemes in this example is compared, where the IELSM region accounts for 49.7% of the total geometric domain area. It can be seen from the table that the number of nodes is reduced by 41.2% – 47.5% and the CPU time by 34.2% – 48.5%. FIGURE 40 shows the CPU time versus nodes curves for each scheme in this example, and the computational efficiency follows the same trend as in FIGURE 27 and FIGURE 34. FIGURE 41 linearly fits the log(CPU time per step) versus log(Nodes) lines for the different schemes based on all data points in Examples 5.3 – 5.5, where the slope b$b$ of the line represents the computational complexity exponent in CPU time per step $\propto$ Nodes$^b$. The $b$ = 0.73 of FEM-CST is close to the $b$ = 0.76 of IELSM-TRI/CST; the $b$ = 0.89 of IELSM-HEX/CST is intermediate; the $b$ = 0.91 of IELSM-SQ/CST follows; and the $b$ = 1.00 of IELSM-SQ exhibits an approximately linear growth with the number of nodes. Overall, the computational efficiency of the IELSM-FEM coupling scheme is not significantly degraded by mixing the two discretization systems: the per-step solution time of IELSM-TRI/CST and its growth rate with the number of nodes are both comparable to those of FEM-CST, making it the most efficient among all schemes; IELSM-HEX/CST is slightly less efficient; while IELSM-SQ/CST and IELSM-SQ have the longest per-step solution times and the highest complexity exponents, so fully tessellated IELSM does not have an advantage in solution efficiency. The above differences arise from the smoothness of the transition between each IELSM element type and the FEM CST element, as well as from differences in the number of degrees of freedom and the sparsity structure of the stiffness matrix. Overall, for fracture problems characterized by damage localization, the IELSM-FEM coupling strategy achieves a good balance between accuracy (see Sections 5.2 – 5.5) and efficiency, providing a feasibility basis for the application of IELSM to relatively large-scale engineering problems.

TABLE 10 DEN: Comparison of computational efficiency between coupled and uncoupled schemes.

| Element size *L* (mm) | Nodes | | Reduction | CPU time (s) | | Reduction |
|---|---|---|---|---|---|---|
| | SQ/CST | Full-domain SQ | | SQ/CST | Full-domain SQ | |
| 5 | 1012 | 1722 | 41.2% | 896 | 1362 | 34.2% |
| 2.5 | 3667 | 6541 | 43.9% | 2551 | 4958 | 48.5% |
| 1.25 | 13545 | 25801 | 47.5% | 9941 | 19303 | 48.5% |

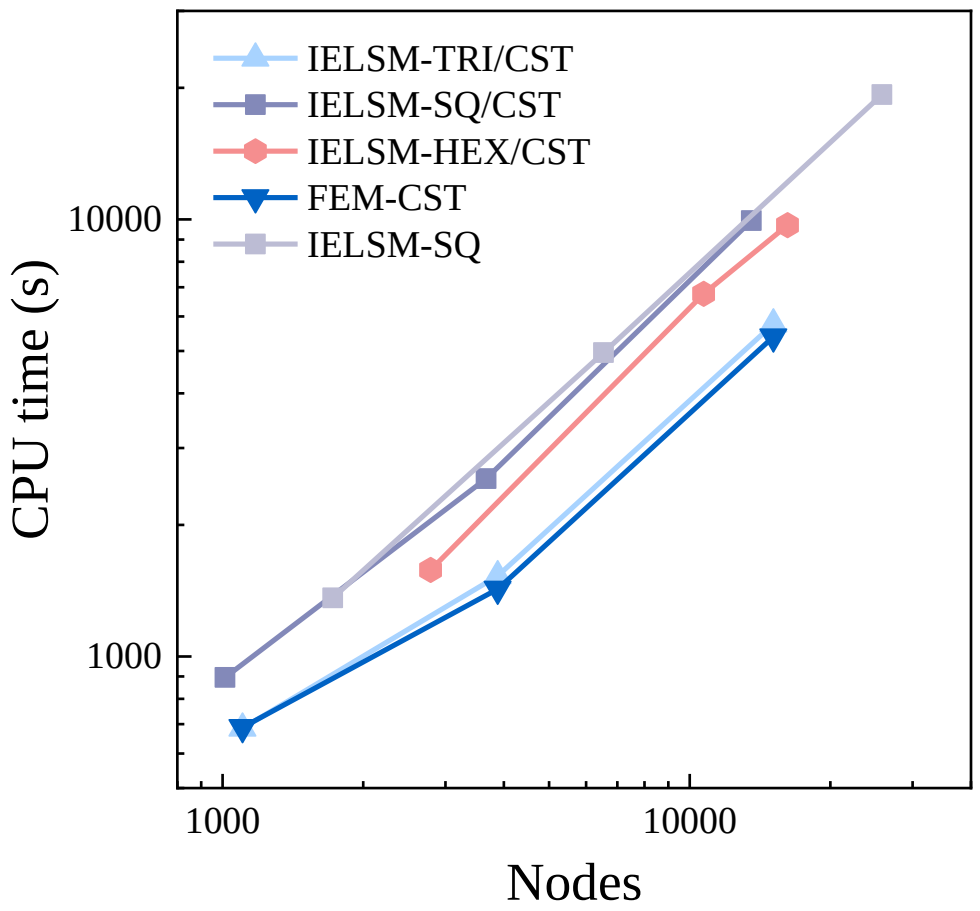

FIGURE 40 DEN: CPU time versus nodes curves for each scheme.

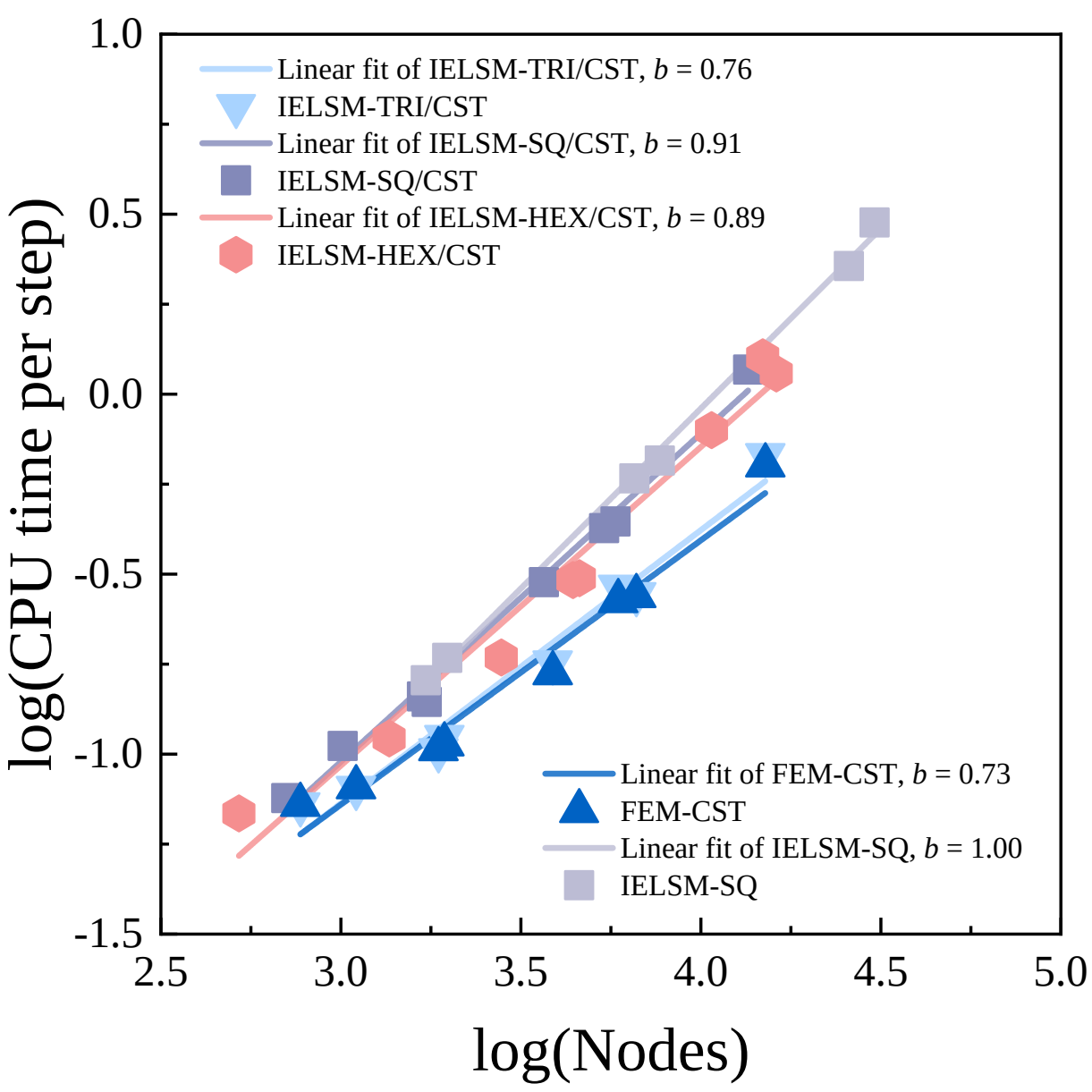

FIGURE 41 DEN: Linear fits of log(CPU time per step) versus log(Nodes) for the different schemes based on all data points in Sections 5.2 – 5.5.

## 6. Conclusion

To address the conflict among accuracy, efficiency, and implementation simplicity in the numerical simulation of quasi-brittle fracture, this paper proposes the Tessellated Isotropic Elastic Lattice Spring Model (IELSM), establishes an admissibility theory applicable to arbitrary polygonal elements, achieves its coupling with standard finite elements and an isotropic damage model, and systematically validates it through multiple benchmark examples. The main conclusions are as follows:

(1) A unified IELSM theory applicable to arbitrary element shapes is established. The macroscopic isotropy requirement is reduced to four rotational invariance conditions and one elastic parameter matching condition, forming five governing equations. The solvability of this equation system serves as the theoretical criterion for the admissibility of an element shape. By introducing symmetry and equal-stiffness conditions, it is proved that equilateral triangles, rectangles satisfying positive-definiteness constraints, regular hexagons, and arbitrary regular $N$-gons ($N \geq 5$) can all serve as admissible elements. An L-shaped concave polygon is used as an example to verify the applicability of the theory to concave elements. As a result, the admissible elements for regular meshes are extended from the conventional equilateral triangle, square, and regular hexagon to arbitrary element shapes and their combinations that can form a seamless tessellation. In addition, for the general underdetermined case, quadratic programming is identified as a feasible approach for solving the spring stiffness distribution of arbitrary polygonal elements.

(2) Coupling IELSM elements with CST elements of FEM greatly optimizes the allocation of computational cost without degrading the accuracy of the coupling scheme. The pure bending convergence test shows that the displacement field error convergence rates of the IELSM-FEM coupling schemes ($R$ = 2.05 – 2.19) are higher than those of full-domain IELSM ($R$ = 2.01) and full-domain FEM ($R$ = 1.88 – 2.01); the displacement field accuracy is close to that of full-domain IELSM and better than that of full-domain FEM. The stress/strain field convergence rates are comparable to those of full-domain IELSM, with local errors only at the coupling boundaries due to the post-processing scheme, but this does not affect the overall displacement field accuracy. In fracture examples such as the L-shaped panel and the double-edge-notched specimen, the load-displacement curves of IELSM-SQ/CST and fully tessellated IELSM-SQ highly coincide, and the damage contours and crack propagation paths are completely consistent, indicating that coupling does not adversely affect fracture prediction.

(3) When using the TRI, SQ, and HEX elements of IELSM to simulate fracture, the HEX element (H15 spring arrangement) performs best. It predicts the smallest damage region and the least dissipated energy, thereby yielding a lower peak force that is more consistent with the experimental range than TRI and SQ. In terms of Force-displacement (or Force-CMOD) curves, the HEX element results agree best with the experiments, followed by TRI; SQ is close to the full-domain FEM scheme and is relatively the worst. In terms of crack paths, the crack contours predicted by HEX are most consistent with the experimental and literature results, followed by TRI, whereas SQ exhibits a certain mesh-orientation dependence, with cracks tending to propagate along the mesh orientation. These differences arise from the different spring topologies and orientation distributions of the elements: HEX-H15 has the largest number of springs and the most uniform orientation distribution, so damage localizes earliest and energy is released in a concentrated manner into the main crack, effectively suppressing spurious damage diffusion and secondary crack branching.

(4) By restricting IELSM to the damage-prone region, compared with fully tessellated IELSM discretization, the number of nodes is reduced by 41.2% – 80.9% and the CPU time by 34.2% – 85.2%, and the reduction in computational cost increases with mesh refinement. Efficiency analysis shows that the computational complexity exponent $b$ of the per-step solution time (the slope of log(CPU time per step) versus log(Nodes)) for each scheme is: FEM-CST $b$ = 0.73, IELSM-TRI/CST $b$ = 0.76, IELSM-HEX/CST $b$ = 0.89, IELSM-SQ/CST $b$ = 0.91, and IELSM-SQ $b$ = 1.00. Among them, the per-step solution time of IELSM-TRI/CST and its growth rate with the number of nodes are comparable to those of FEM-CST, making it the most efficient tier among all schemes; IELSM-HEX/CST is intermediate; IELSM-SQ/CST and IELSM-SQ have the longest per-step solution times and the highest complexity exponents, so fully tessellated IELSM does not have an advantage in solution efficiency. Overall, the IELSM-FEM coupling strategy achieves a good balance between accuracy and efficiency, providing a feasibility basis for the application of IELSM to relatively large-scale engineering problems.

Future work will include, in addition to the already mentioned development of arbitrary polygonal elements via quadratic programming, verification of the numerical performance of IELSM-Q4 coupling, optimization of post-processing procedures, and quantitative comparison of different damage definitions, the development of three-dimensional IELSM, large deformation theory, and adaptive mesh refinement.

**Declaration of competing interest**

The authors declare that they have no known competing financial interests or personal relationships that could have appeared to influence the work reported in this paper.

During the preparation of this work, the authors used DeepSeek to correct grammatical mistakes in the original text, ensure that the paragraph structure and the language were clear and cohesive, and check some algebraic and asymptotic calculations. After using this tool, the authors reviewed and edited the content as needed and take full responsibility for the content of the published article.

**Acknowledgements**

DM thanks the supports of the National Natural Science Foundation of China (12002247) and the Fundamental Research Fund for the Central Universities of China (WUT: 2021IVB013).

## Appendix A

### *A.1.* Proof of infeasibility for arbitrary triangular elements

As shown in FIGURE 2, let the weights of the three edge springs of a triangular element in the total strain energy be $w_j = k_j L_j^2$, with orientation angles $\theta_j$ ($j = 1, 2, 3$). Introducing the complex variables $\xi_j = e^{i2\theta_j} = \cos 2\theta_j + i\sin 2\theta_j$ and $\xi_j^2 = e^{i4\theta_j} = \cos 4\theta_j + i\sin 4\theta_j$, the four real equations in Eq. (9) can be combined pairwise into two complex equations:

$$\begin{cases} \sum_{j=1}^{3} w_j \cos 2\theta_j + i\sum_{j=1}^{3} w_j \sin 2\theta_j = \sum_{j=1}^{3} w_j e^{i2\theta_j} = \sum_{j=1}^{3} w_j \xi_j = 0, \\ \sum_{j=1}^{3} w_j \cos 4\theta_j + i\sum_{j=1}^{3} w_j \sin 4\theta_j = \sum_{j=1}^{3} w_j e^{i4\theta_j} = \sum_{j=1}^{3} w_j \xi_j^2 = 0. \end{cases} \tag{A.1}$$

We examine whether there exists a real sequence $\{w_j\}$ satisfying (A.1) that is not identically zero. From the first equation of (A.1), we obtain:

$$w_3 = -\frac{w_1\xi_1 + w_2\xi_2}{\xi_3}, \tag{A.2}$$

substituting into the second equation gives:

$$w_1\xi_1(\xi_1 - \xi_3) + w_2\xi_2(\xi_2 - \xi_3) = 0. \tag{A.3}$$

Therefore, if a nonzero solution exists, the weight ratios must satisfy:

$$\frac{w_1}{w_2} = -\frac{\xi_2(\xi_2 - \xi_3)}{\xi_1(\xi_1 - \xi_3)}, \tag{A.4}$$

Since $w_j$ are real, the right-hand side of the above expression must be real. A complex number is real if and only if it equals its own conjugate, i.e.,

$$\frac{\xi_2(\xi_2 - \xi_3)}{\xi_1(\xi_1 - \xi_3)} = \frac{\overline{\xi_2}(\overline{\xi_2} - \overline{\xi_3})}{\overline{\xi_1}(\overline{\xi_1} - \overline{\xi_3})} = \frac{\xi_2^{-1}\frac{\xi_2 - \xi_3}{\xi_2\xi_3}}{\xi_1^{-1}\frac{\xi_1 - \xi_3}{\xi_1\xi_3}} = \frac{\xi_2 - \xi_3}{\xi_1 - \xi_3} \cdot \frac{\xi_1^{\ 2}}{\xi_2^{\ 2}}. \tag{A.5}$$

Simplifying yields:

$$\xi_1^{\ 3} = \xi_2^{\ 3}. \tag{A.6}$$

By symmetry, we similarly obtain $\xi_2^{\ 3} = \xi_3^{\ 3}$, and hence:

$$\xi_1^{\ 3} = \xi_2^{\ 3} = \xi_3^{\ 3}, \tag{A.7}$$

that is:

$$e^{i6\theta_1} = e^{i6\theta_2} = e^{i6\theta_3}. \tag{A.8}$$

Therefore, $6\theta_1$, $6\theta_2$, and $6\theta_3$ differ pairwise by integer multiples of $2\pi$, i.e., the three orientation angles $\theta_j$ differ from one another by integer multiples of $\pi/3$. Considering that the three edges of a triangle are connected sequentially, the interior angles are determined by the differences between the edge orientation angles, and the sum of the three interior angles is $\pi$. The only geometry that satisfies this symmetry and can form a closed triangle is the equilateral triangle: in this case, $\theta_3 = \theta_1 + \pi/3$ and $\theta_2 = \theta_3 + \pi/3$, and $\xi_j$ then satisfies (A.1), so the homogeneous equations admit a nonzero solution ($w_1 = w_2 = w_3$). If the triangle is not equilateral, the cubic equality condition cannot be satisfied, and the imaginary part of the right-hand side of (A.4) is nonzero, so the equation system (A.1) has only the zero solution $w_1 = w_2 = w_3 = 0$, which contradicts the fifth condition (14). For any non-equilateral triangle, there exists no real sequence ($w_1$, $w_2$, $w_3$) that simultaneously satisfies the five isotropy equations. Therefore, a triangular element containing only three edge springs cannot be generalized to arbitrary shapes, and the classical LSM must use the equilateral triangle as its basic element. This proof reveals, from the perspective of complex algebra, the necessary role of symmetry in removing the overdeterminacy of the isotropy constraints.

## Appendix B

### *B.1.* Proof of arbitrary regular polygon elements

As shown in FIGURE B1, let the circumradius of a regular $N$-gon be $l_r$, with vertices located sequentially at angles $2\pi j/N$ ($j$ = 0, 1, …, $N$ - 1). All chords inside the polygon connecting any two vertices (including the edges) are grouped by a "step" $\delta$: the step $\delta$ indicates that the two vertices connected by the chord are separated by $\delta$ edges, where $\delta$ = 1, …, $N$/2. In general, the chord length corresponding to step $\delta$ is:

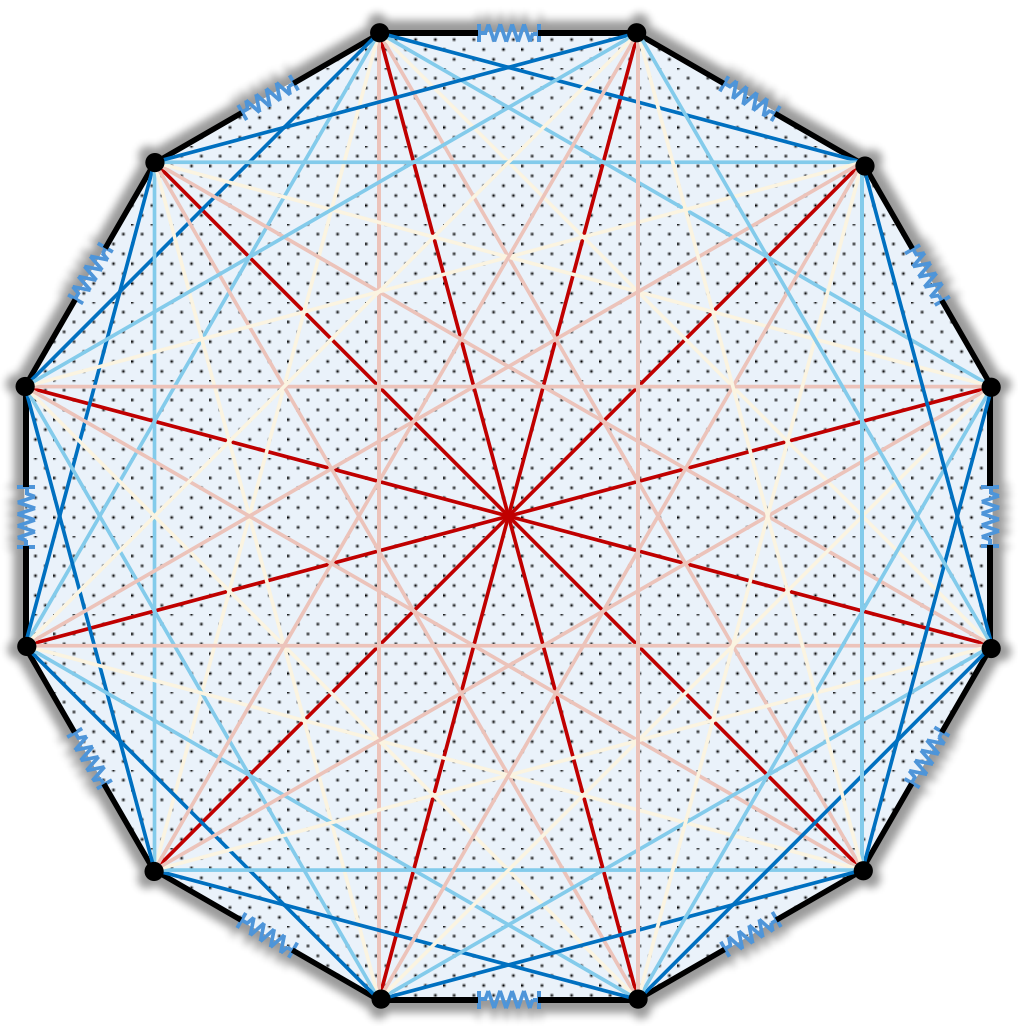

FIGURE B1 Schematic of spring arrangement within a regular polygonal element.

$$L_\delta = 2l_r \sin(\frac{\delta\pi}{N}). \tag{B.1}$$

Springs are arranged on every chord in each group. The spring stiffnesses of all springs within the same group are set equal (the stiffness of the $\delta$-th group is denoted by $k_\delta$), so the weight of each spring in the $\delta$-th group is $w_\delta = k_\delta L_\delta^2$. The number of chords in each group is:

$$n_\delta = \begin{cases} N, & \delta < \dfrac{N}{2}, \\ \dfrac{N}{2}, & \delta = \dfrac{N}{2} \text{ and } N \text{ is even.} \end{cases} \tag{B.2}$$

For a chord with step $\delta$, the indices of its two end vertices are $(j, j+\delta)$. From the vertex coordinates, the orientation angle of this chord is obtained as:

$$\theta_{j,\delta} = \frac{2\pi j}{N} + \frac{\pi\delta}{N} + \frac{\pi}{2} \pmod{2\pi}, \quad j = 0,1,\ldots,n_\delta - 1. \tag{B.3}$$

It can be seen that, for a fixed $\delta$, the orientation angles of the chords are uniformly distributed with an interval of $2\pi/N$ as $j$ varies. Referring to the complex-form rotational invariance conditions (A.1) in Appendix A, and summing by group, the contribution of the $\delta$-th group is

$$w_\delta \sum_{j=0}^{n_\delta - 1} e^{i2\theta_{j,\delta}}. \tag{B.4}$$

Substituting the orientation angle expression:

$$\sum_{j=0}^{n_\delta - 1} e^{i2\theta_{j,\delta}} = e^{i\left(\frac{2\delta}{N}+1\right)\pi} \sum_{j=0}^{n_\delta - 1} e^{i\frac{4\pi j}{N}}. \tag{B.5}$$

According to the root-of-unity summation identity, when $N \neq 2$, $\sum_{j=0}^{N-1} e^{i\frac{4\pi j}{N}} = 0$. For $N \geq 3$, this condition always holds. Therefore,

$$\sum_{j=0}^{n_\delta - 1} e^{i2\theta_{j,\delta}} = 0, \ \ (\delta < \frac{N}{2}) \tag{B.6}$$

when $\delta = N/2$ ($N$ even), the same summation reduces to $\sum_{j=0}^{N/2-1} e^{i\frac{4\pi j}{N}} = 0$ (since $N/2 > 1$) , so it also holds. Similarly, for the $4\theta$ terms:

$$\sum_{j=0}^{n_\delta - 1} e^{i4\theta_{j,\delta}} = e^{i\left(\frac{4\delta}{N}+2\right)\pi} \sum_{j=0}^{n_\delta - 1} e^{i\frac{8\pi j}{N}}. \tag{B.7}$$

When $N \neq 4$, this sum is zero. For $N \geq 5$, this condition always holds, so:

$$\sum_{j=0}^{n_\delta - 1} e^{i4\theta_{j,\delta}} = 0. \tag{B.8}$$

Since the above conclusions hold independently for each group $\delta$, their weighted sum remains zero. Therefore, the four rotational invariance conditions are automatically satisfied for any regular $N$-gon with $N \geq 5$, independently of the stiffness values of the spring groups. The only remaining constraint is the overall modulus equation, Eq. (14):

$$\sum_{\delta=1}^{N/2} n_\delta w_\delta = \frac{4SE}{1+\nu}, \tag{B.9}$$

where $S$ is the area of the regular $N$-gon:

$$S = \frac{N}{2} l_r^2 \sin\left(\frac{2\pi}{N}\right). \tag{B.10}$$

This equation contains $N/2$ unknowns $w_\delta$, so the system is underdetermined. To obtain a unique solution, the equal-stiffness assumption is introduced: all springs are assigned equal stiffness, i.e., $k_\delta = k$. Then $w_\delta = kL_\delta^2$, and the overall modulus condition becomes:

$$k \sum_{\delta=1}^{N/2} n_\delta L_\delta^2 = \frac{4SE}{1+\nu}. \tag{B.11}$$

Substituting (B.1) and using $\sum_{\delta=1}^{N-1}\sin^2\left(\frac{\delta\pi}{N}\right)=\frac{N}{2}$ together with the symmetric pairing of $\delta$ and $N-\delta$ yields (both odd and even $N$ give the same result):

$$k\sum_{\delta=1}^{N/2} n_\delta L_\delta^2 = N^2 l_r^2, \tag{B.12}$$

thus,

$$k=\frac{4SE}{(1+\nu)N^2 l_r^2}=\frac{4\frac{N}{2}l_r^2\sin\left(\frac{2\pi}{N}\right)E}{(1+\nu)N^2 l_r^2}=\frac{2E}{N(1+\nu)}\sin\left(\frac{2\pi}{N}\right). \tag{B.13}$$

Substituting $N = 6$ gives:

$$k=\frac{2E}{6(1+\nu)}\sin\left(\frac{2\pi}{6}\right)=\frac{\sqrt{3}E}{6(1+\nu)}, \tag{B.14}$$

which is consistent with the H15 scheme in Section 2.2.3.

## Appendix C

### *C.1.* The detailed forms of $\mathbf{K}_{\text{IELSM, Tri}}$, $\mathbf{K}_{\text{IELSM, H9}}$, $\mathbf{K}_{\text{IELSM, H12}}$, and $\mathbf{K}_{\text{IELSM, H15}}$ under plane strain

Element stiffness matrix $\mathbf{K}_{\text{IELSM, Tri}}$ of the equilateral triangular element:

$$\mathbf{K}_{\text{IELSM,Tri}} = \sum_{i=1}^{3}\mathbf{K}_{n,i} + \mathbf{K}_{volume}$$

$$= \frac{E}{(1+\nu)(1-2\nu)}\begin{bmatrix} \frac{\sqrt{3}(7-8\nu)}{24} & & & & & \\ \frac{1}{8} & \frac{\sqrt{3}(5-8\nu)}{24} & & & sym & \\ \frac{\sqrt{3}(4\nu-5)}{24} & \frac{(1-4\nu)}{24} & \frac{\sqrt{3}(7-8\nu)}{24} & & & \\ \frac{4\nu-1}{8} & \frac{\sqrt{3}(4\nu-1)}{24} & -\frac{\nu}{8} & \frac{\sqrt{3}(5-8\nu)}{24} & & \\ -\frac{\sqrt{3}(1-2\nu)}{12} & -\frac{1-2\nu}{4} & -\frac{\sqrt{3}(1-2\nu)}{12} & \frac{1-2\nu}{4} & \frac{\sqrt{3}(1-2\nu)}{6} & \\ -\frac{\nu}{2} & \frac{\sqrt{3}(\nu-1)}{6} & \frac{\nu}{2} & \frac{\sqrt{3}(\nu-1)}{6} & 0 & \frac{\sqrt{3}(1-\nu)}{3} \end{bmatrix} \quad \text{(B.1)}$$

Element stiffness matrix $\mathbf{K}_{\text{IELSM, H9}}$ of H9:

$$\mathbf{K}_{\text{IELSM,H9}} = \sum_{i=1}^{9} \mathbf{K}_{n,i} + \mathbf{K}_{volume}$$

$$= \frac{E}{(1+\nu)(1-2\nu)} \begin{bmatrix}
\frac{\sqrt{3}(23-44\nu)}{48} & & & & & & & & & & & \\
\frac{4\nu-1}{16} & \frac{\sqrt{3}(7-12\nu)}{16} & & & & & & & & & & \\
\frac{\sqrt{3}(28\nu-15)}{48} & \frac{1-4\nu}{16} & \frac{\sqrt{3}(23-44\nu)}{48} & & & & & & & & sym & \\
\frac{4\nu-1}{16} & \frac{\sqrt{3}(4\nu-1)}{16} & \frac{1-4\nu}{16} & \frac{\sqrt{3}(7-12\nu)}{16} & & & & & & & & \\
\frac{\sqrt{3}(1-4\nu)}{24} & \frac{1-4\nu}{8} & \frac{\sqrt{3}(8\nu-3)}{24} & -\frac{1}{8} & \frac{\sqrt{3}(5-8\nu)}{12} & & & & & & & \\
0 & 0 & \frac{2\nu-1}{4} & \frac{\sqrt{3}(2\nu-1)}{4} & 0 & \frac{\sqrt{3}(1-2\nu)}{2} & & & & & & \\
\frac{\sqrt{3}(4\nu-3)}{48} & \frac{4\nu-3}{16} & \frac{\sqrt{3}(4\nu-1)}{48} & \frac{1-4\nu}{16} & \frac{\sqrt{3}(8\nu-3)}{24} & \frac{1-2\nu}{4} & \frac{\sqrt{3}(23-44\nu)}{48} & & & & & \\
\frac{4\nu-3}{16} & \frac{\sqrt{3}(4\nu-3)}{48} & \frac{4\nu-1}{16} & \frac{\sqrt{3}(1-4\nu)}{16} & \frac{1}{8} & \frac{\sqrt{3}(2\nu-1)}{4} & \frac{4\nu-1}{16} & \frac{\sqrt{3}(7-12\nu)}{16} & & & & \\
\frac{\sqrt{3}(4\nu-1)}{48} & \frac{4\nu-1}{16} & \frac{\sqrt{3}(4\nu-3)}{48} & \frac{3-4\nu}{16} & \frac{\sqrt{3}(1-4\nu)}{24} & 0 & \frac{\sqrt{3}(28\nu-15)}{48} & \frac{1-4\nu}{16} & \frac{\sqrt{3}(23-44\nu)}{48} & & & \\
\frac{1-4\nu}{16} & \frac{\sqrt{3}(1-4\nu)}{16} & \frac{3-4\nu}{16} & \frac{\sqrt{3}(4\nu-3)}{16} & \frac{4\nu-1}{8} & 0 & \frac{4\nu-1}{16} & \frac{\sqrt{3}(4\nu-1)}{16} & \frac{1-4\nu}{16} & \frac{\sqrt{3}(7-12\nu)}{16} & & \\
\frac{\sqrt{3}(8\nu-3)}{24} & \frac{1}{8} & \frac{\sqrt{3}(1-4\nu)}{24} & \frac{4\nu-1}{8} & \frac{\sqrt{3}(4\nu-3)}{16} & 0 & \frac{\sqrt{3}(1-4\nu)}{24} & \frac{1-4\nu}{8} & \frac{\sqrt{3}(8\nu-3)}{24} & -\frac{1}{8} & \frac{\sqrt{3}(5-8\nu)}{12} & \\
\frac{1-2\nu}{4} & \frac{\sqrt{3}(2\nu-1)}{4} & 0 & 0 & 0 & 0 & 0 & 0 & \frac{2\nu-1}{4} & \frac{\sqrt{3}(2\nu-1)}{4} & 0 & \frac{\sqrt{3}(1-2\nu)}{2}
\end{bmatrix} \quad \text{(B.2)}$$

Element stiffness matrix $\mathbf{K}_{\text{IELSM, H12}}$ of H12:

$$\mathbf{K}_{\text{IELSM,H12}} = \sum_{i=1}^{12} \mathbf{K}_{n,i} + \mathbf{K}_{volume}$$

$$= \frac{E}{(1+\nu)(1-2\nu)} \begin{bmatrix}
\frac{\sqrt{3}(23-44\nu)}{48} & & & & & & & & & & & \\
\frac{4\nu-1}{16} & \frac{\sqrt{3}(7-12\nu)}{16} & & & & & & & & & & \\
\frac{\sqrt{3}(20\nu-11)}{48} & \frac{1-4\nu}{16} & \frac{\sqrt{3}(23-44\nu)}{48} & & & & & & & & & \\
\frac{4\nu-1}{16} & \frac{\sqrt{3}(4\nu-1)}{16} & \frac{1-4\nu}{16} & \frac{\sqrt{3}(7-12\nu)}{16} & & & & & & sym & & \\
\frac{\sqrt{3}(10\nu-7)}{48} & -\frac{1+2\nu}{16} & \frac{\sqrt{3}(14\nu-5)}{48} & -\frac{1+2\nu}{16} & \frac{\sqrt{3}(5-8\nu)}{12} & & & & & & & \\
\frac{3(2\nu-1)}{16} & \frac{\sqrt{3}(2\nu-1)}{16} & \frac{3(2\nu-1)}{16} & \frac{3\sqrt{3}(2\nu-1)}{16} & 0 & \frac{\sqrt{3}(1-2\nu)}{2} & & & & & & \\
\frac{\sqrt{3}(1-4\nu)}{48} & \frac{1-4\nu}{16} & \frac{\sqrt{3}(4\nu-1)}{48} & \frac{1-4\nu}{16} & \frac{\sqrt{3}(14\nu-5)}{48} & \frac{3(1-2\nu)}{16} & \frac{\sqrt{3}(23-44\nu)}{48} & & & & & \\
\frac{1-4\nu}{16} & \frac{\sqrt{3}(1-4\nu)}{16} & \frac{4\nu-1}{16} & \frac{\sqrt{3}(4\nu-3)}{16} & \frac{1+2\nu}{16} & \frac{3\sqrt{3}(2\nu-1)}{16} & \frac{4\nu-1}{16} & \frac{\sqrt{3}(7-12\nu)}{16} & & & & \\
\frac{\sqrt{3}(4\nu-1)}{48} & \frac{4\nu-1}{16} & \frac{\sqrt{3}(1-4\nu)}{48} & \frac{4\nu-1}{16} & \frac{\sqrt{3}(10\nu-7)}{48} & \frac{3(1-2\nu)}{16} & \frac{\sqrt{3}(20\nu-11)}{48} & \frac{1-4\nu}{16} & \frac{\sqrt{3}(23-44\nu)}{48} & & & \\
\frac{1-4\nu}{16} & \frac{\sqrt{3}(4\nu-3)}{16} & \frac{4\nu-1}{16} & \frac{\sqrt{3}(1-4\nu)}{16} & \frac{1+2\nu}{16} & \frac{\sqrt{3}(2\nu-1)}{16} & \frac{4\nu-1}{16} & \frac{\sqrt{3}(4\nu-1)}{16} & \frac{1-4\nu}{16} & \frac{\sqrt{3}(7-12\nu)}{16} & & \\
\frac{\sqrt{3}(14\nu-5)}{48} & \frac{1+2\nu}{16} & \frac{\sqrt{3}(10\nu-7)}{48} & \frac{1+2\nu}{16} & \frac{\sqrt{3}(1-4\nu)}{12} & 0 & \frac{\sqrt{3}(10\nu-7)}{48} & -\frac{1+2\nu}{16} & \frac{\sqrt{3}(14\nu-5)}{48} & -\frac{1+2\nu}{16} & \frac{\sqrt{3}(5-8\nu)}{12} & \\
-\frac{3(2\nu-1)}{16} & \frac{3\sqrt{3}(2\nu-1)}{16} & \frac{3(1-2\nu)}{16} & \frac{\sqrt{3}(2\nu-1)}{16} & 0 & 0 & \frac{3(2\nu-1)}{16} & \frac{\sqrt{3}(2\nu-1)}{16} & \frac{3(2\nu-1)}{16} & \frac{3\sqrt{3}(2\nu-1)}{16} & 0 & \frac{\sqrt{3}(1-2\nu)}{2}
\end{bmatrix} \quad \text{(B.3)}$$

Element stiffness matrix $\mathbf{K}_{\text{IELSM, H15}}$ of H15:

$$\mathbf{K}_{\text{IELSM,H15}} = \sum_{i=1}^{15} \mathbf{K}_{n,i} + \mathbf{K}_{volume}$$

$$= \frac{E}{(1+\nu)(1-2\nu)} \begin{bmatrix}
\frac{\sqrt{3}(17-32\nu)}{48} & & & & & & & & & & & \\
\frac{1}{16} & \frac{\sqrt{3}(19-32\nu)}{48} & & & & & & & & & & \\
\frac{\sqrt{3}(12\nu-7)}{48} & \frac{1-4\nu}{16} & \frac{\sqrt{3}(17-32\nu)}{48} & & & & & & & & & \\
\frac{4\nu-1}{16} & \frac{\sqrt{3}(4\nu-1)}{16} & -\frac{1}{16} & \frac{\sqrt{3}(19-32\nu)}{48} & & & & & sym & & & \\
\frac{\sqrt{3}(\nu-1)}{12} & -\frac{\nu}{4} & \frac{\sqrt{3}(3\nu-1)}{12} & -\frac{\nu}{4} & \frac{\sqrt{3}(5-8\nu)}{12} & & & & & & & \\
\frac{2\nu-1}{8} & \frac{\sqrt{3}(2\nu-1)}{24} & \frac{2\nu-1}{8} & \frac{\sqrt{3}(2\nu-1)}{8} & 0 & \frac{\sqrt{3}(1-2\nu)}{3} & & & & & & \\
-\frac{\sqrt{3}}{48} & -\frac{1}{16} & \frac{\sqrt{3}(4\nu-1)}{48} & \frac{1-4\nu}{16} & \frac{\sqrt{3}(3\nu-1)}{12} & \frac{1-2\nu}{8} & \frac{\sqrt{3}(17-32\nu)}{48} & & & & & \\
-\frac{1}{16} & -\frac{\sqrt{3}}{16} & \frac{4\nu-1}{16} & \frac{\sqrt{3}(4\nu-5)}{48} & \frac{\nu}{4} & \frac{\sqrt{3}(2\nu-1)}{8} & \frac{1}{16} & \frac{\sqrt{3}(19-32\nu)}{48} & & & & \\
\frac{\sqrt{3}(4\nu-1)}{48} & \frac{4\nu-1}{16} & -\frac{\sqrt{3}}{48} & \frac{1}{16} & \frac{\sqrt{3}(\nu-1)}{12} & \frac{1-2\nu}{8} & \frac{\sqrt{3}(12\nu-7)}{48} & \frac{1-4\nu}{16} & \frac{\sqrt{3}(17-32\nu)}{48} & & & \\
\frac{1-4\nu}{16} & \frac{\sqrt{3}(4\nu-5)}{48} & \frac{1}{16} & -\frac{\sqrt{3}}{16} & \frac{\nu}{4} & \frac{\sqrt{3}(2\nu-1)}{24} & \frac{4\nu-1}{16} & \frac{\sqrt{3}(4\nu-1)}{16} & -\frac{1}{16} & \frac{\sqrt{3}(19-32\nu)}{48} & & \\
\frac{\sqrt{3}(3\nu-1)}{12} & \frac{\nu}{4} & \frac{\sqrt{3}(\nu-1)}{12} & \frac{\nu}{4} & -\frac{\sqrt{3}}{12} & 0 & \frac{\sqrt{3}(\nu-1)}{12} & -\frac{\nu}{4} & \frac{\sqrt{3}(3\nu-1)}{12} & -\frac{\nu}{4} & \frac{\sqrt{3}(5-8\nu)}{12} & \\
\frac{1-2\nu}{8} & \frac{\sqrt{3}(2\nu-1)}{8} & \frac{1-2\nu}{8} & \frac{\sqrt{3}(2\nu-1)}{24} & 0 & 0 & -\frac{1-2\nu}{8} & \frac{\sqrt{3}(2\nu-1)}{24} & -\frac{1-2\nu}{8} & \frac{\sqrt{3}(2\nu-1)}{8} & 0 & \frac{\sqrt{3}(1-2\nu)}{3}
\end{bmatrix} \tag{B.4}$$